\documentclass[a4paper,12pt, abstract=true]{scrartcl}

\usepackage[ansinew]{inputenc}
\usepackage[english]{babel}
\usepackage{latexsym,amssymb,amsfonts,amsmath,bbm}
\usepackage{enumitem} 
\usepackage{theorem}
\usepackage{xcolor}
\usepackage{graphicx}
\usepackage{esint}
\usepackage{stmaryrd}
\usepackage[normalem]{ulem} 

\numberwithin{equation}{section}
\theoremstyle{plain}

\newtheorem{dummytheorem}{Dummy-Theorem}[section]
\newcommand{\proofendsign}{$\Box$} 

\newtheorem{lemma}[dummytheorem]{Lemma}
\newtheorem{theorem}[dummytheorem]{Theorem}
\newtheorem{proposition}[dummytheorem]{Proposition}

\newtheorem{corollary}[dummytheorem]{Corollary}

\newtheorem{example}[dummytheorem]{Example}

\newtheorem{remark}[dummytheorem]{Remark}

\newenvironment{proof}{{\noindent \textbf{Proof} }}
 {{\hspace*{\fill}\proofendsign\par\bigskip}}

\newcommand{\Flinks}{F^{\leftarrow}}

\newcommand{\MMM}{\mathbb{M}}
\newcommand{\NNN}{\mathbb{N}}

\newcommand{\RRR}{\mathbb{R}}

\newcommand{\FFF}{\mathbb{F}}
\newcommand{\GGG}{\mathbb{G}}
\newcommand{\GGGH}{\GGG^{G,H}}

\newcommand{\EEE}{\mathbb{E}}

\newcommand{\pr}{\mathbb{P}}
\newcommand{\MP}{\mathbb{P}}
\newcommand{\MQ}{\mathbb{Q}}

\newcommand{\vari}{\mathbb{V}\textrm{ar}}

\newcommand{\eins}{\mathbbm{1}}

\newcommand{\cB}{\mathcal B}
\newcommand{\cC}{\mathcal C}
\newcommand{\cF}{\mathcal F}

\newcommand{\cK}{\mathcal K}

\newcommand{\cM}{\mathcal M}

\newcommand{\cU}{\mathcal U}
\newcommand{\cV}{\mathcal V}

\newcommand{\OFP}{(\Omega,{\cal F},\pr)}

\def\bcswitch{\left\{\renewcommand{\arraystretch}{0.5}\begin{array}{c@{,~}c}}

\def\ecswitch{\end{array}\right.}

\begin{document}
\title{Estimating the errors for solutions of the SAA method to solve compound and risk averse stochastic programs}
\author{
Volker Krätschmer\footnote{Faculty of Mathematics, University of Duisburg--Essen, { volker.kraetschmer@uni-due.de}}
}
\date{~}

\maketitle
\begin{abstract}
This paper is a study on solutions of the Sample Average Approximation Method to solve compound stochastic programs. We derive nonasymptotic upper estimates for probabilities of the approximation errors. The results depend on the sample size with explicit terms instead of unspecified universal constants. They allow to conclude immediately nonasymptotic rates for the optimal solutions, and they may be utilized to construct nonasymptotic confidence regions for unique solutions of the genuine compound stochastic programs. In the special case of classical risk neutral stochastic programs, we end up with upper estimates of deviation probabilities for $M$-estimators, and their nonasymptotic rates. Moreover, we may also demonstrate how to apply the results to sample average approximation of risk averse stochastic programs. In this respect we consider stochastic programs expressed in terms of absolute semideviation risk measures and Average Value at Risk. The investigations are based on concentration inequalities from the recent contribution \cite{Kraetschmer2023a}. The line of reasoning does not rely on pathwise analytical properties of the objectives. In particular, continuity or convexity in the parameter is not imposed in advance as usual in the literature on the Sample Average Approximation method. It is also shown that objectives with H\"older continuous paths meet the requirements of the main results. Moreover, the main results are applied to objectives whose paths are piecewise H\"older continuous, as e.g. in two stage mixed-integer programs.
\end{abstract}

\section{Introduction}
Consider a compound stochastic program 
\begin{equation}
\label{compound optimization}
\min_{\theta\in\Theta}\EEE\big[H\big(\theta,\EEE[G(\theta,Z_{1})],Z_{1}\big)\big],
\end{equation}
where $\Theta$ denotes a compact subset of $\RRR^{m}$, whereas $Z$ stands for a $d$-dimensional random vector with distribution $\MP^{Z}$. Compound stochastic programs extend the class of classical risk neutral stochastic programs, and they also enclose some classes of risk averse stochastic pograms. Outstanding examples are stochastic programs with \textit{absolute semideviation risk measures} and \textit{Average Value at Risk} as objectives.
\par
In general the parameterized distribution of the goal function $G$ is unknown, but some information is available by i.i.d. samples. Using this information, a general device to solve approximately  problem (\ref{compound optimization}) is provided by the so-called \textit{Sample Average Approximation} (\textit{SAA}). For explanation, let us consider a sequence $(Z_{j})_{j\in\NNN}$ of independent $d$-dimensional random vectors on some fixed probability space $\OFP$ which are identically distributed as the $d$-dimensional random vector $Z$.
Then the SAA method approximates the genuine optimization problem (\ref{compound optimization}) by the following one
\begin{equation}
\label{compound SAA}
\min_{\theta\in\Theta}\frac{1}{n}~\sum_{j=1}^{n}H\Big(\theta,\frac{1}{n}~\sum_{k=1}^{n}G(\theta,Z_{k}),Z_{j}\Big)\quad(n\in\NNN).
\end{equation}
In the classical risk neutral case, i.e. $H(\theta,\cdot,z)\equiv G(\theta,z)$ for 
$(\theta,z)\in\Theta\times\RRR^{d}$, an important subject is to analyze asymptotic distributions of the 
stochastic sequence
$$
\left(\sqrt{n}\Big[\inf_{\theta\in\Theta}~\frac{1}{n}~\sum_{k=1}^{n}G(\theta,Z_{k}) - \inf_{\theta\in\Theta}\EEE[G(\theta,Z_{1})]\right)_{n\in\NNN}.
$$
Standard results may be found in the monograph \cite{ShapiroEtAl}. For more general compound programs only a few investigations are known, e.g. by \cite{DentchevaEtAl2017}, \cite{ErmolievNorkin2013} and \cite{GuiguesKraetschmerShapiro2018}. Most of the known contributions are based on analytical path properties of the process $G$ like continuity or convexity. Very recently, in \cite{Kraetschmer2020} the investigations have been extended in respect of allowing for more general goal functions $G$, which make possible to apply the results to stochastic programs whose goal functions are neither continuous nor convex in the parameter.  
\smallskip

In recent years some attention has been paid to the nonasymptotic analysis of the SAA method. The contributions \cite{GuiguesEtAl2017} and \cite{Kraetschmer2023a} deal with upper bounds for the deviation probabilities
\begin{align*}
\MP\Big\{\sqrt{n}\Big|\inf_{\theta\in\Theta}~\frac{1}{n}~\sum_{k=1}^{n}G(\theta,Z_{k}) - \inf_{\theta\in\Theta}\EEE[G(\theta,Z_{1})]\Big| > \varepsilon\Big\}\quad(\varepsilon > 0)
\end{align*}
in terms of the sample size $n$. The study \cite{GuiguesEtAl2017} focusses on processes $\big(G(\theta,Z))_{\theta\in\Theta}$ with convex paths and light-tailed one-dimensional marginal distributions. The authors achieve to find bounds which do not depend on the dimension of the parameter set. In some sense the investigations in \cite{Kraetschmer2023a} are more general, e.g. avoiding to require analytical properties of the paths in advance. Also more general marginal distributions are allowed, going beyond light tail behaviour. In addition the results are extended to compound and risk averse stochastic programs. However, the resulting bounds may increase with the dimension of the paramter set even in the classical risk neutral case.
\smallskip

Another aspect of the SAA is to study the behaviour of their optimal solutions 
$\widehat{\theta}_{n}$. They are known as $M$-estimators in the risk neutral case. 
Based on the assumption that the stochastic program \eqref{compound optimization} has a unique solution $\theta^{*}$ we are interested in the deviation probabilities 
\begin{equation*}
\label{deviation Einleitung}
\MP\big(\big\{\|\widehat{\theta}_{n} - \theta^{*}\|_{m} > \varepsilon\cdot n^{-\gamma}\big\}\big)\quad(n\in\NNN, \varepsilon > 0)
\end{equation*}
for suitable $\gamma\in ]0,\infty[$ dependent on the sample size, where $\|\cdot\|_{m}$ denotes the Euclidean norm on $\RRR^{m}$. Upper bounds of these probabilities might be used to construct nonasymptotic confidence regions for the unique solution $\theta^{*}$.
\smallskip

To the best of our knowledge very little is known about the deviation probabilities, even in the case of $M$-estimation. The literature on $M$-estimation is mainly focussed on criteria to derive convergence rates for the solutions (see e.g. \cite{vanderVaartWellner1996}, \cite{Kosorok2008}) instead of nonasymptotic rates. An exception is provided by the recent contribution \cite{OstrovskiiBach2021} on nonasymptotic estimation of the deviation probabilities. The investigations there rely on objectives $G$ which are of a very special form, in particular they are three times differentiable in the parameter. Let us also refer to the paper \cite{Guigues2017}. There in the risk neutral case the author proposes an alternative way to construct nonasymptotic confidence regions for the unique solution of \eqref{compound optimization}. It is based on a multistep version of the stochastic mirror descent algorithm instead of $M$-estimation. The bounds on the probabilities for the approximation errors are formulated in terms of the number of iterations.
\smallskip

An early comprehensive study of general compound sample average approximation has been presented in \cite{ErmolievNorkin2013}. One of the main results in this paper deals with the deviation probabilities for Hausdorff distances between approximate solution sets of the genuine problem \eqref{compound optimization} and the SAA problems \eqref{compound SAA}. Unfortunately, it may not be utilized directly for the desired upper estimation of the deviation probabilities.
\smallskip

Our contribution is to derive estimates of the deviation probabilities for general objectives without imposing analytical properties of the paths in advance.  Inspired by \cite{Pflug1999} an essential ingredient will be to require some condition on small increments 
$\GGG^{G,H}_{n}(\theta,\cdot) - \GGG^{G,H}_{n}(\theta^{*},\cdot)$ of the \textit{compound empirical process}
\begin{align*}
\mathbb{G}_{n}^{G,H}(\theta,\cdot) 
= 
\frac{1}{\sqrt{n}}~\sum_{j=1}^{n}\left(H\Big(\theta,\frac{1}{n}~\sum_{k=1}^{n}G\big(\theta,Z_{k}\big),Z_{j}\Big) - \EEE\big[H\big(\theta,\EEE[G(\theta,Z_{1})],Z_{1}\big)\big]\right),
\end{align*} 
for almost all sample sizes $n$.  
The key will be to look at upper estimates for the probabilities
\begin{align*}
\MP\Big(\Big\{\sup_{\|\theta - \theta^{*}\|_{m}\leq\delta}\big|\GGGH_{n}(\theta,\cdot) - \GGGH_{n}(\theta^{*},\cdot)\big| > \varepsilon~n^{-\gamma}\Big\}\Big)\quad(\varepsilon > 0),
\end{align*}
dependent on the sample size. 
Our investigations do not rely on analytical properties for the paths of $G$. 
Instead we restrict ourselves to processes $G$, where the 
associated families $\FFF^{\Theta} := \{G(\theta,\cdot)\mid \theta\in\Theta\}$ and $
\FFF^{\Theta}_{\delta} := \{G(\theta,\cdot) - G(\theta^{*},\cdot)\mid \theta\in\Theta, \|\theta - \theta^{*}\|_{m}\leq\delta\}$ ($\delta > 0$ small) of Borel measurable mappings are ``small'' in the sense that they have finite 
uniform entropy integrals. 
\medskip

The paper is organized as follows. In the following section we shall point out the starting conditions on the mappings $G$ and $H$. Furthermore we shall provide the general notations used throughout the paper. Afterwards, we shall develop in Section \ref{deviation probabilities compound SAA} nonasymptotic upper bounds for the deviation probabilities \eqref{deviation Einleitung}. Firstly, we shall clarify the relationship between these probabilities and the concentration inequalities for the increments $\GGG^{G,H}_{n}(\theta,\cdot) - \GGG^{G,H}_{n}(\theta^{*},\cdot)$. This will also motivate our specific additional requirements for the function classes $\FFF^{\Theta}$ and $\FFF^{\Theta}_{\delta}$. Then we shall present our two main results. 
The first one works for all sample sizes. However, it does not take into account tail behaviour of the parameterized random variables $G(\theta,\cdot)$ from the process $G$. If these variables have simultenously sub-Weibull tails (e.g. light tails), then the second main result improves the former one by exponential bounds. All the results provide explicit bounds instead of using unspecified universal constants. As an immediate by product the estimates of the deviation probabilities give nonasymptotic rates for the solutions of the SAA. We shall also illustrate how to construct  nonasymptotic confidence regions for the unique solution $\theta^{*}$ from our main results. This will be done in Section \ref{confidence regions}. 
\par
In the following Section \ref{specific classes} the assumptions of the main results are exemplified by some specific classes of objectives $G$. We shall consider the case that $G$ has H\"older continuous paths. Also processes $G$ are discussed whose paths are piecewise H\"older continuous but not necessarily continuous. Value functions of two stage mixed-integer programs are typical examples for processes of such a kind. 
\par
Our line of reasoning relies on a second order growth condition for the goal function of 
the genuine problem \eqref{compound optimization} and a representation of the function $H$. Both will be discussed in Section \ref{illustration assumptions}.
\par
Simplifications of our main results in the case of $M$-estimation will be pointed out in Section \ref{error estimates m-estimators}. In particular an alternative approach is provided to find convergence rates for $M$-estimators. 
Also the application to risk averse stochastic programs in terms of absolute semideviations will be discussed in more detail. This will be done in Section \ref{error estimates absolut semideviations}. The subsequent Section \ref{estimates under Average Value at Risk} is devoted to the SAA of risk averse stochastic programs under Average Value at Risk. In this situation, we may not apply directly our results on deviation probabilities \eqref{deviation Einleitung}. But it will be demonstrated how to obtain upper estimates for the deviation probabilities from them.
\par 
The mathematical key of our argumentation is the subject of Section \ref{MainResult}. Here we shall present our results on concentration inequalities for increments of the compound empirical process. They are interesting in their own right. Finally, 
Section \ref{Beweise} gathers proofs of several results from the different sections. The main tools are concentration inequalities from the recent paper \cite{Kraetschmer2023a} which will be recalled at the beginning of the section.
%
\section{The set up and notations}
\label{set up}
Let the probability space $(\Omega,\cF, \pr)$ from the introduction be atomless and complete, and let us recall that $\Theta$ is a compact subset of $\RRR^{m}$. We shall denote the diameter of $\Theta$ w.r.t. the Euclidean metric by $\Delta(\Theta)$, assuming $\Delta(\Theta)> 0$. We shall use the symbols $\|\cdot\|_{m}$ and $\|\cdot\|_{m+1}$ for the Euclidean norms on $\RRR^{m}$ and $\RRR^{m+1}$ respectively. The transpose of a vector $y$ from a Euclidean space will be denoted by $y^{T}$. As a further general notation for any $\varepsilon > 0$ we shall write $K_{\varepsilon}$ 
to denote the largest integer that does not exceed $\ln(\varepsilon)/\ln(2)$.
\par
Denoting by $Z$ any $d$-dimensional random vector with distribution $\MP^{Z}$, we shall investigate the following compound form of stochastic program
\begin{equation}
\label{composite optimization}
\min_{\theta\in\Theta}\EEE\big[H\big(\theta,\EEE[G(\theta,Z)],Z\big)\big] =: \min_{\theta\in\Theta} \psi_{H,\Theta}(\theta),
\end{equation}
where the mappings $G: \Theta\times\RRR^{d}\rightarrow\RRR$ and $H:\Theta\times\RRR\times\RRR^{d}\rightarrow\RRR$ satisfy the following properties 
\begin{enumerate}[label=(A \arabic*), ref=(A \arabic*)]
\item \label{assump:A1} $G(\theta,\cdot)$ is $\MP^{Z}$-integrable for every $\theta\in\Theta$, and $G$ is assumed to satisfy the following continuity property
\begin{center}
$G(\theta_{k},\cdot)\to G(\theta,\cdot)$ in $\MP^{Z}$-probability whenever $\theta_{k}\to\theta$ w.r.t. the Euclidean metric.
\end{center}
\item \label{assump:A2} 
There is some strictly positive square $\MP^{Z}$-integrable mapping $\xi^{G}:\RRR^{d}\rightarrow\RRR$ such that 
$$
\sup\limits_{\theta\in\Theta}\big|G(\theta,z)\big|\leq\xi^{G}(z)\quad\mbox{for}~z\in\RRR^{d}.
$$ 
\item \label{measurability H} 
$H(\theta,t,\cdot)$ is Borel measurable for $(\theta,t)\in\Theta\times\RRR$, and $H\big(\theta,\EEE[G\big(\theta,Z_{1})],\cdot\big)$ is $\MP^{Z}$-integrable for $\theta\in\Theta$.
\item \label{assump:A4}
There exist nonnegative square $\MP^{Z}$-integrable mappings $L_{1},L_{2}$ such that
\begin{align*}
&
\big|H\big(\theta,t,z\big) - H\big(\vartheta,s,z\big)\big|^{2}\\
&
\leq L_{1}(z)^{2}~|G(\theta,z) - G(\vartheta,z)|^{2} + L_{2}(z)^{2}~|t - s|^{2}
\quad\mbox{for}~\theta,\vartheta\in\Theta, t,s\in\RRR, z\in\RRR^{d}.
\end{align*}
\end{enumerate}
\begin{remark}
\label{Stetigkeiten}
We may invoke Vitalis theorem (see e.g. \cite[Proposition 21.4]{Bauer2001}) to see that the mapping $\theta\mapsto\EEE[G(\theta,Z_{1})]$ on $\Theta$ is continuous under \ref{assump:A1} and \ref{assump:A2}. Moreover, by assumption \ref{assump:A4} we may also observe that the sequence  
$$
\Big(H\big(\theta_{k},\EEE[G(\theta_{k},Z_{1})],Z_{1}\big) -  
H\big(\theta,\EEE[G(\theta,Z_{1})],Z_{1}\big)\Big)_{n\in\NNN}
$$ 
converges in probability to $0$ for $\|\theta_{k}-\theta\|_{m}\to 0$, and together with \ref{assump:A2} it is dominated by some integrable random variable.  Hence by Vitalis theorem again we may show that the mapping 
$\theta\mapsto \EEE\big[H\big(\theta,\EEE[G(\theta,Z_{1})],Z_{1}\big)\big]$ on $\Theta$ is a continuous mapping w.r.t. the Euclidean norm. In particular, optimization problem 
\eqref{composite optimization} has a solution due to compactness of $\Theta$.
\end{remark}
Henceforth we set for abbreviation
\begin{align*}
\overline{\psi}(\theta) := \EEE[G(\theta,Z_{1})]\quad\mbox{and}\quad\psi_{H,\Theta}(\theta) := \EEE\big[H\big(\theta,\overline{\psi}(\theta),Z_{1}\big)\big]\quad(\theta\in\Theta).
\end{align*}
A general device to solve approximately problem \eqref{composite optimization} is to use the \textit{sample average approximation} (SAA). For explanation let $(Z_{j})_{j\in\NNN}$ be a sequence of independent $d$-dimensional random vectors on the complete atomless probability $\OFP$ which are identically distributed as the random vector $Z$. Then the SAA optimization problem based on the i.i.d. sample $(Z_{1},\ldots,Z_{n})$ associated with \eqref{composite optimization} reads as follows
\begin{equation}
\label{composite SAA}
\min_{\theta\in\Theta}\frac{1}{n}~\sum_{j=1}^{n}H\Big(\theta,\frac{1}{n}~\sum_{k=1}^{n}G(\theta,Z_{k}),Z_{j}\Big)\quad(n\in\NNN).
\end{equation}
Note that in view of assumptions \ref{assump:A1} - \ref{assump:A4} the optimization problems \eqref{composite SAA} always have finite infimal values.
\smallskip

We always assume that the genuine optimization problem \eqref{composite optimization} has a unique solution $\theta^{*}$. 
\begin{align}
\label{unique solution}
\mbox{Optimization problem}~\eqref{composite optimization}~\mbox{has a unique solution}~\theta^{*}\in\Theta. 
\end{align}
Furthermore we fix any sequence $(\widehat{\theta}_{n})_{n\in\NNN}$ of random vectors $\widehat{\theta}_{n}$ which minimize the SAA problem w.r.t. the sample size $n$. Our aim is to investigate the 
probabilities
\begin{align}
\label{deviation}
\MP\big(\big\{\|\widehat{\theta}_{n} - \theta^{*}\|_{m} > \varepsilon\cdot n^{-\gamma}\big\}\big)\quad(n\in\NNN, \varepsilon > 0)
\end{align}
for suitable $\gamma\in ]0,\infty[$. Such bounds will enable us to construct nonasymptotic confidence regions for the solution $\theta^{*}$.
\smallskip

The essential part in our line of reasoning is to find concentration inequalities for small increments of the \textit{compound empirical process} $\mathbb{G}_{n}^{G,H}:\Theta\times\Omega\rightarrow\RRR$, defined by
\begin{align*}
\mathbb{G}_{n}^{G,H}(\theta,\omega) 
= 
\frac{1}{\sqrt{n}}~\sum_{j=1}^{n}\left(H\Big(\theta,\frac{1}{n}~\sum_{k=1}^{n}G\big(\theta,Z_{k}(\omega)\big),Z_{j}(\omega)\Big) - \EEE\big[H\big(\theta,\EEE[G(\theta,Z_{1})],Z_{1}\big)\big]\right).
\end{align*}
This means to provide upper estimates of the following probabilities
\begin{align}
\label{deviations compound}
\MP\Big(\Big\{\sup_{\theta\in\cU_{\delta}}\big|\GGGH_{n}(\theta,\cdot) - \GGGH_{n}(\theta^{*},\cdot)\big| > \varepsilon~n^{-\gamma}\Big\}\Big)\quad(\varepsilon > 0)
\end{align} 
for small $\delta > 0$ and suitable $\gamma > 0$, where $\cU_{\delta}:= 
\{\theta\in\Theta\mid \|\theta - \theta^{*}\|_{m}\}$. For this purpose we shall use methods from empirical process theory. As usual in this theory, we define the \textit{empirical process operator}
\begin{equation*}
(\MP_{n}-\MP)(f) := \frac{1}{n}\sum_{j=1}^{n}\Big(f(Z_{j}) - \EEE[f(Z_{j})]\Big)
\end{equation*}
for $\MP^{Z}$-integrable mappings $f:\RRR^{d}\rightarrow\RRR$. We shall use this notation throughout the paper. The devices from empirical process theory will mainly be applied to the following function classes 
\begin{align*}
\FFF^{\Theta} := \{G(\theta,\cdot)\mid\theta\in\Theta\}\quad\mbox{and}\quad
 \FFF_{\delta}^{\Theta} := \{G(\theta,\cdot) - G(\theta^{*},\cdot)\mid\theta\in\cU_{\delta}\}~(\delta > 0). 
\end{align*}
\smallskip

\section{Deviation probabilities for solutions of compound SAA}
\label{deviation probabilities compound SAA}
Throughout this section we want to derive upper estimations for the probabilities \eqref{deviation}. In \cite{Pflug1999} the author illustrates how to achieve such estimations via variational inequalities. We shall adapt this idea to our situation. As a starting point the unique solution $\theta^{*}$ 
will be required to fulfill the \textit{second order growth condition}.
\begin{enumerate}[label=(A \arabic*), ref=(A \arabic*)]
\setcounter{enumi}{4}
\item \label{assump:A7} There exists $M_{1} > 0$ such that the goal function $\psi_{H,\Theta}$ of optimization \eqref{composite optimization} satisfies
\begin{align*}
\psi_{H,\Theta}(\theta) - \psi_{H,\Theta}(\theta^{*})\geq M_{1}~\|\theta - \theta^{*}\|_{m}^{2}\quad\mbox{for}~\theta\in\Theta.
\end{align*}
\end{enumerate}
Often local versions of the second order condition are used. Actually, in our setting local versions are equivalent with the global version \ref{assump:A7}. This will be discussed in Section \ref{illustration assumptions}.
\medskip

Now, the adaption of the idea from \cite{Pflug1999} leads to the following result in terms of the deviation probabilities \eqref{deviations compound} for the compound empirical processes. It will be the basic step of our investigations. It is in parts known from textbook proofs of results on convergences rates for $M$-estimators (see \cite{Kosorok2008}, \cite{vanderVaartWellner1996}). 
\begin{proposition}
\label{variational inequalities}
Let \eqref{unique solution} and the assumptions \ref{assump:A1} - \ref{assump:A7} be fulfilled. Then with $M_{1} > 0$ from \ref{assump:A7}, 
and with the mapping $L_{2}$ from \ref{assump:A4}
\begin{align*}
& 
\MP\big(\big\{\|\widehat{\theta}_{n} - \theta^{*}\|_{m} > \varepsilon\cdot n^{-\gamma}\big\}\cap\overline{\Omega}\big)\\
&\leq 
\sum_{k= K_{\varepsilon} + 1\atop 2^{k-1}/n^{\gamma}\leq\Delta(\Theta)}^{\infty}\hspace*{-0.25cm}\MP^{*}\Big(\Big\{\sup_{\theta\in U_{\Delta(\Theta)\wedge (2^{k}/n^{\gamma})}}\big|\GGG_{n}^{G,H}(\theta,\cdot) - \GGG_{n}^{G,H}(\theta^{*},\cdot) \big| > 
M_{1}~2^{2 (k - 1)}~n^{- 2\gamma + 1/2} 
\Big\}\cap\overline{\Omega}\Big)
\end{align*}
holds for $\overline{\Omega}\in\mathcal{F}$, $n\in\NNN$, $\delta\in ]0,\Delta(\Theta)]$ and $\varepsilon, \gamma\in ]0,\infty[$. Here $\MP^{*}$ denotes the outer probability w.r.t. $\MP$. 
\end{proposition}
Proposition \ref{variational inequalities} will be shown in Subsection \ref{Proof of Proposition variational inequalities}.
\medskip

According to Proposition \ref{variational inequalities} we may upper estimate the deviation probabilities \eqref{deviation} by suitable concentration inequalities for locally small increments of the compound empirical processes $\GGG_{n}^{G,H}$. In deriving such inequalities we want to avoid subtleties of measurability by  
imposing the following separability condition on the objective $G$.
\begin{enumerate}[label=(A \arabic*), ref=(A \arabic*)]
\setcounter{enumi}{5}
\item \label{assump:A3} For $\delta > 0$ there exist some at most countable subset $\cC(\cU_{\delta})$ of the set $\cU_{\delta}$ and $(\MP^{Z})^{n}$-null sets $N_{n}^{\delta}$ $(n\in\NNN)$ such that 
$$
\inf_{\vartheta\in\cC(\cU_{\delta})}\left\{\EEE[|G(\theta,Z_{1}) - G(\vartheta,Z_{1})|] + \max_{j\in\{1,\ldots,n\}}\big|G(\theta,z_{j}) - G(\vartheta,z_{j})\big|\right\} = 0\quad
$$
if $\theta\in\cU_{\delta}$, $n\in\NNN$ and $(z_{1},\ldots,z_{n})\in\RRR^{d n}\setminus N_{n}^{\delta}$.
\end{enumerate}
By compactness, $\Theta = \cU_{\delta}$ holds for every $\delta$ larger than the diameter of $\Theta$. In particular the separability property in \ref{assump:A3} is also satisfied for the entire parameter set $\Theta$. Hence assumption \ref{assump:A3} will allow us to invoke results from \cite{Kraetschmer2023a}, gathered in Theorem \ref{upper estimation} below, to develop upper estimations for the probabilities on the right hand side of the inequality in the statement of Proposition \ref{variational inequalities}. The essential ingredient of these results is to replace analytical conditions on the paths with requirements which intuitively make the function classes $\FFF^{\Theta}$ and $\FFF^{\Theta}_{\delta}$ small in a certain sense. The precise formulation of this intuition is based on covering numbers for classes of Borel measurable mappings from $\RRR^{d}$ into $\RRR$ w.r.t. $L^{p}$-norms.  To recall these concepts adapted to our situation, let us fix any nonvoid set $\FFF$ of Borel measurable mappings from $\RRR^{d}$ into $\RRR$ and any probability measure $\MQ$ on $\cB(\RRR^{d})$ with metric $d_{\MQ,p}$ induced by the $L^{p}$-norm $\|\cdot\|_{\MQ,p}$ for $p\in [1,\infty[$.  
\begin{itemize}
\item \textit{Covering numbers for $\FFF$}\\
We use $N\big(\eta,\FFF,L^{p}(\MQ)\big)$ to denote the minimal number to cover $\FFF$ by closed $d_{\MQ,p}$-balls of radius $\eta > 0$ with centers in $\FFF$. We define $N\big(\eta,\FFF,L^{p}(\MQ)\big) := \infty$ if no finite cover is available. 
\item An \textit{envelope} of $\FFF$ is defined to mean some Borel measurable mapping $C_{\FFF}:\RRR^{d}\rightarrow\RRR$ satisfying $\sup_{h\in\FFF}|h|\leq C_{\FFF}$. If an envelope $C_{\FFF}$ has strictly positive outcomes, we shall speak of a \textit{positive envelope}.
\item $\cM_{\textrm{\tiny fin}}$ denotes the set of all probability measures on $\cB(\RRR^{d})$ with finite support.
\end{itemize}
A function class will be considered as ``small'' if it has finite uniform entropy integrals.
With a slight abuse of convention, 
we shall call 
\begin{align}
&
\label{Entropie-Integral I}
J(\FFF,C_{\FFF},\delta) := \int_{0}^{\delta}\sup_{\MQ\in \cM_{\textrm{\tiny fin}}}\sqrt{\ln\big(2~N\big(\varepsilon~\|C_{\FFF}\|_{\MQ,2},\FFF,L^{2}(\MQ)\big)\big)}~d\varepsilon
\end{align}
the \textit{uniform entropy integral (up to $\delta$)} of the function class $\FFF$ with positive envelope $C_{\FFF}$. 
\medskip

As already mentioned the function classes $\FFF^{\Theta}$ and $\FFF^{\Theta}_{\delta}$ are the relevant ones to develop via the results from \cite{Kraetschmer2023a} the desired upper estimations for the deviation probabilities \eqref{deviation}. 
Note that \ref{assump:A2} means nothing else but the existence of some square $\MP^{Z}$-integrable positive envelope of $\FFF^{\Theta}$. 
 We assume that the function classes $\FFF^{\Theta}$ and $\FFF^{\Theta}_{\delta}$ are ``small'' in the sense that they have finite uniform entropy integrals:
\begin{enumerate}[label=(A \arabic*), ref=(A \arabic*)]
\setcounter{enumi}{6}
\item \label{assump:A5} Under \ref{assump:A2} with function $\xi^{G}$ and \ref{assump:A4} with function $L_{1}$, there exist $\beta\in ]0,1]$, $\overline{M}_{1}, \overline{M}^{1} > 0$, and a family $(\xi^{G}_{\delta})_{\delta\in ]0,\Delta(\Theta)]}$ of square $\MP^{Z}$-integrable positive envelopes $\xi^{\Theta}_{\delta}$ of $\FFF_{\delta}^{\Theta}$ satisfying
$$
\|(L_{1}\vee 1)\cdot \xi_{\delta}^{G}\|_{\MP^{Z},2}\leq \overline{M}_{1}~\delta^{\beta}\quad\mbox{and}\quad J(\FFF_{\delta}^{\Theta},\xi_{\delta}^{G},1/8)\vee 
J(\FFF^{\Theta},\xi^{G},1/8) \leq \overline{M}^{1} 
$$
for $\delta\in ]0,\Delta(\Theta)]$.
\end{enumerate}
Condition \ref{assump:A5} plays the role to find a suitable nonasymptotic rate $n^{\gamma}$ in Proposition \ref{variational inequalities}, and thus also for the deviation probabilities in \eqref{deviation}. It will turn out that the choice 
$n^{\gamma}$ with $\gamma = 1/(4-2\beta)$ fits, where $\beta\in ]0,1]$ as in \ref{assump:A5}.
%
\smallskip
 
Finally, for simplification we additionally require the following representation of the functions $\EEE[H(\theta,\cdot,Z_{1})]$: 
\begin{enumerate}[label=(A \arabic*), ref=(A \arabic*)]
\setcounter{enumi}{7}
\item \label{assump:A6} 
There exist $\delta_{1} > 0, K_{\delta_{1}}\geq 0$
and a Borel measurable mapping 
$m_{\Theta\times\RRR}: \Theta\times\RRR\rightarrow\RRR$ which 
satisfies 
\begin{align*}
&
\big|m_{\Theta\times\RRR}\big(\theta,\EEE[G(\theta,Z_{1})] + x\big) - m_{\Theta\times\RRR}\big(\theta^{*},\EEE[G(\theta^{*},Z_{1})] + x\big)\big|\\
&
\leq K_{\delta_{1}}~\big\|\big(\theta - \theta^{*},\EEE[G(\theta,Z_{1})]-\EEE[G(\theta^{*},Z_{1})]\big)\big\|_{m+1}\quad\mbox{for}~(\theta,x)\in\Theta\times [-\delta_{1},\delta_{1}],
\end{align*}
and
$$
\EEE[H(\theta,t,Z_{1})] - \EEE[H(\theta,s,Z_{1})] = \int_{s}^{t}m_{\Theta\times\RRR}(\theta,u)~du\quad\mbox{for}~(\theta,t), (\theta,s)\in\Theta\times\RRR, t > s.
$$
\end{enumerate}
A simplier condition ensuring together with \ref{assump:A5} the assumption \ref{assump:A6} is provided in Lemma \ref{simplier A6} below.
\medskip

Before further developping our main results we want to present a simple example where all the assumptions \ref{assump:A1} - \ref{assump:A6} are fulfilled simultaneously.
\begin{example}
\label{simple example}
Let $m = d$ and let $\Theta$ be a convex compact subset of $\RRR^{m}\setminus\{0\}$. Furthermore $Z$ is a centered $d$-dimensional, normally distributed random vector with positive definite covariance matrix $\Sigma$, and we fix a random variable $W$ which is chi-squared distributed with $1$ degree of freedom. By compactness of $\Theta$ we may find some open convex subset $\cU$ of $\RRR^{m}\setminus\{0\}$ enclosing $\Theta$. We consider the following stochastic program
$$
\min_{\theta\in\Theta}\EEE\big[\big((\theta^{T}Z)^{2} - \EEE[(\theta^{T}Z)^{2}]\big)^{+}\big].
$$
It is a specialization of the optimization problem \eqref{composite optimization} with mappings $G$ and $H$ defined by
\begin{align*}
G(\theta,z) := (\theta^{T}z)^{2}\quad\mbox{and}\quad H(\theta,t,z) := \big((\theta^{T}Z)^{2} - t\big)^{+}.
\end{align*}
\begin{itemize}
\item [1)] $G$ is continuous, implying immediately \ref{assump:A1}. 
A square $\MP^{Z}$-integrable positive envelope $\xi^{G}$ of $\FFF^{\Theta}$ is defined by $\xi^{G}(z) := \max_{\theta\in\Theta}\|\theta\|_{m}^{2}~[\lambda^{\Sigma}~z^{T}\Sigma^{-1} z +\eins_{\{0\}}(z)]$, where $\lambda^{\Sigma}$ stands for the maximal eigenvalue of $\Sigma$. Then by the continuity of $G$ along with the dominated convergence theorem the condition \ref{assump:A3} follows easily.
\item [2)] $H(\theta,t,\cdot)$ is $\MP^{Z}$-integrable for $\theta\in\Theta$ and $t\in\RRR$ due to square integrability of $Z_{1}^{T}\Sigma^{-1} Z_{1}$. Furthermore property \ref{assump:A4} holds with $L_{1} = L_{2} :\equiv 2$.
\item [3)] We have $\EEE\big[H\big(\theta,\EEE[G(\theta,Z)],Z\big)\big] = 
\EEE\big[\big|G(\theta,Z) - \EEE[G(\theta,Z)]\big|\big]/2$ for $\theta\in\Theta$ 
(see \cite[Proposition 6.1]{ShapiroEtAl}). 
\begin{itemize}
\item [a)] Since $\theta^{T}Z$ is a centered normally distributed random variable with variance $\theta^{T}\Sigma\theta$, we obtain via the change of variable formula
\begin{align*}
\EEE\big[H\big(\theta,\EEE[G(\theta,Z)],Z\big)\big] 
&= 
\int_{\RRR}|x^{2} - \theta^{T}\Sigma\theta|~\dfrac{1}{2~\sqrt{2\pi~\theta^{T}\Sigma\theta}}~\exp\left(-\dfrac{x^{2}}{2 \theta^{T}\Sigma\theta}\right)~dx\\
&= 
\dfrac{\theta^{T}\Sigma\theta}{2}~
\int_{\RRR}|w^{2} - 1|~\dfrac{1}{\sqrt{2\pi}}~\exp\left(-\dfrac{w^{2}}{2}\right)~dw\\ 
&= 
\theta^{T}\Sigma\theta~\dfrac{\EEE[|W - 1|]}{2}\quad\mbox{for}~\theta\in\Theta.
\end{align*}
\item [b)] The mapping $\widehat{\psi}:\cU\rightarrow \RRR$, defined by $\widehat{\psi}(y) := y^{T}~\Sigma~y~\EEE[|W - 1|]/2$, is strictly convex and thus also continuous because $\Sigma$ is positive definite. Hence the goal function $\psi_{H,\Theta}$ of the stochastic program has a unique minimum $\theta^{*}$.
\item [c)] The mapping $\widehat{\psi}$ is also twice continuously differentiable so that by Taylor approximation
\begin{align*}
\psi_{H,\Theta}(\theta) - \psi_{H,\Theta}(\theta^{*}) 
&= 
\widehat{\psi}(\theta) - \widehat{\psi}(\theta^{*})\\ 
&= \EEE[|W - 1|]~\left[(\theta - \theta^{*})^{T}\Sigma\theta^{*} + (\theta - \theta^{*})^{T}\Sigma (\theta -\theta^{*})\right]
\end{align*}
holds for $\theta\in\Theta$. Moreover, $\theta^{*} + t~(\theta- \theta^{*})\in\Theta$ for $\theta\in\Theta$ and $t\in [0,1]$ due to convexity of $\Theta$. Since $\theta^{
*}$ minimizes $\psi_{H,\Theta}$ we may further observe
\begin{align*}
\EEE[|W - 1|]~(\theta - \theta^{*})^{T}\Sigma\theta^{*} 
&= 
\lim_{t\searrow 0}\dfrac{\widehat{\psi}(\theta^{*} + t [\theta - \theta^{*}]) - \widehat{\psi}(\theta^{*})}{t}\\ 
&= 
\lim_{t\searrow 0}\dfrac{\psi_{H,\Theta}(\theta^{*} + t [\theta - \theta^{*}]) -\psi_{H,\Theta}(\theta^{*}) }{t}\geq 0 
\end{align*}
for $\theta\in\Theta$. As a consequence the second order growth condition \ref{assump:A7} is fulfilled with $M_{1} := \EEE[|W - 1|]~\lambda_{\Sigma}$, where $\lambda_{\Sigma}$ stands for the minimal Eigenvalue of $\Sigma$.
\end{itemize}
\item [4)]
Let $\pi_{1},\ldots,\pi_{m}$ denote the coordinate mappings from $\RRR^{m}$ onto $\RRR$.
\begin{itemize}
\item [a)] The linear hulls of $\FFF^{\Theta}$ and $\FFF^{\Theta}_{\delta}$ $(\delta\in ]0,\Delta(\Theta)])$ are generated by the mappings $\pi_{i}\cdot\pi_{j}$, where $i\in\{1,\ldots,m\}$ and $j\in\{i,\ldots,m\}$. In particular, all these linear hulls have dimension of at most $m(m+1)/2$. Then by Example \ref{finite dimension} in Subsection \ref{Beweis von Proposition important covering numbers}
\begin{align*}
&
J(\FFF^{\Theta},\xi^{G},\eta)\vee J(\FFF^{\Theta}_{\delta},\xi_{\delta},\eta)\\
&\leq 2~\eta~\sqrt{\ln(e~m[m +1] + 4~e) + [m(m+1)/2 + 1]~\ln(16\cdot e)}
\end{align*}
for $\eta\in ]0,1[, \delta\in ]0,\Delta(\Theta)]$, and any positive envelope $\xi_{\delta}$ of $\FFF^{\Theta}_{\delta}$.
\item [b)] Denoting the maximal eigenvalue of the matrix $\Sigma$ by $\lambda^{\Sigma}$, we may define via $\xi^{G}_{\delta}(z) := 2~\lambda^{\Sigma}~\max_{\theta\in\Theta}\|\theta\|_{m}~[\lambda^{\Sigma}~z^{T}\Sigma^{-1} z + \eins_{\{0\}}(z)]~\delta$ a family $\big(\xi^{G}_{\delta}\big)_{\delta\in ]0,\Delta]}$ of positive envelopes as required in assumption \ref{assump:A5} with $\beta = 1$ and the constant $\overline{M}_{1} := 4~\lambda^{\Sigma}~\max_{\theta\in\Theta}\|\theta\|_{m}~\sqrt{m (m + 2)}$.  
\item [c)] In view of a) we may choose in \ref{assump:A5} the constant 
$$
\overline{M}^{1} := \sqrt{\ln(e~m[m +1] + 4~e) + [m(m+1)/2 + 1]~\ln(16\cdot e)}/4
$$
\end{itemize}
\item [5)] For condition \ref{assump:A6} we may choose the function $m_{\Theta\times\RRR}:\Theta\times\RRR\rightarrow\RRR$ defined by $m_{\Theta\times\RRR}(\theta,x) :=  F_{\theta}(x) - 1$.
\begin{itemize}
\item [a)] The random variable $(\theta^{T}Z)^{2}$ has the same distribution as 
$\theta^{T}\Sigma\theta~W$ for every $\theta\in \Theta$. Hence, denoting the distribution function of $W$ by $F_{W}$, we obtain
$
F_{\theta}(x) - 1 =  F_{W}\big(x/\theta^{T}\Sigma\theta\big) - 1\quad\mbox{for}~\theta\in\Theta, x\in\RRR
$.
\item [b)] Since $F_{W}$ is continuously differentiable, and since $y^{T}~\Sigma~y > 0$ for $y\in\cU$, the mapping 
$$
m_{\cU\times\RRR}:\cU\times\RRR\rightarrow\RRR,~(y,x)\mapsto F_{W}\big(x/y^{T}~\Sigma~y\big) - 1
$$
is continuously differentiable. In particular by Taylor approximation theorem, for any nonvoid compact interval $\mathcal{I}$ in $\RRR$, the restriction of $m_{\cU\times\RRR}$ to 
$\Theta\times\mathcal{I}$ is Lipschitz continuous with Lipschitz constant
$$
L_{\mathcal{I}} := \max_{(\theta,x)\in\Theta\times\mathcal{I}}\|\nabla_{(\theta,x)}m_{\cU\times\RRR}\|_{m+1},
$$
where $\nabla_{(\theta,x)}m_{\cU\times\RRR}$ denotes the gradient of $m_{\cU\times\RRR}$ at $(\theta,x)$.
\item [c)] The mapping $\theta\mapsto\EEE[G(\theta,Z)] = \theta^{T}\Sigma\theta$ on $\Theta$ is continuous. Therefore, the interval $I_{\delta_{1}} := [\inf_{\theta\in\Theta}\theta^{T}\Sigma\theta - \delta_{1},\sup_{\theta\in\Theta}\theta^{T}\Sigma\theta + \delta_{1}]$ is compact for $\delta_{1} > 0$. We may conclude that condition \ref{assump:A6} is fulfilled for any $\delta_{1} > 0$ with $K_{\delta_{1}} = L_{\mathcal{I}_{\delta_{1}}}$.
\end{itemize}
\end{itemize}
\end{example}
\medskip

As one main idea we shall use a two step truncation procedure to derive upper estimates for the deviation probabilities \eqref{deviation}. Defining
\begin{align}
& \label{Hilfsereignis 1}
\Omega_{n,a} := \big\{\sup_{\theta\in\Theta}|(\MP_{n} - \MP)\big(G(\theta,\cdot)\big)| \leq a/\sqrt{n}\big\}\quad(n\in\NNN, a > 0),
\end{align}
the aim is to provide separately upper bounds for 
$
\MP(\{\|\widehat{\theta}_{n} - \theta^{*}\| > \varepsilon\cdot n^{-\gamma}\}\cap \Omega_{n,a})
$
and $\MP(\Omega\setminus\Omega_{n,a})$. Concerning the truncated deviations we shall invoke Propositions \ref{variational inequalities}, and we consider a further truncation, restricting the increments of the compound empirical processes to the following auxiliary events
\begin{align}
&\label{Hilfsereignis 2}
\Omega^{n}_{\delta,b} := \left\{\sup_{\theta\in\cU_{\delta}}
|(\MP_{n} - \MP)\big(G(\theta,\cdot) - G(\theta^{*},\cdot)\big)| \leq \frac{b~\delta^{\beta}}{\sqrt{n}}\right\} 
\quad(n\in\NNN, \beta\in ]0,1]; b, \delta > 0).
\end{align}
So actually, what we shall do is to upper estimate separately 
\begin{align*}
\MP^{*}\Big(\Big\{\sup_{\theta\in U_{\widehat{\delta}_{nk}}}\big|\GGG_{n}^{G,H}(\theta,\cdot) - \GGG_{n}^{G,H}(\theta^{*},\cdot) \big| > 
M_{1}~2^{2 (k - 1)}~n^{- 2\gamma + 1/2} 
\Big\}\cap\Omega_{n,a}\cap\Omega^{n}_{\widehat{\delta}_{nk},b_{k}}\Big)
\end{align*}
and $\MP(\Omega\setminus \Omega^{n}_{\widehat{\delta}_{nk},b_{k}})$ for arbitrary $a > 0$ and $k$ varying over a certain finite set of integers, where 
$\widehat{\delta}_{nk} := \Delta(\Theta)\wedge (2^{k}/n^{\gamma})$, and $b_{k}$ is a specific positive constant dependent on $k$.
\par
Note that in view of \ref{assump:A3} along with the completeness of $\OFP$ any $\Omega_{n,a}$ and $\Omega_{\delta,b}^{n}$ belong to $\cF$ (see Lemma \ref{messbare Hilfsereignisse} below).
\medskip

Recalling mapping $L_{2}$ from \ref{assump:A4}, the constants $\overline{M}_{1}, \overline{M}^{1}$ from 
\ref{assump:A5} and 
$\delta_{1} > 0$, $K_{\delta_{1}}\geq 0$ from \ref{assump:A6}, 
the following constants will play an import role in our main results:
\begin{align}
&\label{Mab}
M_{a,b} := K_{\delta_{1}}~a~[\Delta(\Theta)^{1-\beta} + \overline{M}_{1}] + ~b~\EEE[L_{2}(Z_{1})]\quad(a, b > 0),\\
&\label{Mb}
\overline{M}_{b} := \overline{M}_{1}~(1 + \|L_{2}\|_{\MP^{Z},2}) + b~\|L_{2}\|_{\MP^{Z},2}\quad(b > 0),\\
&\label{frakg}
\mathfrak{g}(t) := 
\frac{t^{2}}{8(t + 1)(5 t + 28)}\quad(t > 0).
\end{align}
For $a, b,\delta > 0, n\in\NNN$ we also define the following abbreviations  
\begin{align}
&\label{etantab}
\overline{\eta}_{n}(a,b,\delta) := 
64 ~\sqrt{2}~ \overline{M}^{1} + 32 
\sqrt{
\ln(24) + \big[\ln (a/b) - \ln(\sqrt{n}~\delta^{\beta})\big]^{+}
},\\
&\label{etadachepsilon}
\widehat{\eta} := 64~\sqrt{2}~(\overline{M}^{1} + 0,7),\\
& \label{etaabbeta}
n(a,b,\beta) := \frac{\overline{M}_{b}^{2}~\Delta(\Theta)^{2\beta}}{2}\vee \frac{\eins_{]0,\infty[}\big(\EEE[L_{2}(Z_{1})]\big)~a^{2}}{\delta_{1}^{2}}\vee \big[\overline{M}_{b}^{2}~\Delta(\Theta)^{2\beta + 1}]^{\frac{4 - 2\beta}{3 - 2\beta}}.
\end{align}

Now, we are ready to formulate our results on upper estimates of the probabilities in \eqref{deviation}. They will be described in terms of explicit constants. We start with a rough estimation.
\begin{theorem}
\label{simpliest version rough}
Let \eqref{unique solution} and assumptions \ref{assump:A1} - \ref{assump:A6} be fulfilled with 
\begin{itemize}
\item $L_{2}$ denoting the square $\MP^{Z}$-integrable mapping from \ref{assump:A4},
\item the constant $M_{1} > 0$ as in \ref{assump:A7},
\item constant $\beta\in ]0,1]$ from \ref{assump:A5}.
\end{itemize}
Finally, let $a,b,\varepsilon > 0$ with $a~2^{\beta + 1/2}\leq b\varepsilon^{\beta + 1/2}$, 
and let $M_{a,b},\overline{M}^{b},\widehat{n},n(a,b,\beta)$ be the constants introduced respectively in \eqref{Mab}, \eqref{Mb}, \eqref{etadachepsilon}, \eqref{etaabbeta}.


Then, setting $\delta_{nk} := 2^{k}n^{-1/(4- 2\beta)}$ for $k\in\mathbb{Z}, n\in\NNN$, and using notations for the auxiliary events defined in \eqref{Hilfsereignis 1} and \eqref{Hilfsereignis 2}, 
\begin{align*}
\MP\big(\big\{n^{1/(4 - 2\beta)}~ \|\widehat{\theta}_{n} - \theta^{*}\|_{m} > \varepsilon\big\}\big) \leq P^{\varepsilon}_{1} + P^{\varepsilon}_{n 2} + \MP\big(\Omega\setminus\Omega_{n,a}\big)
\end{align*}
for $n\in\NNN$ with $n\geq n(a,b,\beta)$, where
\begin{align*}
&P^{\varepsilon}_{1} := \frac{10 [\overline{M}_{b}~\widehat{\eta} + \eins_{]0,\infty[}\big(\EEE[L_{2}(Z_{1})]\big)~M_{a,b}]}{M_{1}}~\left(\frac{2}{\varepsilon}\right)^{3/2 - \beta} \sqrt{\frac{2}{\varepsilon\wedge 2}},\\
&P^{\varepsilon}_{n 2} := \sum_{k= K_{\varepsilon} + 1\atop \delta_{n (k-1)}\leq\Delta(\Theta)}^{\infty}\hspace*{-0.35cm}\MP\big(\Omega\setminus\Omega^{n}_{\delta_{nk}\wedge\Delta(\Theta),\sqrt{2}^{k} b}\big).
\end{align*}
\end{theorem}
We may improve the estimations in Theorem \ref{simpliest version rough} if the  square $\MP^{Z}$-integrable mappings $\xi^{G}_{\delta}$ from assumption assumption \ref{assump:A5} have weak tails.  
For preparation, we shall associate any square $\MP^{Z}$-integrable mapping $\xi$ with the following events
\begin{equation}
\label{Restmenge}
B_{n}^{\xi} := \Big\{\frac{1}{n}\sum_{j=1}^{n}\xi(Z_{j})^{2}\leq  2\EEE[\xi(Z_{1})^{2}]\Big\}\quad(n\in\NNN).
\end{equation}
The idea now is to further restrict the increments of $\GGG^{G,H}_{n}$ to events of the form $B_{n}^{\xi^{G}_{a,b,\delta}}$ for $a,b > 0,\delta\in ]0,\Delta(\Theta)])$, where 
\begin{equation}
\label{neue envelope}
\xi^{G}_{a, b,\delta} := \sqrt{(L_{1}(\cdot)\vee 1)^{2}\cdot \xi^{G}_{\delta}(\cdot)^{2} + L_{2}(\cdot)^{2}\cdot (b + \overline{M}_{1})^{2}\cdot \delta^{2\beta}},
\end{equation}
and $L_{1}, L_{2}$ are the $\MP^{Z}$-square integrable mappings from \ref{assump:A4}, whereas $\overline{M}_{1} > 0, \beta\in ]0,1]$ are as in \ref{assump:A5}.
\begin{theorem}
\label{simpliest version}
Let \eqref{unique solution} and assumptions \ref{assump:A1} - \ref{assump:A6} be fulfilled with
\begin{itemize}
\item $\xi^{G}$ standing for the positive envelope of $\FFF^{\Theta}$ as in \ref{assump:A2},
\item $L_{1}, L_{2}$ denoting the square $\MP^{Z}$-integrable mappings from condition \ref{assump:A4}
\item constant $M_{1} > 0$ as in \ref{assump:A7},
\item $\big(\xi^{G}_{\delta}\big)_{\delta\in ]0,\Delta(\Theta)]}$ being the family of positive envelopes from \ref{assump:A5},
\item constants $\overline{M}_{1} > 0$ as well as $\beta\in ]0,1]$ from \ref{assump:A5}.
\end{itemize} 
Moreover, let $a, b, \varepsilon > 0$ with $a~2^{\beta + 1/2}\leq b~\varepsilon^{\beta + 1/2}$, 
and let $M_{a,b}, \overline{M}_{b},\mathfrak{g}(t)$, $\widehat{n}$, $n(a,b,\beta)$ be the constants introduced respectively in \eqref{Mab}, \eqref{Mb}, \eqref{frakg}, \eqref{etadachepsilon}, \eqref{etaabbeta}. Finally, $\xi^{G}_{a,b,\delta}$ denotes for $\delta > 0$ the 
square $\MP^{Z}$-integrable mapping from \eqref{neue envelope}.
\smallskip

Then, setting $\delta_{nk} := 2^{k}n^{-1/(4- 2\beta)}$ for $k\in\mathbb{Z},n\in\NNN$, and using the notations for the auxiliary events defined in \eqref{Hilfsereignis 1}, \eqref{Hilfsereignis 2} and \eqref{Restmenge}, 
the following statement holds:
\begin{itemize}
\item [~] If $n\in\NNN$ with $n\geq n(a,b,\beta)$, and if for $t > 0$ the inequality 
$$
2^{K_{\varepsilon} (2 - \beta) - K_{\varepsilon}^{+}/2} > 4 ~\big[M_{a,b}~\eins_{]0,\infty[}\big(L_{2}(Z_{1})\big) + \overline{M}_{b}~\big(1 + 2 (t + 1)~\widehat{\eta}\big)\big]/M_{1}
$$ 
holds, then
\begin{align*}
&
\MP\big(\big\{n^{1/(4 - 2\beta)}~ \|\widehat{\theta}_{n} - \theta^{*}\|_{m} > \varepsilon\big\}\big)\leq P^{t,\varepsilon}_{1} + P^{\varepsilon}_{n 2} + P^{\varepsilon}_{n 3} + \MP\big(\Omega\setminus  \Omega_{n,a}\big),
\end{align*} 
where
\begin{align*}
&P^{t,\varepsilon}_{1} :=
\exp\left(\frac{\eins_{]0,\infty[}\big(\EEE[L_{2}(Z_{1})]\big)~M_{a,b}~\mathfrak{g}(t)}{(\sqrt{\varepsilon}\wedge 1)~b~\|L_{2}\|_{\MP^{Z},2}}\right)\cdot \frac{2^{4- \beta}~\overline{M}_{b}\cdot\exp\left(\dfrac{-\varepsilon^{3/2- \beta}\sqrt{\varepsilon\wedge 2}~M_{1}~\mathfrak{g}(t)}{2^{4-\beta}~\overline{M}_{b}}\right)}{(3/2-\beta)\ln(2)~M_{1}~\varepsilon^{3/2-\beta}\sqrt{\varepsilon\wedge 2}}\\
&P^{\varepsilon}_{n 2} := 
\sum_{k= K_{\varepsilon} + 1\atop \delta_{n (k-1)}\leq\Delta(\Theta)}^{\infty}\hspace*{-0.65cm}\MP\big(\Omega\setminus\Omega^{n}_{\delta_{nk}\wedge\Delta(\Theta),\sqrt{2}^{k} b}\big)\quad\mbox{and}\quad
P^{\varepsilon}_{n 3} := 
\sum_{k= K_{\varepsilon} + 1\atop \delta_{n (k-1)}\leq\Delta(\Theta)}^{\infty}\hspace*{-0.65cm}\MP\big(\Omega\setminus B_{n}^{\xi^{G}_{a,\sqrt{2}^{k} b,[\delta_{nk}\wedge\Delta(\Theta)]}}\big).
\end{align*}
\end{itemize}

\end{theorem}
In the next step we want to complete the results on the deviation probabilities \eqref{deviation} by upper estimations for the probabilities $\MP(\Omega\setminus\Omega_{n,a}), \MP(\Omega\setminus\Omega^{n}_{\delta,b})$, and the values $P_{n 2}^{\varepsilon}$ in Theorems \ref{simpliest version rough}, \ref{simpliest version}.
\begin{remark}
\label{erste Restwahrscheinlichkeiten}
Let assumptions \ref{assump:A1}, \ref{assump:A2}, \ref{assump:A3} and \ref{assump:A5} be fulfilled with 
\begin{itemize}
\item square $\MP^{Z}$- integrable mapping $\xi^{G}$ from \ref{assump:A2};
\item family $\big(\xi^{G}_{\delta})_{\delta\in ]0,\Delta(\Theta)]}$ of square $\MP^{Z}$- integrable mappings from  \ref{assump:A5};
\item constants $\overline{M}_{1},\overline{M}^{1} > 0$ and $\beta\in ]0,1]$ as in \ref{assump:A5}.
\end{itemize}
Furthermore, let $\mathfrak{g}(t)$ be the constant defined in \eqref{frakg}. Then $\Omega_{n,a}, \Omega^{n}_{\delta,b}$ belong to $\cF$ for $a,b,\delta > 0$. Moreover, using notations for auxiliary events introduced in \eqref{Hilfsereignis 1}, \eqref{Hilfsereignis 2}, and \eqref{Restmenge}, we obtain for any sample size the following upper bounds of the probabilities $\MP(\Omega\setminus\Omega_{n,a}), \MP(\Omega\setminus\Omega^{n}_{\delta,b})$ and the values $P^{\varepsilon}_{n 2}$ in Theorems \ref{simpliest version rough}, \ref{simpliest version}.
\begin{itemize}
\item [1)]$\MP\big(\Omega\setminus\Omega_{n,a}\big)\leq 64\sqrt{2}~\overline{M}^{1}~\|\xi^{G}\|_{\MP^{Z},2}/a$ for $n\in\NNN$, 
\item [2)]$\MP\big(\Omega\setminus\Omega^{n}_{\delta,b}\big)\leq 64\sqrt{2}~\overline{M}_{1}~\overline{M}^{1}/b$ for $n\in\NNN$ and $\delta\in ]0,\Delta(\Theta)]$.
\item [3)] 
For $\varepsilon > 0$, and setting $\delta_{nk} := 2^{k}~n^{-1/(4-2\beta)}$
$$
\sum\limits_{k= K_{\varepsilon} + 1\atop \delta_{n (k-1)}\leq\Delta(\Theta)}^{\infty}\hspace*{-0.35cm}\MP\big(\Omega\setminus\Omega^{n}_{\delta_{nk}\wedge\Delta(\Theta),\sqrt{2}^{k} b}\big)\leq \frac{128~\overline{M}_{1}~\overline{M}^{1}}{(\sqrt{2} - 1)~b~\sqrt{\varepsilon}}.
$$
\end{itemize}
These bounds may be tightened in the following way.
\begin{itemize}
\item [4)] If $t > 0$ such that $a > \|\xi^{G}\|_{\MP,2}~[1 + 64\sqrt{2}~(t + 1)~\overline{M}^{1}]$, and if $n\in\NNN$ satisfies $n\geq\|\xi^{G}\|_{\MP,2}^{2}/2$, then
$$
\MP(\Omega\setminus\Omega_{n,a})\leq 
\exp\Big(-\frac{ a\cdot\mathfrak{g}(t)}{\|\xi^{G}\|_{\MP^{Z},2}}\Big) + 
\MP\big(\Omega\setminus B_{n}^{\xi^{G}}\big).
$$
\item [5)] For $t > 0, \delta\in ]0,\Delta(\Theta)]$ and $n\in\NNN$ with $n\geq \overline{M}_{1}^{2}~\delta^{2\beta}/2$
$$
\MP(\Omega\setminus\Omega^{n}_{\delta,b})\leq \exp\Big(- \frac{b\cdot\mathfrak{g}(t)}{\overline{M}_{1}}\Big) + \MP\big(\Omega\setminus B_{n}^{\xi^{G}_{\delta}}\big)
$$
is valid whenever $b > \overline{M}_{1}~[1 + 64\sqrt{2}~(t + 1)~\overline{M}^{1}]$.
\item [6)]If $n\in\NNN$ with $n\geq \overline{M}_{1}^{2}~\Delta(\Theta)^{2\beta}/2$, then for $t, \varepsilon > 0$ satisfying the inequality $b > \overline{M}_{1}~[1 + 64\sqrt{2}~(t + 1)~\overline{M}^{1}]/\sqrt{\varepsilon}$ 
\begin{align*}
&\sum_{k= K_{\varepsilon} + 1\atop \delta_{n (k-1)}\leq\Delta(\Theta)}^{\infty}\hspace*{-0.35cm}\MP\big(\Omega\setminus\Omega^{n}_{\delta_{nk}\wedge\Delta(\Theta),\sqrt{2}^{k} b}\big)
\\
&\leq 
\frac{2~\sqrt{2}~\overline{M}_{1}}{\sqrt{\varepsilon}~\ln(2)~b~\mathfrak{g}(t)}
~\exp\left(-\frac{\sqrt{\varepsilon}~b~\mathfrak{g}(t)}{\sqrt{2}~\overline{M}_{1}}\right)
+ \sum_{k= K_{\varepsilon} + 1\atop \delta_{n (k-1)}\leq\Delta(\Theta)}^{\infty}\hspace*{-0.35cm}\MP\big(\Omega\setminus B_{n}^{\xi^{G}_{\delta_{nk}\wedge\Delta(\Theta)}}\big).
\end{align*}
\end{itemize}
The proof may be found in Subsection \ref{Restwahrscheinlichkeiten}.
\end{remark}
In the following we want to discuss possible further upper estimations for the probabilities 
$\MP\big(\Omega\setminus B_{n}^{\xi^{G}}\big),\MP\big(\Omega\setminus B_{n}^{\xi^{G}_{\delta}}\big), \MP\big(\Omega\setminus B_{n}^{\xi^{\Theta}_{a,b,\delta}}\big)$ in Theorem \ref{simpliest version} and Remark \ref{erste Restwahrscheinlichkeiten}.
\begin{remark}
\label{Simplifications}
Let $\xi$ denote any square $\MP^{Z}$-integrable mapping with event $B_{n}^{\xi}$ as defined in \eqref{Restmenge}. The following upper estimates of $\MP\big(\Omega\setminus B_{n}^{\xi}\big)$ might be used to make the upper estimations in Theorems \ref{simpliest version rough}, \ref{simpliest version} more explicit.
\begin{itemize}
\item [1)] If the function $\xi$ is bounded by some positive constant $L$, then $\xi \equiv L$ may be chosen so that $\Omega\setminus B_{n}^{\xi} = \emptyset$.
\item [2)] If $\xi$ is $\MP^{Z}$-integrable of order $4$, we may apply Cantelli's inequality to conclude
$$
\MP\big(\Omega\setminus B_{n}^{\xi}\big) \leq \dfrac{\vari[\xi(Z_{1})^{2}]}{  n~\EEE[\xi(Z_{1})^{2}]^{2} + \vari[\xi(Z_{1})^{2}]}\quad\mbox{for}~n\in\NNN.
$$ 
\item [3)]
The upper estimate of the probability $\MP\big(\Omega\setminus B_{n}^{\xi}\big)$ may be further improved if the random variable $\exp\big(\lambda\cdot\xi^{2\alpha }\big)$ is $\MP^{Z}$-integrable for some $\lambda > 0$ and $\alpha\in ]0,1]$. This means that $\xi(Z_{1})^{2}$ has a so called \textit{sub-Weibull} distribution of order $\alpha$ (see \cite{KuchibhotlaChakrabortty2022}). In this case, setting 
$|\xi^{2}|_{\alpha} := \inf\{c > 0\mid \EEE\big[\exp\big(\xi(Z_{1})^{2\alpha}/c^{\alpha}\big)\big]\leq 2\}$ and 
$$
c(\alpha) := 2^{2/\alpha}~\big[1 + \big(d_{\alpha}~\ln(2)\big)^{-1/\alpha}\big]~\sqrt{8}~e^{3}~(2 \pi)^{1/4}~e^{1/24}~\big(\exp(2/e)/\alpha\big)^{1/\alpha}
$$
with $d_{\alpha} := (\alpha~e)^{1/\alpha}/2$,
the inequality
$$
\EEE\Big[\Big|\frac{1}{n}~\sum_{j=1}^{n}\xi(Z_{j})^{2} - \EEE\big[\xi(Z_{1})^{2}\big]\Big|^{p}\Big]
\leq 
|\xi^{2}|_{\alpha}^{p}~c(\alpha)^{p}~\left[\sqrt{\frac{p}{n}} + \frac{p^{1/\alpha}}{n}\right]^{p}
$$
holds for $p\in [1,\infty[$ (see \cite[Proof of Theorem 3.1]{KuchibhotlaChakrabortty2022} along with \cite[(A.6) + (A.1)]{GoetzeEtAl2021}).
Then we may invoke Proposition 3.3 from \cite{GoetzeEtAl2021} to obtain
$$
\MP\big(\Omega\setminus B_{n}^{\xi}\big) \leq 2~\exp\left(- \frac{\ln(2)}{8~e^{2}}~\left[\dfrac{\EEE[\xi(Z_{1})^{2}]^{2}~n}{|\xi^{2}|_{\alpha}^{2}~c(\alpha)^{2}}\wedge \dfrac{\EEE[\xi(Z_{1})^{2}]^{\alpha}~n^{\alpha}}{|\xi^{2}|_{\alpha}^{\alpha}~c(\alpha)^{\alpha}}\right]\right)\quad\mbox{for}~n\in\NNN.
$$
Note additionally the bounds 
\begin{align*}
\dfrac{(\alpha~e)^{1/\alpha}}{2}\sup_{p\in ]0,1]}~\dfrac{\big(\EEE[\xi(Z_{1})^{2 p}]\big)^{1/p}}{p^{\alpha}}\leq 
|\xi^{2}|_{\alpha}\leq (2~e)^{1\alpha}~\sup_{p\in ]0,1]}~\dfrac{\big(\EEE[\xi(Z_{1})^{2 p}]\big)^{1/p}}{p^{\alpha}}
\end{align*}
for $|\xi^{2}|_{\alpha}$ (see \cite[Lemma A.2]{GoetzeEtAl2021}).
\item [4)]
Let $\xi := \sqrt{t_{1} \xi_{1}^{2} + t_{2}\xi_{2}^{2}}$ for some $t_{1},t_{2} >0$ and square $\MP^{Z}$-integrable mappings $\xi_{1}, \xi_{2}$. Examples of such $\xi$ are provided by the mappings $\xi^{\Theta}_{a,b,\delta}$ in \eqref{neue envelope}. Then 
$$
\MP\big(\Omega\setminus B_{n}^{\xi}\big)\leq\sum_{i=1}^{2}\MP\big(\Omega\setminus B_{n}^{\xi_{i}}\big)\quad\mbox{for}~n\in\NNN.
$$
In particular, we might apply the results from 1) - 3) separately to $\xi_{1}$ and $\xi_{2}$. 
\end{itemize}
\end{remark}
\begin{remark}
\label{unangenehmer Rest}
In order to control in Theorem \ref{simpliest version} and Remark \ref{erste Restwahrscheinlichkeiten} the probabilities 
\begin{align*}
\sum_{k= K_{\varepsilon} + 1\atop \delta_{n (k-1)}\leq\Delta(\Theta)}^{\infty}\hspace*{-0.35cm}\MP\big(\Omega\setminus B_{n}^{\xi^{G}_{a,\sqrt{2}^{k} b,[\delta_{nk}\wedge\Delta(\Theta)]}}\big)\quad\mbox{and}\quad \sum_{k= K_{\varepsilon} + 1\atop \delta_{n (k-1)}\leq\Delta(\Theta)}^{\infty}\hspace*{-0.35cm}\MP\big(\Omega\setminus B_{n}^{\xi^{G}_{\delta_{nk}\wedge\Delta(\Theta)}}\big)
\end{align*}
devices from Remark \ref{Simplifications} may be utilized. In particular with $\sum_{\emptyset} := 0$ the observation
\begin{align*}
&
\sum_{k= K_{\varepsilon} + 1\atop \delta_{n (k-1)}\leq\Delta(\Theta)}^{\infty}\hspace*{-0.35cm}\MP\big(\Omega\setminus B_{n}^{\xi^{G}_{a,\sqrt{2}^{k} b,[\delta_{nk}\wedge\Delta(\Theta)]}}\big)\\
&\leq 
\sum_{k= K_{\varepsilon} + 1\atop \delta_{n (k-1)}\leq\Delta(\Theta)}^{\infty}\hspace*{-0.35cm}\MP\big(\Omega\setminus B_{n}^{(L_{1}\vee 1)~\xi^{G}_{\delta_{nk}\wedge\Delta(\Theta)}}\big) 
+ \frac{\big[(4 - 2\beta)~\ln\big(4 \Delta(\Theta)/\varepsilon\big) + \ln(n)\big]^{+}}{(4 - 2\beta)~\ln(2)}~\MP\big(\Omega\setminus B_{n}^{L_{2}}\big)
\end{align*} 
might be useful.
\end{remark}
In principle we may construct nonasymptotic confidence regions for fixed sample size by reassembling the inequalities in Theorems \ref{simpliest version rough}, \ref{simpliest version} together with Remarks \ref{erste Restwahrscheinlichkeiten} - \ref{unangenehmer Rest}. A general procedure is provided in Section \ref{confidence regions}.
\bigskip

As another application of Theorem \ref{simpliest version rough} we may draw on Remark \ref{erste Restwahrscheinlichkeiten} to derive the following nonasymptotic rates for the minimizers of the SAA problems.
\begin{corollary}
\label{allgemeine Konvergenzraten zweiter Teil}
Let (\ref{unique solution}) and assumptions \ref{assump:A1} - \ref{assump:A6} be fulfilled with $\beta\in ]0,1]$ from \ref{assump:A5}, and let $\big(\widehat{\theta}_{n}\big)_{n\in\NNN}$ be a sequence of minimizers of the SAA problems \eqref{composite SAA}. 
Then  
\begin{align*}
\lim_{\varepsilon\to\infty}~\sup_{n\in\NNN}\MP\big(\big\{n^{1/(4 - 2\beta)}~ \|\widehat{\theta}_{n} - \theta^{*}\|_{m} > \varepsilon\big\}\big) = 0.
\end{align*}
\end{corollary}
\begin{proof}
Let $\eta > 0$. We shall use notations \eqref{Hilfsereignis 1}, \eqref{Mab}, \eqref{Mb}.
\par
By statement 1) from Remark \ref{erste Restwahrscheinlichkeiten} we may find some $a = b > 0$ such that for any $n\in\NNN$ the inequality $\MP\big(\Omega\setminus\Omega_{n,a}\big) + 128 \overline{M}_{1} \overline{M}^{1}/[(\sqrt{2}-1)~b]\leq\eta/2$ holds. 
Then by Theorem \ref{simpliest version rough} along with statement 3) from Remark \ref{erste Restwahrscheinlichkeiten} we obtain
\begin{align*}
& 
\MP\big(\big\{n^{1/(4 - 2\beta)}~ \|\widehat{\theta}_{n} - \theta^{*}\|_{m} > 2^{K}\big\}\big)\\ 
&\leq
\frac{640~\sqrt{2}~\overline{M}_{a}[\overline{M}^{1} + 0.7] + 10\cdot\eins_{]0,\infty[}(\EEE[L_{2}(Z_{1})])\cdot M_{a,a}}{M_{1}}~2^{(K - 1)~(\beta - 3/2)} + \frac{\eta}{2}
\end{align*}
for $n, K\in\NNN$, where $M_{1}$ denotes the constant from assumption \ref{assump:A7}. Since $\beta < 3/2$, we may find some $K_{0}\in\NNN$ such that
$$
\MP\big(\big\{n^{1/(4 - 2\beta)}~ \|\widehat{\theta}_{n} - \theta^{*}\|_{m} > 2^{K_{0}}\big\}\big) < \eta\quad\mbox{for}~n\in\NNN.
$$
This completes the proof.
\end{proof}
At the end of this section we want to direct the interest to the contribution \cite{ErmolievNorkin2013}. There the authors investigate the deviation probabilities 
\begin{align*}
\MP\Big(\Big\{n^{\gamma}~d_{H}\big(\widehat{\Theta}_{n,\delta},\Theta_{\delta}\big) > \frac{\varepsilon}{\delta}\Big\}\Big)\quad(\delta > 0)
\end{align*}
instead of the deviation probabilities \eqref{deviation}. Here $d_{H}$ denotes the Hausdorff metric, and $\Theta_{\delta}, \widehat{\Theta}_{n,\delta}$ stand respectively for the 
approximate $\delta$-solution sets of the genuine problem \eqref{composite optimization} and the SAA problems \eqref{composite SAA}. They achieve to find upper estimates in terms of the sample sizes based on the starting assumptions
\begin{itemize}
\item $G$ is uniformly bounded and continuous in the parameter.
\item $H$ Borel measurable as well as continuous in $\theta$, and $H(\theta,\cdot,z)$ is H\"older continuous with some H\"older exponent $\beta$ as well as some H\"older constant, both independent of 
$(\theta,z)\in\Theta\times\RRR^{d}$.
\end{itemize}
The main difference to our paper is that our condition \ref{assump:A5} is replaced with direct requirements for the empirical Rademacher complexities of $G$ and $H$, namely 
\begin{align*}
&R_{n}(G|Z_{1},\ldots,Z_{n}) := \EEE_{\overline{e}}\Big[\sup_{\theta\in\Theta}\Big|\frac{1}{n}\sum_{j=1}^{n}\overline{e}_{j}~G(\theta,Z_{j})\Big|\Big]\leq N_{G}~n^{-\alpha_{G}}~\MP-\mbox{a.s.},\\
&R_{n}(H|Z_{1},\ldots,Z_{n}) := \EEE_{\overline{e}}\Big[\sup_{(\theta,t)\in\Theta\times\RRR}\Big|\frac{1}{n}\sum_{j=1}^{n}\overline{e}_{j}~H(\theta,t,Z_{j})\Big|\Big]\leq N_{H}~n^{-\alpha_{H}}~\MP-\mbox{a.s.}
\end{align*}
for some $\alpha_{G}, \alpha_{H}\in ]0,1/2]$ as well $N_{G}, N_{H} > 0$.
Here $\overline{e} := (\overline{e}_{1},\ldots\overline{e}_{n})$ denotes some $n$-dimensional random vector consisting of independent random variables which are uniformly distributed on $\{-1,1\}$. In view of the symmetrization method from empirical process 
theory (see e.g. \cite[Theorem 3.1.21]{GineNickl2016}, \cite[Lemma 2.3.6]{vanderVaartWellner1996}), condition \ref{assump:A5} together with sepability assumption \ref{assump:A3} relax the first requirement to
$$
\EEE\big[R_{n}(G|Z_{1},\ldots,Z_{n})\big] \leq 16~\sqrt{2}~\|\xi^{G}\|_{\MP^{Z},2}~\overline{M}^{1}/\sqrt{n}
$$
(see also \cite[Remark 3.5.5]{GineNickl2016}). The second requirement cannot be compared directly with assumptions on the functions classes $\FFF^{\Theta}_{\delta}$ in \ref{assump:A5}. These different assumptions are both tailored to different methods employed to obtain the desired deviation probabilities. The line of reasoning in \cite{ErmolievNorkin2013} does not rely on estimations for the deviation probabilities of local increments of the compound empirical process which we pointed out to be the starting point due to Proposition \ref{variational inequalities}. Instead the authors base their argumentation on the 
deviation probabilities
\begin{align*}
\MP\big(\big\{\sup_{\theta\in\Theta}\big|\GGGH(\theta,\cdot)\big| > \varepsilon~\big\}\big)
\end{align*}
(see \cite[proof of Theorem 4.1]{ErmolievNorkin2013}).
\par
The results from \cite{ErmolievNorkin2013} may not be utilized directly for the deviation probabilities \eqref{deviation}. Moreover, in principle by a suitable adaption of the proof for Proposition \ref{variational inequalities} we may also start with deviation probabilities for the compound empirical process $\GGGH$ to  derive deviation probabilities \eqref{deviation} (see e.g. \cite[proof of Lemma 1]{Pflug1999}). However, it was pointed out in \cite{Pflug1999} that nonasymptotic rates may be improved if the deviation probabilities for the increments of the compound empirical process are taken into account.
\section{Nonasymptotic confidence regions}
\label{confidence regions}
In this section we want to illustrate how to construct nonasymptotic confidence regions from our main results on the deviation probabilities \eqref{deviation}. Recalling, assumptions and notations from Theorem \ref{simpliest version} and Remarks \ref{erste Restwahrscheinlichkeiten} - \ref{unangenehmer Rest}, we fix in the first step $t > 0$ and $a > \big(\|\xi^{G}\|_{\MP^{Z},2}\vee \overline{M}_{1}\big)~\big[1 + 64\sqrt{2} (t + 1)\overline{M}^{1}\big]$ such that 
\begin{align*}
\overline{P}_{t,a} := \exp\Big(-\frac{ a\cdot\mathfrak{g}(t)}{\|\xi^{G}\|_{\MP^{Z},2}}\Big) + \frac{2~\overline{M}_{1}}{\ln(2)~a~\mathfrak{g}(t)}
~\exp\left(-\frac{a~\mathfrak{g}(t)}{\overline{M}_{1}}\right)
\end{align*}
is small. In view of statements 4) and 6) in Remark \ref{erste Restwahrscheinlichkeiten}, the value $\overline{P}_{t,a}$ partly controls the probabilities 
$$
\sum\limits_{k= K_{\varepsilon} + 1\atop \delta_{n (k-1)}\leq\Delta(\Theta)}^{\infty}\hspace*{-0.35cm}\MP\big(\Omega\setminus\Omega^{n}_{\delta_{nk}\wedge\Delta(\Theta),\sqrt{2}^{k} a}\big) + \MP\big(\Omega\setminus\Omega_{n,a}\big)
$$
in the result of Theorem \ref{simpliest version} if $a = b$ and $\varepsilon \geq 2$. Note that in case of $a = b$ we have to choose $\varepsilon \geq 2$ to fulfill the required relation between $a,b$ in the display of Theorem \ref{simpliest version}. In the second step we consider any sample size $n\geq n(a,a,\beta)\vee (\|\xi^{G}\|_{\MP^{Z},2}/2)$. We want to apply Theorem \ref{simpliest version} to lower estimate the probabilities
$$
\MP\big(\big\{\|\widehat{\theta}_{n} - \theta^{*}\|_{m}\leq\mathfrak{r}\big\}\big)
$$
for suitable $\mathfrak{r} > 0$ with fixed $a = b$ and $t$. In order to avoid trivialities the positive numbers $\mathfrak{r}$ should be smaller than the diameter $\Delta(\Theta)$ of $\Theta$. We restrict ourselves to $\varepsilon := \mathfrak{r}~n^{1/(4 - 2\beta)}\geq 2$. Such an $\varepsilon$ meets then all requirements of Theorem \ref{simpliest version} if the inequality
\begin{align*}
\mathfrak{r} \geq 2~\left(M^{*}_{t,a}\right)^{2/(3 - 2\beta)}~n^{-1/(4- 2\beta)}
\end{align*}
holds, where
\begin{align*}
M^{*}_{t,a} := \dfrac{4~M_{a,a}~\eins_{]0,\infty[}\big(\EEE[L_{2}(Z_{1})]\big) + 4~\overline{M}_{a}~[1 + 2~(t + 1)~\widehat{\eta})]}{M_{1}}.
\end{align*}
In particular, we have the additional constraint
$$
2~\left(1\vee M^{*}_{t,a}\right)^{2/(3 - 2\beta)}~n^{-1/(4- 2\beta)}\leq \mathfrak{r} < \Delta(\Theta).
$$
Hence we furthermore impose on the sample size the additional condition 
$$
n > \left(\dfrac{2}{\Delta(\Theta)}\right)^{4 - 2\beta}\cdot\left(1\vee M^{*}_{t,a}\right)^{(8 - 4\beta)/(3 - 2\beta)}.
$$
Then we may apply Theorem \ref{simpliest version} to $\varepsilon = \mathfrak{r}~n^{1/(4 - 2\beta)}$ with $\mathfrak{r}\geq 2~(1\vee M_{t,a}^{*})^{2/(3- 2\beta)}~n^{-1/(4-2\beta)}$. We have $\varepsilon\geq 2$ and
$$
\dfrac{2^{4-\beta}~\overline{M}_{a}}{(3/2 - \beta)~M_{1}~\ln(2)~\varepsilon^{3/2 - \beta}~\sqrt{2}}\leq 
\dfrac{2^{4-\beta}~\overline{M}_{a}}{(3/2 - \beta)~M_{1}\ln(2)~2^{2 - \beta}~M_{t,a}^{*}}\leq 
\dfrac{2}{(3 - 2\beta)~\ln(2)}.
$$ 
Hence the value $P^{t,\varepsilon}_{1}$ in the result of Theorem \ref{simpliest version} may be bounded by 
\begin{align*}
P^{t,\varepsilon}_{1}
&\leq 
\exp\left(\frac{\eins_{]0,\infty[}\big(\EEE[L_{2}(Z_{1})]\big)~M_{a,a}~\mathfrak{g}(t)}{a~\|L_{2}\|_{\MP^{Z},2}}\right)\cdot \frac{\exp\left(\dfrac{-\mathfrak{r}^{3/2- \beta}~M_{1}~\mathfrak{g}(t)~n^{(3/2 - \beta)/(4 - 2\beta)}}{2^{7/2-\beta}~\overline{M}_{a}}\right)}{(3/2-\beta)\ln(2)}\\
&=: \widehat{P}^{t,a}_{n}(\mathfrak{r}).
\end{align*}
Combining this with statements 4), 6) in Remark \ref{erste Restwahrscheinlichkeiten}, and Remark \ref{unangenehmer Rest}, we obtain by Theorem \ref{simpliest version} 
\begin{align*}
\MP\big(\big\{\|\widehat{\theta}_{n} - \theta^{*}\|_{m}\leq\mathfrak{r}\big\}\big)\geq 
1 - \overline{P}_{t,a} - \widehat{P}^{t,a}_{n}(\mathfrak{r}) - \overline{P}^{\varepsilon}_{n 3} - \widehat{P}^{\varepsilon}_{n 3} - \MP\big(\Omega\setminus B_{n}^{\xi^{G}}\big)
\end{align*}
for $2~(1\vee M_{t,a}^{*})^{2/(3- 2\beta)}~n^{-1/(4-2\beta)}\leq \mathfrak{r} < \Delta(\Theta)$, where 
\begin{align*}
&\overline{P}^{\varepsilon}_{n 3} := \sum_{k= K_{\varepsilon} + 1\atop \delta_{n (k-1)}\leq\Delta(\Theta)}^{\infty}\hspace*{-0.35cm}\MP\big(\Omega\setminus B_{n}^{\xi^{G}_{\delta_{nk}\wedge\Delta(\Theta)}}\big)\\
&
\widehat{P}^{\varepsilon}_{n 3} := \hspace*{-0.25cm}\sum_{k= K_{\varepsilon} + 1\atop \delta_{n (k-1)}\leq\Delta(\Theta)}^{\infty}\hspace*{-0.35cm}\MP\big(\Omega\setminus B_{n}^{(L_{1}\vee 1)~\xi^{G}_{\delta_{nk}\wedge\Delta(\Theta)}}\big) 
+ \frac{\ln\big(4 \Delta(\Theta)/\frak{r}\big)}{\ln(2)}~\MP\big(\Omega\setminus B_{n}^{L_{2}}\big).
\end{align*}
Depending on the tail behaviour of the random variables $\xi^{G}(Z_{1})^{2}, \xi^{G}_{\delta}(Z_{1})^{2}$, $L_{2}(Z_{1})$ and $(L_{1}(Z_{1})^{2}\vee1)~\xi^{G}_{\delta}(Z_{1})^{2}$ we may further estimate the summands $\overline{P}^{\varepsilon}_{n 3}, \widehat{P}^{\varepsilon}_{n 3}$ and $\MP\big(\Omega\setminus B_{n}^{\xi^{G}}\big)$ using devices from Remark \ref{Simplifications}. Let us illustrate the sketched construction of confidence regions by Example \ref{simple example}.
\begin{example}
\label{simple example cont.}
We want to continue Example \ref{simple example}.
\begin{itemize}
\item [1)] Since $Z$ is by assumption centered and normally distributed with positive definite covariance matrix $\Sigma$, the random variable $Z^{T}\Sigma^{-1}Z^{T}$ is chi-squared distributed with $m$ degrees of freedom. Then we obtain for the envelopes 
$\xi^{G}$ and $\xi^{G}_{\delta}$
as introduced in Example \ref{simple example}
\begin{align*}
&
\|\xi^{G}\|_{\MP^{Z},2} = \max_{\theta\in\Theta}\|\theta\|_{m}^{2}~\lambda^{\Sigma}~\sqrt{m (m + 2)},\\
&
\|\xi^{G}_{\delta}\|_{\MP^{Z},2} = 2~\max_{\theta\in\Theta}\|\theta\|_{m}~\lambda^{\Sigma}~\sqrt{m (m + 2)}~\delta
\quad\mbox{for}~\delta\in ]0,\Delta(\Theta)]
\end{align*}
with $\lambda^{\Sigma}$ denoting the maximal eigenvalue of $\Sigma$. Furthermore, it is well known that $\EEE[\exp(t~Z_{1}^{T}\Sigma^{-1} Z_{1})]$ is finite for some $t\in\RRR$ if and only if $t < 2$, and in this case it has $(1 - 2 t)^{-m/2}$ as its value. Then
$$
\inf\big\{c > 0\mid \EEE[\exp(\lambda~Z_{1}^{T}\Sigma^{-1} Z_{1})/\sqrt{c}]\big\} = 4~\lambda^{2}~2^{4/m}/(2^{2/m} -1)^{2}\quad\mbox{for}~\lambda > 0.
$$ 
In particular, using for $\alpha = 0.5$ the notations $|(\xi^{G})^{2}|_{\alpha}$ and $|(\xi^{G}_{\delta})^{2}|_{\alpha}$ as in statement 3) of Remark \ref{Simplifications}
\begin{align*}
&|(\xi^{G})^{2}|_{0.5} 
= 4~\max_{\theta\in\Theta}\|\theta\|_{m}^{4}~(\lambda^{\Sigma})^{2}~2^{4/m}/[2^{2/m} - 1]^{2}\\
& 
|(\xi^{G}_{\delta})^{2}|_{0.5} 
= 16~\max_{\theta\in\Theta}\|\theta\|_{m}^{2}~(\lambda^{\Sigma})^{2}~\delta^{2}2^{4/m}/[2^{2/m} - 1]^{2}\quad\mbox{for}~\delta\in ]0,\Delta(\Theta)].
\end{align*}
\item [2)] Statement 1) allows us to apply statement 3) from Remark \ref{Simplifications} for $\alpha = 0.5$. Then with constant $c(0.5)$ from there
\begin{align*}
\MP\big(\Omega\setminus B_{n}^{\xi^{G}}\big)\leq 
2~\exp\left(- \frac{\ln(2)}{8~e^{2}}~\left[\big(\varphi(m)^{2}~n\big)\wedge \big(\sqrt{\varphi(m)}~\sqrt{n}\big)\right]\right)
\end{align*}
for $n\in\NNN$, where
\begin{align*}
\varphi(m) := \dfrac{m (m + 2)~(2^{2/m}-1)^{2}}{4^{1 + 2/m}~c(0.5)}
\end{align*}
\item [3)]
Since the mappings $L_{1}, L_{2}$ are constants we have $\Omega\setminus B_{n}^{(L_{1}\vee 1)~\xi^{G}_{\delta}} = \Omega\setminus B_{n}^{\xi^{G}_{\delta}}$ for $\delta\in ]0,\Delta(\Theta)]$, and $\Omega\setminus B_{n}^{L_{2}} =\emptyset$. In particular  $\overline{P}^{\varepsilon}_{n 3} = \widehat{P}^{\varepsilon}_{n 3}$ for every $n\in\NNN$. Moreover, by statement 1) along with statement 3) in Remark \ref{Simplifications} again
\begin{align*}
\MP\big(\Omega\setminus B_{n}^{\xi^{G}_{\delta}}\big)\leq 
2~\exp\left(- \frac{\ln(2)}{8~e^{2}}~\left[\big(\varphi(m)^{2}~n\big)\wedge \big(\sqrt{\varphi(m)}~\sqrt{n}\big)\right]\right)
\end{align*}
for $n\in\NNN$ and $\delta\in ]0,\Delta(\Theta)]$. According to statement 4) in Example \ref{simple example} the condition \ref{assump:A5} is fulfilled with $\beta = 1$. Then for $\mathfrak{r} > 0$, $\varepsilon := \mathfrak{r}~\sqrt{n}$ and $n\in\NNN$
\begin{align*}
&
\overline{P}^{\varepsilon}_{n 3} + \widehat{P}^{\varepsilon}_{n 3}\\ 
&\leq 
4~\dfrac{2~\ln\big(\frac{4~\Delta(\Theta)}{\mathfrak{r}~\sqrt{n}}\big) + \ln(n)}{2~\ln(2)}~
\exp\left(- \frac{\ln(2)}{8~e^{2}}~\left[\big(\varphi(m)^{2}~n\big)\wedge \big(\sqrt{\varphi(m)}~\sqrt{n}\big)\right]\right)\\[0.2cm]
&\leq 
4~\dfrac{\ln\big(\frac{4~\Delta(\Theta)}{\mathfrak{r}}\big) }{\ln(2)}~
\exp\left(- \frac{\ln(2)}{8~e^{2}}~\left[\big(\varphi(m)^{2}~n\big)\wedge \big(\sqrt{\varphi(m)}~\sqrt{n}\big)\right]\right).
\end{align*}
\item [4)] In order to calculate the constants $\overline{P}_{t,a}$ and $\widehat{P}^{t,a}_{n}(\mathfrak{r})$ we may invoke the value of $\|\xi^{G}\|_{\MP^{Z},2}$ from statement 1) as well as the constants $M_{1}, \overline{M}_{1}, \overline{M}^{1} > 0$ and 
$K_{\delta_{1}}\geq 0$ for $\delta_{1} > 0$ as specified in Example \ref{simple example}. Note also that 
$\|L_{2}\|_{\MP^{Z},2} = 2$ holds.
\end{itemize}
\end{example} 


\section{Examples of specific classes of objectives}
\label{specific classes}
Throughout this section we shall point out some classes of objectives $G$ which meet the requirements of the main results, namely Theorems \ref{simpliest version rough}, \ref{simpliest version}. Firstly, we consider the finiteness of uniform entropy integrals, as required implicitely in Theorems \ref{simpliest version rough}, \ref{simpliest version}. It is always guaranteed if the involved covering numbers have polynomial rates. Indeed this relies on the observation, that by using change of variable formula several times along with integration by parts, we obtain
\begin{equation}
\label{Integralabschaetzung}
\int_{0}^{1}\sqrt{v\ln(K/\varepsilon)}~d\varepsilon\leq 2\sqrt{v \ln(K)}\quad\mbox{for}~v\geq 1, K\geq e.
\end{equation}
This observation is the starting point to provide explicit upper estimates of the uniform entropy integrals involved in the result of Theorems, \ref{simpliest version rough}, \ref{simpliest version} if the objective $G$ satisfies specific analytical properties.
\medskip

We start with the following condition
\begin{enumerate}
\item[(H)] There exist some $\beta\in ]0,1]$ and a square $\MP^{Z}$-integrable strictly positive mappings $C:\RRR^{d}\rightarrow ]0,\infty[$ such that 
$$
\big|G(\theta,z) - G(\vartheta,z)\big|\leq C(z)~\|\theta - \vartheta\|_{m}^{\beta}\quad\mbox{for}~z\in\RRR^{d}, \theta, \vartheta\in\Theta.
$$
\end{enumerate}
Under (H) we may construct on the one hand explicitly square $\MP^{Z}$-integrable envelopes $\xi^{G}$ and $\xi^{G}_{\delta}$ of $\FFF^{\Theta}$ and $\FFF^{\Theta}_{\delta}$ respectively. On the other hand we may also provide explicit upper estimates for the associated uniform entropy integrals $J(\FFF^{\Theta},\xi^{G},\varepsilon)$ and $J(\FFF^{\Theta}_{\delta},\xi^{G}_{\delta},\varepsilon)$. The results are gathered in the following proposition.
\begin{proposition}
\label{Hoelder-Bedingung auxiliary}
Let condition (H) be fulfilled with $\beta\in ]0,1]$ and square $\MP^{Z}$-integrable strictly positive mapping $C$. If $G(\theta,\cdot)$ be Borel measurable for every $\theta \in\Theta$, and if $G(\overline{\theta},\cdot)$ is square $\MP^{Z}$-integrable for some $\overline{\theta}\in\Theta$, then the following statements are valid
\begin{itemize}
\item [1)] The mapping $\xi^{G} := C(\cdot)~\Delta(\Theta)^{\beta} + |G(\overline{\theta},\cdot)|$ is a positive square $\MP^{Z}$-integrable envelope of $\FFF^{\Theta}$ with 
$$
J(\FFF^{\Theta},\xi^{G},\varepsilon)\leq 2\varepsilon\sqrt{(3m + 1)\ln(2) + \frac{m}{\beta}\ln(2/\varepsilon)}\quad\mbox{for}~\varepsilon \in ]0,1/2].
$$
In particular by dominated convergence theorem separability condition \ref{assump:A3} is fulfilled for 
every at most countable dense subset $\cC(\cU_{\delta})$ of $\cU_{\delta}$ $(\delta > 0)$.
\item [2)] For $\delta > 0$, the mapping $\xi^{G}_{\delta} := C(\cdot)~2^{\beta}~\delta^{\beta}$ defines a positive square $\MP^{Z}$-integrable envelope of $\FFF^{\Theta}_{\delta}$ satisfying
$$
J(\FFF^{\Theta}_{\delta},\xi^{G}_{\delta},\varepsilon)\leq 2\varepsilon\sqrt{(3m + 1)\ln(2) + \frac{m}{\beta}\ln(2/\varepsilon)}\quad\mbox{for}~\varepsilon \in ]0,1/2].
$$
\end{itemize}
\end{proposition}
\begin{proof}
Note that $|[G(\theta,z) - G(\theta^{*},z)] - [G(\vartheta,z) - G(\theta^{*},z)]|\leq C(z)~\|\theta - \vartheta\|_{m}^{\beta}$ holds for $\theta,\vartheta\in \cU_{\delta}$, and that $0\in\FFF^{\Theta}_{\delta}$ is square $\MP^{Z}$-integrable for $\delta > 0$. Then the statements of Proposition \ref{Hoelder-Bedingung auxiliary} are immediate consequences of separate applications of Proposition 2.6 from \cite{Kraetschmer2023a} to the function classes 
$\FFF^{\Theta}$ and $\FFF^{\Theta}_{\delta}$ ($\delta > 0$). In particular \ref{assump:A3} may be verified for any at most countable dense subset $\cC(\cU_{\delta})$ of $\cU_{\delta}$ $(\delta > 0)$ by the dominated convergence theorem.
\end{proof}
\medskip
\begin{remark}
\label{Remark Hoelder-Bedingung}
Proposition \ref{Hoelder-Bedingung auxiliary} tells us that under (H) the results of Theorems \ref{simpliest version rough}, \ref{simpliest version} may be concluded from assumptions \ref{assump:A1}, \ref{assump:A4} and \ref{assump:A6} together with the conditions
\begin{itemize}
\item $G(\overline{\theta},\cdot)$ is square $\MP^{Z}$-integrable for some $\overline{\theta}\in\Theta$.
\item $(L_{1}\vee 1)~C(\cdot)$ is square $\MP^{Z}$-integrable.
\end{itemize}
In particular we may choose in \ref{assump:A5} the 
constants $\overline{M}_{1} := 2^{\beta} \|(L_{1}\vee 1)~C(\cdot)\|_{\MP^{Z},2}$ and $\overline{M}^{1} := \sqrt{(3m + 1)\ln(2) + m\ln(16)/\beta}/4$. 
\end{remark}
Next, let us consider objective $G$ having the following kind of structure of piecewise H\"older continuity. 
\begin{enumerate}
\item [(PH)] $
G(\theta,z) = \sum\limits_{i = 1}^{r}\eins_{\bigcap_{l=1}^{s_{i}}\{\Lambda_{il}(\theta, \cdot) + a^{i}_{l}\in I_{il}\}}(z)\cdot G^{i}(\theta,z),
$ 
where
\begin{itemize}
\item $r, s_{1},\dots,s_{r}\in\NNN$,
\item $G^{i}$ satisfies \ref{assump:A1}, and (H) with $\beta_{i}\in ]0,1]$ as well as strictly positive square $\MP^{Z}$-integrable $C^{i}:\RRR^{d}\rightarrow\RRR$ for $i\in\{1,\ldots,r\}$,
\item 
$\Lambda_{il}:\RRR^{m}\times\RRR^{d}\rightarrow\RRR$ Borel measurable with $\Lambda_{il}(\cdot,z)$ affine linear for $z\in\RRR^{d}$ ($i\in\{1,\ldots,r\}$, $l\in\{1,\ldots,s_i\}$),
\item $a^{i}_{l}\in\RRR$ for $i\in\{1,\dots,r\}, l\in\{1,\dots,s_{i}\}$,
\item $I_{il} = ]0,\infty[$ or $I_{il} = [0,\infty[$ for $i\in\{1,\dots,r\}$ and $l\in\{1,\dots,s_{i}\}$,
\item The set
$$
\Big\{\bigcap\limits_{l=1}^{s_{i}}\big\{\Lambda_{il}(\theta, \cdot) + a^{i}_{l}\in I_{il}\}\mid i\in\{1,\ldots,r\}, l\in\{1,\ldots,s_{i}\big\}\Big\}
$$
is a partition of $\RRR^{d}$.
\end{itemize}
\end{enumerate}
In two stage mixed-integer programs the goal functions typically may be represented in this way if the random vector $Z$ has compact support (see \cite[p. 121]{EichhornRoemisch2007} along with \cite{Kraetschmer2023a}).
\par
Note that if $G$ satisfies condition (PH), it does not have continuity in $\theta$ in advance.
\par
For abbreviation we set 
$B_{i}(\theta) := \bigcap_{l=1}^{s_{i}}\{z\in\RRR^{d}\mid\Lambda_{il}(\theta, z) + a^{i}_{l}\in I_{il}\}$
for arbitrary $i$ from $\{1,\ldots,r\}$, and we introduce the associated function classes
$$
\FFF_{\textrm{\tiny PH}}^{i} := \big\{\eins_{B_{i}(\theta)}\mid \theta\in\Theta\big\}\quad\mbox{and}\quad\overline{\FFF}_{\textrm{\tiny PH}}^{i} := \big\{G^{i}(\theta,\cdot)\mid \theta\in\Theta\big\}\quad i\in\{1,\ldots,r\}.
$$
Note that the classes $\FFF_{\textrm{\tiny PH}}^{i}$ are uniformly bounded by $1$. We borrow from \cite{Kraetschmer2023a} (Proposition 2.8 there) the following result concerning construction of envelope $\xi^{G}$ of $\FFF^{\Theta}$ and upper estimation of the associated uniform entropy integrals $J(\FFF^{\Theta},\xi^{G},\varepsilon)$.
\begin{proposition}
\label{startingpoint}
The set $B_{i}(\theta)$ is a Borel subset of $\RRR^{d}$, and the mapping $G^{i}(\theta,\cdot)$ is Borel measurable for $\theta\in\Theta$ and $i\in\{1,\ldots,r\}$. In particular $G(\theta,\cdot)$ is Borel measurable for every $\theta\in\Theta$. Moreover, if $G^{1}(\overline{\theta},\cdot),\ldots,G^{r}(\overline{\theta},\cdot)$ are square $\MP^{Z}$-integrable for some $\overline{\theta}\in\Theta$, and if $\xi_{1},\ldots,\xi_{r}$ denote bounded positive envelopes of the classes $\FFF_{\textrm{\tiny PH}}^{1},\ldots,\FFF_{\textrm{\tiny PH}}^{r}$ respectively, then the mapping $\xi^{G} := \sum_{i=1}^{r}\xi_{i}\cdot\big(\Delta(\Theta)^{\beta_{i}}~C_{i}(\cdot) + |G^{i}(\overline{\theta},\cdot)|\big)$ is a positive square $\MP^{Z}$-integrable envelope of $\FFF^{\Theta}$ satisfying  
\begin{align*}
&
J(\FFF^{\Theta},\xi^{G},\varepsilon)
\\
&
\leq 2\varepsilon\sqrt{r + 2 r~m~\ln(3) +  m \ln(4 r/\varepsilon) \sum_{i=1}^{r} 1/\beta_{i} + \ln(2) + [5 + 2 \ln(4 r/\varepsilon)]~ (m + 2) \sum_{i=1}^{r}s_{i}}
\end{align*}
for $\varepsilon\in ]0,1]$.
\end{proposition}
Let us turn over to consider the function classes $\FFF^{\Theta}_{\delta}$ under representation (PH). We set $A\Delta B := A\setminus B~\cup~B\setminus A$ for sets $A, B$, and we shall use the following notation
\begin{equation}
\label{alle Differenzmengen}
\overline{B}_{i\delta} = \bigcup_{\theta\in\cU_{\delta}}B_{i}(\theta)\Delta B_{i}(\theta^{*})\quad\mbox{for}~\delta > 0, i\in\{1,\ldots,r\}.
\end{equation}
\begin{proposition}
\label{important covering numbers}
Let $G$ fulfill representation (PH) with mappings $G^{1}(\theta,\cdot),\ldots,G^{r}(\theta,\cdot)$ being Borel measurable for $\theta\in\Theta$. Furthermore, let $C^{1},\ldots C^{r}:\RRR^{d}\rightarrow\RRR$ as well as $\beta_{1},\ldots,\beta_{r}\in ]0,1]$ be from representation (PH), and let for $\delta > 0$ denote by $\overline{B}_{1 \delta},\ldots,\overline{B}_{r \delta}$ the sets as in \eqref{alle Differenzmengen}. Then, for $\delta > 0$, fixing Borel subsets $\widehat{B}_{1\delta},\ldots \widehat{B}_{r\delta}$ of $\RRR^{d}$ with 
$\widehat{B}_{i\delta}\supseteq \overline{B}_{i\delta}$ $(i\in\{1,\ldots,r\})$, a positive envelope $\xi^{G}_{\delta}$ of the function class $\FFF^{\Theta}_{\delta}$ is 
defined by $\xi^{G}_{\delta}(z) := \sum_{i=1}^{r}[2~\delta^{\beta_{i}}~C^{i}(z) + |G(\theta^{*},z)|~\eins_{\widehat{B}^{\delta}_{i}} + (\delta\wedge 1)^{2}]$, 
and 
\begin{align*}
J(\FFF^{\Theta}_{\delta},\xi^{G}_{\delta},\varepsilon)
&\leq 
2\varepsilon \sqrt{r+\ln(2) + c_{mr\varepsilon} + \big[5 + 2 \ln\big(8 (r + 1)/\varepsilon\big)\big]~ d_{mr}}\\
&\quad 
+ 
2\varepsilon \sqrt{r  + 4\big[1 + \ln\big(8 (r + 1)/\varepsilon\big)\big]~ d_{mr} }
\end{align*}
holds for every $\varepsilon\in]0,1[$, where $c_{mr\varepsilon} := 2 r~m~\ln(3) +  m \ln(8 (r + 1)/\varepsilon) \sum_{i=1}^{r} 1/\beta_{i}$ and $d_{mr} := (m + 2) \sum_{i=1}^{r} s_{i}$.
\end{proposition}
The proof of Proposition \ref{important covering numbers} is subject of Subsection \ref{Beweis von Proposition important covering numbers}.
\begin{remark}
\label{Remark PH}
Let us discuss how to apply Theorem \ref{simpliest version rough}, \ref{simpliest version} if $G$ has representation (PH), and the sets $B_{i}(\theta)$ satisfy the following property
\begin{itemize}
\item [(*)] The set $\left\{z\in\RRR^{d}\mid \Lambda_{il}(\theta,z) = - a^{i}_{l}~\mbox{for some}~ \theta\in\Theta\right\}$ is a $\MP^{Z}$-null set for all indices $i \in \{1,\ldots,r\}$ and $l\in\{1,\ldots,s_{i}\}$ with $I_{il} = [0,\infty[$. 
\end{itemize}
 \begin{itemize}
\item [1)] By Proposition \ref{startingpoint} we know how to find a square $\MP^{Z}$-integrable positive envelope of $\FFF^{\Theta}$ for \ref{assump:A2}. Furthermore under condition (*) the continuity property required in \ref{assump:A1} and also condition \ref{assump:A3} may be verified by routine procedures. 
\item [2)] Since $\Theta$ is separable, property (*) implies that any set $\overline{B}_{i\delta}$ may be described up to some $\MP^{Z}$-null set as an at most countable union of Borel subsets of $\RRR^{d}$. Hence for $\delta > 0$, $i\in\{1,\ldots,r\}$, the set $\overline{B}_{i\delta}$ is $\MP^{Z}$-measurable, and thus $\eins_{\overline{B}_{i\delta}}(Z_{1})$ is a random variable on the complete probability space $\OFP$.
\item [3)] Under \ref{assump:A4} with square $\MP^{Z}$-integrable mapping $L_{1}$ we may replace \ref{assump:A5} with the following condition:
\begin{itemize}
\item [(**)] The mappings $(L_{1}\vee 1)~C^{1},\ldots,(L_{1}\vee 1)~C^{r}$ are square $\MP^{Z}$-integrable, and there exist $\beta\in ]0,1]$ and $\widehat{M}_{1} > 0$  such that for every $i\in\{1,\ldots,r\}$
\begin{align*}
\beta\leq\beta_{i}\quad\mbox{and}\quad\sup_{\delta\in ]0,\Delta(\Theta)]}~\EEE\big[\big(L_{1}(Z_{1})\vee 1\big)^{2}~G^{i}(\theta^{*},Z_{1})^{2}~\eins_{\overline{B}_{i\delta}}(Z_{1})\big]/\delta^{2\beta}\leq\widehat{M}_{1}.
\end{align*}
\end{itemize}
In view of 2) we may find for any $\overline{B}_{i\delta}$ some Borel subset $\widehat{B}_{i\delta}$ of $\RRR^{d}$ enclosing $\overline{B}_{i\delta}$ such that $\eins_{\widehat{B}_{i\delta}}(Z_{1})= \eins_{\overline{B}_{i\delta}}(Z_{1})$ $\MP$-a.s.. Using the envelopes from Propositions \ref{startingpoint}, \ref{important covering numbers} with the sets $\widehat{B}_{i\delta}$, assumption (**) implies \ref{assump:A5} with $\beta$ and 
\begin{align*}
&\overline{M}_{1} := r~\big(\widehat{M}_{1}^{1/2} + [\Delta(\Theta)\wedge 1]^{2-\beta}~\|L_{1}\vee 1\|_{\MP^{Z},2}\big)\\ 
&\qquad\qquad+ 2~\sum_{i=1}^{r}\|(L_{1}\vee 1)~C^{i}\|_{\MP^{Z},2}~[\Delta(\Theta)^{\beta_{i}-\beta}\vee 1]\\
& 
\overline{M}^{1} 
:= 
\sqrt{r + \overline{c}_{mr} + \ln(2) + [5 + 2 \ln(\overline{r})]~ d_{mr}}/4 
+ 
\sqrt{r + 4 [1 + \ln(\overline{r})]~ d_{mr} }/4, 
\end{align*}
where $\overline{r} := 64 (r+1)$, and $\overline{c}_{mr} := 2 r~m~\ln(3) +  m \ln(\overline{r}) \sum_{i=1}^{r} 1/\beta_{i}$ as well as $d_{mr} := (m + 2) \sum_{i=1}^{r} s_{i}$.
\end{itemize}
\end{remark}
\section{Illustration of some assumptions}
\label{illustration assumptions}
In this section we shall comment on the second order growth condition \ref{assump:A7} and assumption \ref{assump:A6}. We start with the second order growth condition.
\par
In the first view this condition seems to be more restrictive than the more usual 
\textit{local second order growth condition}.
\begin{enumerate}[label=(A \arabic*'), ref=(A \arabic*')]
\setcounter{enumi}{4}
\item \label{assump:A7'}
For some $\widehat{\delta}, \widehat{M} > 0$ the goal function $\psi_{H,\Theta}$ of optimization \eqref{composite optimization} satisfies
\begin{align*}
\psi_{H,\Theta}(\theta) - \psi_{H,\Theta}(\theta^{*})\geq \widehat{M}~\|\theta - \theta^{*}\|_{m}^{2}\quad\mbox{for}~\theta\in\Theta~\mbox{with}~\|\theta - \theta^{*}\|_{m}\leq\widehat{\delta}.
\end{align*}
\end{enumerate}
Of course the second order growth condition implies the local second order growth condition. Actually, it will turn out that both conditions are even equivalent within our setting. We define for every $\delta > 0$ 
the nonnegative number 
$$
M^{H}(\delta) := \inf\{\psi_{H,\Theta}(\theta) - \psi_{H,\Theta}(\theta^{*})\mid\theta\in\Theta\setminus\cU_{\delta}\}.
$$
\begin{lemma}
\label{lokal gleich global}
Let \eqref{unique solution} and the assumptions \ref{assump:A1}, \ref{assump:A2} and \ref{assump:A4}, \ref{assump:A7'} be fulfilled.
Then $\psi_{H,\Theta}$ is a continuous mapping w.r.t. the Euclidean norm. Moreover,
if \ref{assump:A7'} holds with 
$\widehat{\delta}, \widehat{M} > 0$, then the number $M^{H}(\widehat{\delta})$ is strictly positive and 
\ref{assump:A7} is satisfied with $M_{1} = \widehat{M}\wedge [M^{H}(\widehat{\delta})/ \Delta(\Theta)^{2}]$.
\end{lemma}
\begin{proof}
First of all, $\psi_{H,\Theta}$ is a continuous mapping w.r.t. the Euclidean norm due to Remark \ref{Stetigkeiten}.
\par
Now, let us assume that \ref{assump:A7'} holds with $\widehat{\delta}, \widehat{M} > 0$. By continuity of the goal function $\psi_{H,\Theta}$ along with compactness of $\Theta$, and since $\theta^{*}$ is the unique minimizer of $\psi_{H,\Theta}$, any minimizing sequence $(\theta_{k})_{k\in\NNN}$ converges to $\theta^{*}$. This implies $M^{H}(\widehat{\delta}) > 0$.
\par
For $\theta\in\Theta\setminus\cU_{\widehat{\delta}}$ the definition of $M^{H}(\widehat{\delta})$ implies
\begin{align*}
\psi_{H,\Theta}(\theta) - \psi_{H,\Theta}(\theta^{*})\geq [M^{H}(\widehat{\delta})/\|\theta - \theta^{*}\|_{m}^{2}]~\|\theta - \theta^{*}\|_{m}^{2}\geq [M^{H}(\widehat{\delta})/\Delta(\Theta)^{2}]~\|\theta - \theta^{*}\|_{m}^{2}.
\end{align*} 
Now, we may conclude immediately from \ref{assump:A7'} the second order growth condition with 
$M_{1} = \widehat{M}\wedge [M^{H}(\widehat{\delta})/\Delta(\Theta)^{2}]$. This completes the proof. 
\end{proof}
\begin{remark}
\label{well-separated}
Let \eqref{unique solution} and \ref{assump:A1}, \ref{assump:A2} as well as \ref{assump:A4} be fulfilled. 
We already know from Lemma \ref{lokal gleich global} that $\psi_{H,\Theta}$ is a continuous mapping w.r.t. the Euclidean norm. Then, if in addition the unique solution 
$\theta^{*}$ belongs to the topological interior $int(\Theta)$ of $\Theta$ w.r.t. the standard topology of $\RRR^{m}$, the continuity of $\psi_{H\Theta}$ together with Lemma \ref{lokal gleich global} imply that the growth conditions \ref{assump:A7} and \ref{assump:A7'} are both equivalent with the following property
$$
\liminf_{t\searrow 0, \underline{y}'\to \underline{y}}\frac{\psi_{H,\Theta}(\theta^{*} + t~\underline{y}') - \psi_{H,\Theta}(\theta^{*})}{t^{2}/2} > 0\quad\mbox{for all}~\underline{y}\in\RRR^{m}\setminus\{0\}
$$
(see \cite[Proposition 3.100]{BonnansShapiro2000}). This condition in turn is valid if $\psi_{H,\Theta}$ is twice continuously differentiable at $\theta^{*}$ with positive definite Hessian matrix.
\end{remark}
Next we want to find a simplier condition to ensure assumption \ref{assump:A6}. We shall use notation $\cV_{\delta} := \big\{(\theta,t)\in\Theta\times\RRR\mid \|(\theta,t) - (\theta^{*},\EEE[G(\theta^{*},Z_{1})])\big\|_{m+1}\leq\delta\}$ for $\delta > 0$.
\begin{enumerate}[label=(A \arabic*'), ref=(A \arabic*')]
\setcounter{enumi}{7}
\item \label{assump:A6'} There exist $\overline{\delta} >0, \overline{K} \geq 0$ and a Borel measurable mapping 
$m_{\Theta\times\RRR}: \Theta\times\RRR\rightarrow\RRR$ which is bounded by $\overline{K}$, Lipschitz continuous on $\cV_{\overline{\delta}}$ with Lipschitz-constant $K_{\overline{\delta}}$, and satisfies
$$
\EEE[H(\theta,t,Z_{1})] - \EEE[H(\theta,s,Z_{1})] = \int_{s}^{t}m_{\Theta\times\RRR}(\theta,u)~du\quad\mbox{for}~(\theta,t), (\theta,s)\in\Theta\times\RRR, t > s.
$$
\end{enumerate}
Indeed, under condition \ref{assump:A5} the assumption \ref{assump:A6} may be derived from \ref{assump:A6'}. 
\begin{lemma}
\label{simplier A6}.
Let \ref{assump:A1}, \ref{assump:A5} and \ref{assump:A6'} be fulfilled with 
$\overline{M}_{1} > 0, \beta\in ]0,1]$ from condition \ref{assump:A5}, and mapping $m_{\Theta\times\RRR}$ from \ref{assump:A6'}. Furthermore let $\overline{\delta}$ and 
$\overline{K}, K_{\overline{\delta}}\geq 0$ be as in \ref{assump:A6'}. Then, setting
$\overline{\delta}_{1} := \Delta(\Theta)\wedge\overline{\delta}\wedge \big(\overline{\delta}/[2~\Delta(\Theta)^{1-\beta} + 2~\overline{M}_{1}]\big)^{1/\beta}$, the assumption \ref{assump:A6} is satisfied with $\delta_{1} =\overline{\delta}/2$ and 
$K_{\delta_{1}} := K_{\overline{\delta}}\vee (2 \overline{K}/\overline{\delta}_{1})$.
\end{lemma}
\begin{proof}
By \ref{assump:A5} we may observe for any $\theta\in\cU_{\overline{\delta_{1}}}$ and every $x\in [-\overline{\delta}/2,\overline{\delta}/2]$
\begin{align*}
\big\|\big(\theta - \theta^{*},[\overline{\psi}(\theta) - x] - \overline{\psi}(\theta^{*})  \big)\big\|_{m+1}^{2}
&\leq 
\|\theta - \theta^{*}\|_{m}^{2} + 
\big(|\overline{\psi}(\theta) - \overline{\psi}(\theta^{*}) | + |x|\big)^{2}\\
&\leq 
\overline{\delta}_{1}^{2} + \big(\overline{M}_{1}~\overline{\delta}_{1}^{\beta} + \overline{\delta}/2\big)^{2}\leq \overline{\delta}^{2}.
\end{align*} 
This means $(\theta,\overline{\psi}(\theta) - x)\in\cV_{\overline{\delta}}$ for 
$(\theta,x)\in\cU_{\overline{\delta}_{1}}\times [-\overline{\delta}/2,\overline{\delta}/2]$. Hence in view of \ref{assump:A6'}
\begin{align*}
\big|m_{\Theta\times\RRR}(\theta,\overline{\psi}(\theta) - x) - m_{\Theta\times\RRR}(\theta^{*},\overline{\psi}(\theta^{*}) - x)\big|\leq K_{\overline{\delta}}~\big\|\big(\theta - \theta^{*},\overline{\psi}(\theta)- \overline{\psi}(\theta^{*})\big)\big\|_{m+1}
\end{align*}
for $\theta\in\cU_{\overline{\delta}_{1}}$ and $x\in [-\overline{\delta}/2,\overline{\delta}/2]$. Moreover, by boundedness of 
$m_{\Theta\times\RRR}$ according to \ref{assump:A6'} we may further conclude for 
$\theta\in\Theta\setminus\cU_{\overline{\delta}_{1}}$ and $x\in\RRR$
\begin{align*}
&
\big|m_{\Theta\times\RRR}(\theta,\overline{\psi}(\theta) - x) - m_{\Theta\times\RRR}(\theta^{*},\overline{\psi}(\theta^{*}) - x)\big|\\ 
&= 
\dfrac{\big|m_{\Theta\times\RRR}(\theta,\overline{\psi}(\theta) - x) - m_{\Theta\times\RRR}(\theta^{*},\overline{\psi}(\theta^{*}) - x)\big|}{\big\|\big(\theta - \theta^{*},\overline{\psi}(\theta)- \overline{\psi}(\theta^{*})\big)\big\|_{m+1}}~\big\|\big(\theta - \theta^{*},\overline{\psi}(\theta)- \overline{\psi}(\theta^{*})\big)\big\|_{m+1}\\
&\leq
\dfrac{2~\overline{K}}{\overline{\delta}_{1}}~\big\|\big(\theta - \theta^{*},\overline{\psi}(\theta)- \overline{\psi}(\theta^{*})\big)\big\|_{m+1}.
\end{align*}
This completes the proof.
\end{proof}
\section{Error estimates for M-estimators}
\label{error estimates m-estimators}
Throughout this section we consider the classical risk neutral optimization problem
\begin{align} 
\label{optimization risk neutral} 
\min_{\theta\in\Theta}\EEE[G(\Theta,Z)]
\end{align}
and the associated empirical counterparts 
\begin{align}
\label{risk neutral SAA}
\min_{\theta\in\Theta}\frac{1}{n}\sum_{j=1}^{n}G(\theta,Z_{j})\quad(n\in\NNN).
\end{align} 
This setting may be viewed as a special case of \eqref{composite optimization} with \eqref{composite SAA}, choosing the mapping $H:\Theta\times\RRR\times\RRR^{d}\rightarrow\RRR$, defined by 
$H(\theta,t,z) = G(\theta,z)$. In this situation we have the following specializations of particular assumptions of Theorems \ref{simpliest version rough}, \ref{simpliest version}.
\begin{itemize}
\item [(7i)] The condition \ref{measurability H} holds under \ref{assump:A1} together with \ref{assump:A2}.
\item [(7ii)] Assumption \ref{assump:A4} is satisfied with $L_{1} :\equiv 1$ as well as $L_{2} :\equiv 0$,
\item [(7iii)] In \ref{assump:A6} we may use $\delta_{1} > 0$ arbitrarily together with $K_{\delta_{1}} := 0$ as well as $m_{\Theta\times\RRR} :\equiv 0$. 
\item [(7iv)] The mapping $\xi^{G}_{a,b,\delta}$ defined in \eqref{neue envelope} boils down to the mapping $\xi^{G}_{a,b,\delta} = \xi^{G}_{\delta}$ from assumption \ref{assump:A5} for $a,b > 0$ and $\delta\in ]0,\Delta(\Theta)]$. By \ref{assump:A5} again this implies
$\|\xi^{G}_{a,b,\delta}\|_{\MP^{Z},2}\leq \overline{M}_{1}~\delta^{\beta}$ with $\overline{M}_{1} > 0$ and $\beta\in ]0,1]$ from \ref{assump:A5}.
\item [(7v)] The terms, introduced in \eqref{Mab}, \eqref{Mb} and \eqref{etaabbeta}, satisfy $M_{a,b} = 0$, $\overline{M}_{b} = \overline{M}_{1}$, and
  for $a,b, t, \delta > 0$,  and
$$
n(a,b,\beta) 
=
[\overline{M}_{1}^{2}~\Delta(\Theta)^{2\beta}/2]\vee [\overline{M}_{1}^{2}~\Delta(\Theta)^{2\beta +1}]^{\frac{4-2\beta}{3 - 2\beta}} 
$$
for $a,b,\delta > 0$, $\beta\in ]0,1]$. 
\end{itemize}
We assume that \eqref{unique solution} holds, i.e. \eqref{optimization risk neutral} has a unique minimizer $\theta^{*}\in \Theta$. In addition we consider any sequence $\big(\widehat{\theta}_{n}\big)_{n\in\NNN}$ of $m$-dimensional random vectors $\widehat{\theta}_{n}$ which minimize problem \eqref{risk neutral SAA} w.r.t. the sample size $n$. Such random vectors are also known as \textit{M-estimators}. The application of Theorems \ref{simpliest version rough}, \ref{simpliest version} reads as follows.
\begin{theorem}
\label{deviation m-estimation}
Let assumptions \ref{assump:A1}, \ref{assump:A2} and \ref{assump:A7} - \ref{assump:A5} be fulfilled with 
\begin{itemize}
\item the constant $M_{1} > 0$ as in \ref{assump:A7},
\item the constant $\beta\in ]0,1]$ from \ref{assump:A5}.
\end{itemize}
Furthermore let $a,b,\varepsilon > 0$ with $a~2^{\beta + 1/2}\leq b\varepsilon^{\beta + 1/2}$, and let $\mathfrak{g}(t), \widehat{\eta}$ be the constants introduced in \eqref{frakg}, \eqref{etadachepsilon} respectively.  
Then, setting $\delta_{nk} := 2^{k}n^{-1/(4- 2\beta)}$ for $k\in\mathbb{Z}, n\in\NNN$, and using notation for the auxiliary events defined in \eqref{Hilfsereignis 1},   \eqref{Hilfsereignis 2} and \eqref{Restmenge} the following statements hold.
\begin{itemize}
\item[1)] For $n\in\NNN$ with $n\geq [\overline{M}_{1}^{2}~\Delta(\Theta)^{2\beta}/2]\vee [\overline{M}_{1}^{2}~\Delta(\Theta)^{2\beta +1}]^{\frac{4-2\beta}{3 - 2\beta}} $
\begin{align*}
&
\MP\big(\big\{n^{1/(4 - 2\beta)}~ \|\widehat{\theta}_{n} - \theta^{*}\|_{m} > \varepsilon\big\}\big)
\leq P^{\varepsilon}_{1} + P^{\varepsilon}_{n 2} + \MP\big(\Omega\setminus\Omega_{n,a}\big),
\end{align*}
where
\begin{align*}
&
P^{\varepsilon}_{1} := \frac{10~ \overline{M}_{1}~\widehat{\eta}}{M_{1}}\left(\frac{2}{\varepsilon}\right)^{3/2 - \beta} \sqrt{\frac{2}{\varepsilon\wedge 2}}
\quad\mbox{and}\quad
P^{\varepsilon}_{n 2} := \sum_{k= K_{\varepsilon} + 1\atop \delta_{n (k-1)}\leq\Delta(\Theta)}^{\infty}\hspace*{-0.35cm}\MP\big(\Omega\setminus\Omega^{n}_{\delta_{nk}\wedge\Delta(\Theta),\sqrt{2}^{k} b}\big).
\end{align*}
\item[2)]If $n\in\NNN$ with $n\geq [\overline{M}_{1}^{2}~\Delta(\Theta)^{2\beta}/2]\vee [\overline{M}_{1}^{2}~\Delta(\Theta)^{2\beta +1}]^{\frac{4-2\beta}{3 - 2\beta}} $, and if for $t > 0$ the inequality $2^{K_{\varepsilon} (2 - \beta) - K_{\varepsilon}^{+}/2} > 4 ~\overline{M}_{1}~\big[1 + 2 (t + 1)~\widehat{\eta}\big]/M_{1}$ holds, then
\begin{align*}
&
\MP\big(\big\{n^{1/(4 - 2\beta)}~ \|\widehat{\theta}_{n} - \theta^{*}\|_{m} > \varepsilon\big\}\big)\leq P^{t,\varepsilon}_{1} + P^{\varepsilon}_{n 2} + P_{n 3} + \MP\big(\Omega\setminus\Omega_{n,a}\big),
\end{align*} 
where
\begin{align*}
&P^{t,\varepsilon}_{1} :=
\frac{2^{4- \beta}~\overline{M}_{1}\cdot\exp\left(-\dfrac{\varepsilon^{3/2- \beta}\sqrt{\varepsilon\wedge 2}~M_{1}~\mathfrak{g}(t)}{2^{4-\beta}~\overline{M}_{1}}\right)}{(3/2-\beta)\ln(2)~M_{1}~\varepsilon^{3/2-\beta}\sqrt{\varepsilon\wedge 2}}\\
&P^{\varepsilon}_{n 2} := \sum_{k= K_{\varepsilon} + 1\atop \delta_{n (k-1)}\leq\Delta(\Theta)}^{\infty}\hspace*{-0.35cm}\MP\big(\Omega\setminus\Omega^{n}_{\delta_{nk}\wedge\Delta(\Theta),\sqrt{2}^{k} b}\big)\quad\mbox{and}\quad
%
P^{\varepsilon}_{n 3} :=  
\sum_{k= K_{\varepsilon} + 1\atop \delta_{n (k-1)}\leq\Delta(\Theta)}^{\infty}\hspace*{-0.35cm}\MP\big(\Omega\setminus B_{n}^{\xi^{G}_{\delta_{nk}\wedge\Delta(\Theta)}}\big).
\end{align*}
\end{itemize}
In particular the sequence $\big(n^{1/(4 - 2\beta)}~ [\widehat{\theta}_{n} - \theta^{*}]\big)_{n\in\NNN}$ is uniformly tight.
\end{theorem}
\begin{proof}
In view of (7i) - (7iii) we may apply Theorems \ref{simpliest version rough}, \ref{simpliest version} and Corollary \ref{allgemeine Konvergenzraten zweiter Teil}. Then additionally taking into account (7iv), (7v) statement 1) follows immediately from Theorem \ref{simpliest version rough}, whereas statement 2) is a direct consequence of Theorem \ref{simpliest version}. Finally, the tightness $\big(n^{1/(4 - 2\beta)}~ [\widehat{\theta}_{n} - \theta^{*}]\big)_{n\in\NNN}$ may be concluded from Corollary \ref{allgemeine Konvergenzraten zweiter Teil}. The proof is complete.
\end{proof}
In the recent study \cite{OstrovskiiBach2021} on deviation probabilities for $M$-estimators the authors consider convex $\Theta$ and objectives $G$ of the following form
$$
G(\theta,z) = \mathfrak{l}\big(\pi_{2}(z),\pi_{1}(z)^{T}\theta\big),
$$
where
\begin{itemize}
\item $\pi_{1}(x_{1},\ldots,x_{d}) := (x_{1},\ldots,x_{d-1})$ and $\pi_{2}(x_{1},\ldots,x_{d}) := x_{d}$ for $(x_{1},\ldots,x_{d})\in\RRR^{d}$;
\item $\mathfrak{l}: \RRR^{2}\rightarrow\RRR$ Borel measurable, and  $\mathfrak{l}(y,\cdot)$ continuously differentiable for $y\in\RRR$.
\end{itemize}
Hence, under obvious conditions on integrability, such objectives satisfy property (H) with $\beta = 1$ and function $C:\RRR^{d}\rightarrow\RRR$, defined by
$$
C(z) = 1\vee \big[\|\pi_{1}(z)\|_{m}\cdot\sup_{\theta\in\Theta}\big|\mathfrak{l}\big(\pi_{2}(z),\cdot\big)'\big(\pi_{1}(z)^{T}\theta\big)\big|\big].
$$
Then, invoking Remark \ref{Remark Hoelder-Bedingung}, we obtain immediately upper bounds for the deviation probabilities of $M$-estimators from Theorem \ref{deviation m-estimation} along with Remarks \ref{erste Restwahrscheinlichkeiten} - \ref{unangenehmer Rest}.
\par
The authors in \cite{OstrovskiiBach2021} follow a different line of reasoning, additionally 
assuming that for any $y\in\RRR$ the mapping $\mathfrak{l}(y,\cdot)$ is three times differentiable with bounds of the third derivatives in terms the second derivatives, a condition also known as \textit{self-concordance}. They also require the random vector $\pi_{1}(Z)$ to have finite subgaussian norm, defined to mean that there exists some $\lambda > 0$ such that
$$
\sup_{x\in\RRR^{d-1}}\EEE\Big[\exp\Big(\lambda~\big(x^{T}\pi_{1}(Z)\big)^{2}\Big)\Big]\leq 2.
$$
In this special situation the authors achieved to boil down the investigations of the deviation probabilities for $M$-estimators to concentration inequalities for quadratic forms of subgaussian random vectors. As a special feature, their result, namely Theorem 3.1, provides bounds which do not depend on the dimension of the parameter set $\Theta$. In contrast the estimates derived from Theorem \ref{deviation m-estimation}, rely on the uniform entropy integrals for the function classes $\FFF^{\Theta}, \FFF^{\Theta}_{\delta}$ which often increase with the dimension (see e.g. Propositions \ref{Hoelder-Bedingung auxiliary}, \ref{startingpoint}, \ref{important covering numbers}). 
\section{Error estimates under absolute semidevation risk measures}
\label{error estimates absolut semideviations}
In this section we study risk averse 
stochastic programs
\begin{align*}
\min_{\theta\in\Theta}\rho\big(G(\theta,Z_{1})\big),
\end{align*}
where in the objective the functional $\rho$ is an \textit{absolute semideviation risk measure}. This means that for $\lambda\in ]0,1]$ the functional $\rho = \rho_{1,\lambda}$ is defined as follows
$$
\rho_{1,\lambda}:L^{1}\OFP\rightarrow\RRR,~X\mapsto\EEE[X] + \lambda~ \EEE\big[\big(X - \EEE[X]\big)^{+}\big],
$$
where $L^{1}\OFP$ stands for the ordinary $L^{1}$-space on $\OFP$, tacitely identifying random variables which differ on $\MP$-null sets only. It is well-known that absolute semideviation risk measures are increasing w.r.t. the increasing convex order (cf. e.g. \cite[Theorem 6.51 along with Example 6.23 and Proposition 6.8]{ShapiroEtAl}). 

Introducing the notation 
\begin{equation}
\label{Hilfsgoals}
H_{\lambda}:\Theta\times\RRR\times\RRR^{d}\rightarrow\RRR,~(\theta,t,z)\mapsto G(\theta,z) + \lambda \big(G(\theta,z) - t\big)^{+}\quad(\lambda \in ]0,1]), 
\end{equation}
the stochastic program under absolute semidevation risk measure may also be viewed a compound 
stochastic program \eqref{composite optimization} with $H = H_{\lambda}$. 
\par
Let us fix $\lambda\in ]0,1]$. We assume that the stochastic program has a unique solution $\theta^{*}$, and we consider a sequence $\big(\widehat{\theta}_{n}\big)_{n\in\NNN}$ of $m$-dimensional random vectors $\widehat{\theta}_{n}$ which solve the empirical counterpart of the stochastic program w.r.t. the sample size $n$. The aim is to find upper estimations of the deviation probabilities by applying the results Theorems \ref{simpliest version rough}, \ref{simpliest version}. For this purpose we want to specialize some of the general assumptions in Section \ref{deviation probabilities compound SAA} to the choice $H = H_{\lambda}$.
\begin{itemize}
\item [(8i)] Requirement \ref{measurability H} is met under \ref{assump:A1} together with \ref{assump:A2}.
\item [(8ii)] Condition \ref{assump:A4} is already fulfilled with $L_{1} = L_{1}^{\lambda} :\equiv 2 (1 + \lambda)$ and $L_{2} = L_{2}^{\lambda} :\equiv 2 \lambda$.
\item[(8iii)] We may simplify \ref{assump:A5} by the following condition.
\begin{enumerate}[label=(A \arabic*+), ref=(A \arabic*+)]
\setcounter{enumi}{6}
\item \label{assump:A5+}Under \ref{assump:A2} with function $\xi^{G}$, there exist $\beta\in ]0,1]$, $\widehat{\overline{M}}_{1}, \overline{M}^{1} > 0$, and a family $(\xi^{G}_{\delta})_{\delta\in ]0,\Delta(\Theta)]}$ of square $\MP^{Z}$-integrable positive envelopes $\xi^{\Theta}_{\delta}$ of $\FFF_{\delta}^{\Theta}$ satisfying
$$
\|\xi_{\delta}^{G}\|_{\MP^{Z},2}\leq \widehat{\overline{M}}_{1}~\delta^{\beta}\quad\mbox{and}\quad J(\FFF_{\delta}^{\Theta},\xi_{\delta}^{G},1/8)\vee 
J(\FFF^{\Theta},\xi^{G},1/8) \leq \overline{M}^{1} 
$$
for $\delta\in ]0,\Delta(\Theta)]$.
\end{enumerate}
Under \ref{assump:A5+} assumption \ref{assump:A5} is fulfilled with the family $(\xi^{G}_{\delta})_{\delta\in ]0,\Delta(\Theta)]}$ of envelopes, the constant $\overline{M}^{1}$, and $\overline{M}_{1} = 2 (1 + \lambda)~\widehat{\overline{M}}_{1}$.
\item [(8iv)] Under \ref{assump:A5+} we have $\xi^{G}_{a,b,\delta} = 2 \sqrt{(1 + \lambda)^{2}~\xi^{G}_{\delta}(\cdot)^{2} + \lambda^{2}~(b + 2 (1 + \lambda)~\widehat{\overline{M}}_{1})^{2}~\delta^{2 \beta}}$, where $\xi^{G}_{a,b,\delta}$ is the mapping defined in \eqref{neue envelope}, and $\xi^{G}_{\delta}, \widehat{\overline{M}}_{1}, \beta$ are from \ref{assump:A5+}. This means $B_{n}^{\xi^{G}_{\delta}}\subseteq B_{n}^{\xi^{G}_{a,b,\delta}}$, using notation \eqref{Restmenge}. 
\item [(8v)] In \ref{assump:A6} the mapping $m_{\Theta\times\RRR}$ may be defined by $m_{\Theta\times\RRR}(\theta,x) := -\lambda [1 - F_{\theta}(x)]$. Hence it suffice to require the following property.
\begin{enumerate}[label=(A \arabic*+), ref=(A \arabic*+)]
\setcounter{enumi}{7}
\item \label{assump:A6+} There exist $\delta_{1} > 0, \widehat{K}_{\delta_{1}}\geq 0$
such that
\begin{align*}
&
\big|F_{\theta}\big(\theta,\EEE[G(\theta,Z_{1})] + x\big) - F_{\theta^{*}}\big(\theta^{*},\EEE[G(\theta^{*},Z_{1})] + x\big)\big|\\
&
\leq \widehat{K}_{\delta_{1}}~\big\|\big(\theta - \theta^{*},\EEE[G(\theta,Z_{1})]-\EEE[G(\theta^{*},Z_{1})]\big)\big\|_{m+1}\quad\mbox{for}~(\theta,x)\in\Theta\times [-\delta_{1},\delta_{1}].
\end{align*}
\end{enumerate}
In this case assumption \ref{assump:A6} is fulfilled with $\delta_{1}$ and $K_{\delta_{1}} = \lambda~\widehat{K}_{\delta_{1}}$.
\item [(8vi)] Under \ref{assump:A5+} with $\widehat{\overline{M}}_{1} > 0, \beta\in ]0,1]$, and \ref{assump:A6+} with $\delta_{1} > 0$, $\widehat{K}_{\delta_{1}}\geq 0$ the terms in \eqref{Mab}, \eqref{Mb} and \eqref{etaabbeta} may be rewritten by 
\begin{align}
&\label{6.2}
M_{a,b} = M_{a,b}^{\lambda} := \lambda \widehat{K}_{\delta_{1}}~a~[\Delta(\Theta)^{1-\beta} + 2(1 + \lambda)\widehat{\overline{M}}_{1}] + 2\lambda b\\
&\label{6.3}
\overline{M}_{b} = \overline{M}_{b}^{\lambda} := 2~\widehat{\overline{M}}_{1} (1 + \lambda) (1 + 2 \lambda) + 2\lambda b\\
&\label{6.4}
n(a,b,\beta) = \widehat{\widehat{n}}(a,b,\beta,\lambda) := \max\Big\{\overline{M}_{b}^{\lambda}/2, a^{2}/\delta_{1}^{2}, [(\overline{M}_{b}^{\lambda})^{2}~\Delta(\Theta)^{2\beta +1}]^{\frac{4-2\beta}{3-2\beta}}\Big\}
\end{align}
\end{itemize}
The application of the results from Section \ref{deviation probabilities compound SAA} leads to the following error estimates for the solutions $\widehat{\theta}_{n}$.
\begin{theorem}
Let $\lambda\in ]0,1]$, 
and let assumptions \ref{assump:A1}, \ref{assump:A2}, \ref{assump:A7}, \ref{assump:A3} as well as \ref{assump:A5+}, \ref{assump:A6+} be fulfilled with 
\begin{itemize}
\item positive envelope $\xi^{G}$ of $\FFF^{\Theta}$ as in \ref{assump:A2},
\item constant $M^{\lambda}_{1} = M_{1} > 0$ as in \ref{assump:A7},
\item constants $\widehat{\overline{M}}_{1} > 0, \beta\in ]0,1]$ from \ref{assump:A5+},
\item the family $\big(\xi^{G}_{\delta}\big)_{\delta\in ]0,\Delta(\Theta)]}$ of positive envelopes assumed in \ref{assump:A5+}.
\end{itemize}
Furthermore let $a,b,\varepsilon > 0$ with $a~2^{\beta + 1/2}\leq b\varepsilon^{\beta + 1/2}$,  
and let $M_{a,b}^{\lambda}, \overline{M}^{\lambda}_{b}, \widehat{\widehat{n}}(a,b,\beta,\lambda)$ be the constants introduced respectively in \eqref{6.2}, \eqref{6.3}, \eqref{6.4}, whereas $\mathfrak{g}(t), \widehat{\eta}$ stand for the constants in \eqref{frakg}, \eqref{etadachepsilon} respectively. Then, setting $\delta_{nk} := 2^{k}n^{-1/(4- 2\beta)}$ for $k\in\mathbb{Z}$, $n\in\NNN$, and using notations for the auxiliary events defined in \eqref{Hilfsereignis 1}, \eqref{Hilfsereignis 2} and \eqref{Restmenge} the following statements are valid.
\begin{itemize}
\item [1)] For $n\in\NNN$ with $n\geq \widehat{\widehat{n}}(a,b,\beta,\lambda)$ 
\begin{align*}
&
\MP\big(\big\{n^{1/(4 - 2\beta)}~ \|\widehat{\theta}_{n} - \theta^{*}\|_{m} > \varepsilon\big\}\big)\leq P^{\varepsilon}_{\lambda 1} + P^{\varepsilon}_{n 2} + \MP\big(\Omega\setminus\Omega_{n,a}\big), 
\end{align*}
where 
\begin{align*} 
&P^{\varepsilon}_{\lambda 1} := \frac{10 [
\overline{M}^{\lambda}_{b}~\widehat{\eta} + M_{a,b}^{\lambda}]}{M^{\lambda}_{1}}~\left(\frac{2}{\varepsilon}\right)^{3/2 - \beta} \sqrt{\frac{2}{\varepsilon\wedge 2}}
\\
& P^{\varepsilon}_{n 2} := \sum_{k= K_{\varepsilon} + 1\atop \delta_{n (k-1)}\leq\Delta(\Theta)}^{\infty}\hspace*{-0.35cm}\MP\big(\Omega\setminus\Omega^{n}_{\delta_{nk}\wedge\Delta(\Theta),\sqrt{2}^{k} b}\big).
\end{align*}
\item[2)]If $n\in\NNN$ with $n\geq \widehat{\widehat{n}}(a,b,\beta,\lambda)$, and if for $t > 0$ the inequality 
$$
2^{K_{\varepsilon} (2 - \beta) - K_{\varepsilon}^{+}/2} > 4 ~\big[M^{\lambda}_{a,b}+ \overline{M}^{\lambda}_{b}~\big(1 + 2 (t + 1)~\widehat{\eta}\big)\big]/M^{\lambda}_{1}
$$ 
holds, then
\begin{align*}
&
\MP\big(\big\{n^{1/(4 - 2\beta)}~ \|\widehat{\theta}_{n} - \theta^{*}\|_{m} > \varepsilon\big\}\big)\leq P^{t,\varepsilon}_{\lambda 1} + P^{\varepsilon}_{n 2} + P^{\varepsilon}_{n 3} + \MP\big(\Omega\setminus\Omega_{n,a}\big), 
\end{align*}
where 
\begin{align*}
&P^{t,\varepsilon}_{\lambda 1} :=
\exp\left(\frac{M^{\lambda}_{a,b}~\mathfrak{g}(t)}{(\sqrt{\varepsilon}\wedge 1)~2~\lambda~b}\right)\cdot \frac{2^{4- \beta}~\overline{M}^{\lambda}_{b}\cdot\exp\left(-\dfrac{\varepsilon^{3/2- \beta}\sqrt{\varepsilon\wedge 2}~M^{\lambda}_{1}~\mathfrak{g}(t)}{2^{4-\beta}~\overline{M}^{\lambda}_{b}}\right)}{(3/2-\beta)\ln(2)~M^{\lambda}_{1}~\varepsilon^{3/2-\beta}\sqrt{\varepsilon\wedge 2}},\\
&P^{\varepsilon}_{2 n} :=  \sum_{k= K_{\varepsilon} + 1\atop \delta_{n (k-1)}\leq\Delta(\Theta)}^{\infty}\hspace*{-0.35cm}\MP\big(\Omega\setminus\Omega^{n}_{\delta_{nk}\wedge\Delta(\Theta),\sqrt{2}^{k} b}\big)\quad\mbox{and}\quad
P^{\varepsilon}_{n 3} :=  \sum_{k= K_{\varepsilon} + 1\atop \delta_{n (k-1)}\leq\Delta(\Theta)}^{\infty}\hspace*{-0.35cm}\MP\big(\Omega\setminus B_{n}^{\xi^{G}_{\delta_{nk}\wedge\Delta(\Theta)}}\big).
\end{align*}
\end{itemize}
In particular the sequence $\big(n^{1/(4 - 2\beta)}~ [\widehat{\theta}_{n} - \theta^{*}]\big)_{n\in\NNN}$ is uniformly tight.
\end{theorem}
\begin{proof}
By (8i) - (8iii) together with (8v) we may apply Theorems \ref{simpliest version rough}, \ref{simpliest version} and Corollary \ref{allgemeine Konvergenzraten zweiter Teil}. Then statement 1) may be concluded immediately from Theorem \ref{simpliest version rough} along with (8ii), (8iii), (8v) and (8vi). Statement 2) follows from Theorem \ref{simpliest version} by (8ii) - (8vi), note in particular (8iv). The tightness of $\big(n^{1/(4 - 2\beta)}~ [\widehat{\theta}_{n} - \theta^{*}]\big)_{n\in\NNN}$ is a direct consequence of Corollary \ref{allgemeine Konvergenzraten zweiter Teil}. The proof is complete.
\end{proof}
\section{Error estimates under Average Value at Risk}
\label{estimates under Average Value at Risk}
As before 
$L^{1}\OFP$ denotes the usual $L^{1}$-space on $\OFP$, identifying the random variables which differ on $\MP$-null sets only. We consider the risk averse stochastic program 
\begin{align}
\label{AVAR-program}
\min_{\theta\in\Theta}\rho_{\alpha}\big(G(\theta,Z_{1})\big),
\end{align}
where the objective $\rho_{\alpha}$ is an \textit{Average Value at Risk}. More precisely, for some fixed $\alpha\in ]0,1[$ the functional $\rho_{\alpha}$ is defined by
$$
\rho_{\alpha}:L^{1}\OFP\rightarrow\RRR,~X\mapsto \int_{0}^{1}\eins_{]\alpha,1[}(u)~\Flinks_{X}(u)~du,
$$
where $\Flinks_{X}$ stands for the left-continuous quantile function of the distribution function $F_{X}$ of $X$. The method of sample average approximation may be adapted to optimization \eqref{AVAR-program} replacing each distribution function $F_{\theta}$ of $G(\theta,Z)$ with its empirical counterparts $\widehat{F}_{n,\theta}$ based on the samples $(Z_{1},\ldots,Z_{n})$. The resulting approximating optimization problems read as follows 
\begin{align}
\label{AVAR-SAA}
\min_{\theta\in\Theta}\int_{0}^{1}\eins_{]\alpha,1[}(u)~\big(\widehat{F}_{n,\theta}\big)^{\leftarrow}(u)~du\quad(n\in\NNN),
\end{align}
where $\big(\widehat{F}_{n,\theta}\big)^{\leftarrow}$ denotes the left-continuous quantile function of the empirical distribution function $\widehat{F}_{n,\theta}$. We assume that \eqref{AVAR-program} has a unique solution $\theta^{\alpha,*}$ and consider any sequence $(\widehat{\theta}^{\alpha}_{n})_{n\in\NNN}$ of random vectors $\widehat{\theta}^{\alpha}_{n}$ which solve optimization \eqref{AVAR-SAA} w.r.t. the sample size $n\in\NNN$. The aim is to upper estimate the deviation probabilities
\begin{align}
\label{deviations}
\MP\big(\big\{\|\widehat{\theta}^{\alpha}_{n} - \theta^{\alpha,*}\|_{m}  > \varepsilon~n^{-\gamma}\big\}\big)\quad(n\in\NNN,\varepsilon > 0)
\end{align}
for suitable $\gamma > 0$. The initial idea of our argumenation is a transformation procedure in two steps.
\begin{itemize}
\item Step 1: reformulation of the objective in \eqref{AVAR-program}
\item Step 2: definition of auxiliary problems
\end{itemize}
In step 1 we represent the stochastic program \eqref{AVAR-program} as a risk neutral stochastic program w.r.t. some auxiliary new parameter set. Afterwards we apply in step 2 the SAA method to the reformulated stochastic program. This will lead to new auxiliary SAA problems. 
\smallskip

Concerning step 1 it is well-known that the Average Value at Risk may also be represented by
\begin{align}
\label{AVAR-representation}
\rho_{\alpha}(X) = \inf_{x\in\RRR}\EEE\big[(X + x)^{+}/(1-\alpha) - x\big],
\end{align}
where the set of minimizers consists of all numbers $-x_{\alpha}$, denoting by $x_{\alpha}$ any $\alpha$-quantile of $X$ (see \cite{KainaRueschendorf2007}). In particular the Average Value at Risk is nondecreasing w.r.t. the increasing convex order, and 
\eqref{AVAR-program} may be reformulated by the auxiliary risk neutral stochastic program  
\begin{equation}
\label{auxiliary program}
\min_{(\theta,x)\in\Theta\times\RRR}\EEE\big[\big(G(\theta,Z_{1}) + x\big)^{+}/(1 - \alpha) - x\big].
\end{equation}
For step 2 we may also observe by representation \eqref{AVAR-representation} 
\begin{equation}
\label{optimization approximativ}
\int_{0}^{1}\eins_{]\alpha,1[}(u)~\big(\widehat{F}_{n,\theta}\big)^{\leftarrow}(u)~du
= \inf_{x\in\RRR}\frac{1}{n}\sum_{j = 1}^{n}\big[\big(G(\theta,Z_{j}) + x\big)^{+}/(1 - \alpha) - x\big]
\end{equation}
for $n\in\NNN$ and $\theta\in\Theta$. Then the empirical counterparts \eqref{AVAR-SAA} of \eqref{AVAR-program} may be described in terms of the following auxiliary SAA problems 
\begin{equation}
\label{auxiliary SAA}
\inf_{(\theta,x)\in\Theta\times\RRR}\frac{1}{n}\sum_{j = 1}^{n}\big[\big(G(\theta,Z_{j}) + x\big)^{+}/(1 - \alpha) - x\big]\quad(n\in\NNN).
\end{equation}

We strengthen the unique solvability of \eqref{AVAR-program} by 
\begin{equation}
\label{unique solution erweitert}
\eqref{AVAR-program}~\mbox{has a unique solution}~\theta^{\alpha,*}\quad\mbox{and}\quad F_{\theta^{\alpha,*}}~\mbox{has a unique}~\alpha-\mbox{quantile}.
\end{equation}
This requirement implies that the auxiliary optimization \eqref{auxiliary program}
has $\big(\theta^{\alpha,*},-\Flinks_{\theta^{\alpha,*}}(\alpha)\big)$ as its unique solution. For simplification we set $x^{\alpha,*} := -\Flinks_{\theta^{\alpha,*}}(\alpha)$, and we introduce the mapping 
\begin{equation}
\label{new goal}
\widehat{G}_{\alpha}:(\Theta\times\RRR)\times\RRR^{d}\rightarrow\RRR,~ \big((\theta,x),z\big) \mapsto [G(\theta,z) + x]^{+}/(1-\alpha) - x.
\end{equation}
Our intention is to study for suitable $\gamma > 0$ the probabilities 
\begin{equation}
\label{auxiliary probabilities}
\MP\Big(\Big\{\big\|\big(\widehat{\theta}_{n}^{\alpha} - \theta^{\alpha,*},\widehat{x}_{n}^{\alpha} - x^{\alpha,*}\big)\big\|_{m+1} > \varepsilon\cdot n^{-\gamma}\Big\}\Big)\quad(n\in\NNN, \varepsilon > 0)
\end{equation}
of the stochastic minimizers $(\widehat{\theta}_{n}^{\alpha},\widehat{x}_{n}^{\alpha})$ according to the problems \eqref{auxiliary SAA} instead of the deviation probabilities
for stochastic minimizers $\widehat{\theta}_{n}^{\alpha}$ of the problems \eqref{AVAR-SAA}. In view of the following auxiliary result we may restrict ourselves to the probabilities \eqref{auxiliary probabilities} if the objective $G$ is Borel measurable .
\begin{lemma}
\label{completing minimizers}
Let $n\in\NNN$, and let assumption \ref{assump:A2}, \ref{assump:A3} be fulfilled, where in addition $G$ is assumed to be measurable w.r.t. the product $\cB(\Theta)\otimes\cB(\RRR^{d})$ of the Borel $\sigma$-algebra $\cB(\Theta)$ on $\Theta$ and the Borel $\sigma$-algebra $\cB(\RRR^{d})$ on $\RRR^{d}$. Then for any random vector $\widehat{\theta}^{\alpha}_{n}$ on $\OFP$ which is a stochastic minimizer of problem \eqref{AVAR-SAA} w.r.t. the sample size $n$, there is some random variable $\widehat{x}_{n}^{\alpha}$ on $\OFP$ such that $(\widehat{\theta}^{\alpha}_{n},\widehat{x}^{\alpha}_{n})$ is a stochastic minimizer of the auxiliary problem \eqref{auxiliary SAA} w.r.t. the sample size $n$.
\end{lemma}
The proof of Lemma \ref{completing minimizers} is provided in Subsection \ref{Proof of Lemma completing minimizers}.
\medskip

The stochastic minimizers of the auxiliary problems \eqref{auxiliary SAA} may be viewed as $M$-estimators with objective $\widehat{G}_{\alpha}$. This suggests 
to utilize Theorem \ref{deviation m-estimation} to derive upper estimates for the probabilities \eqref{auxiliary probabilities}. Unfortunately the parameter space $\Theta\times\RRR$ in the optimization problems \eqref{auxiliary program}, \eqref{auxiliary SAA} is not totally bounded so that we may not apply this result directly. We already know that under \ref{assump:A1}, \ref{assump:A2} together with lower semicontinuity of $G$ in the parameter $\theta$, the optimization problems \eqref{auxiliary SAA} have bounded solution sets (see \cite{Kraetschmer2023a}). However, they may depend on the realizations of the samples. In \cite{Kraetschmer2023a} a kind of compactification was suggested which allows to restrict the parameter set in the optimizations to suitable compact subsets. The idea is to show that with arbitrarily high probability we may find for large sample sizes events from $\mathcal{F}$ on which all solution sets of the SAA problems \eqref{auxiliary SAA} are contained in a common compact superset. The following result from \cite{Kraetschmer2023a} gives a precise formulation of this idea. 
For preparation let for $\omega\in\Omega$ denote by $S_{n}(\omega)$ the set of solutions of the optimization problem
$$
\inf_{(\theta,x)\in\Theta\times\RRR}\Big\{\frac{1}{n}\sum_{j=1}^{n}\Big(G\big(\theta,Z_{j}(\omega)\big)+ x\Big)^{+}/(1-\alpha) - x\Big\}.
$$
Moreover, we introduce for any mapping $\xi^{G}$ as in \ref{assump:A2} and arbitrary $\alpha\in ]0,1[$ as well as $\eta > 0$ the auxiliary events
\begin{align}
\label{auxiliary events AVAR}
A_{n,\eta}^{\alpha,\xi^{G}} := \Big\{\frac{1}{n}\sum_{j=1}^{n}\xi^{G}(Z_{j})\leq \EEE[\xi^{G}(Z_{1})] + (1 - \alpha)~ \eta\Big\}\quad(n\in\NNN),
\end{align}
and the compact interval
\begin{equation}
\label{compact interval}
\mathcal{I}^{\alpha}_{\xi^{G},\eta} := \Big[-\eta - \frac{\EEE[\xi^{G}(Z_{1})]}{1-\alpha},\frac{(2-\alpha) \eta}{\alpha} + \frac{(2 - \alpha + |2- 3\alpha|)\EEE[\xi^{G}(Z_{1})]}{2\alpha (1-\alpha)} \Big].
\end{equation} 
The following result has been shown in \cite{Kraetschmer2023a} (Theorem 5.7 and Lemma 5.8 along with Example 4.1).
\begin{proposition}
\label{stochastic equicontinuity}
Let \ref{assump:A1}, \ref{assump:A2} be fulfilled with mapping $\xi^{G}$ from \ref{assump:A2}, and let $\alpha\in ]0,1[$. 
If $G(\cdot,z)$ is lower semicontinuous for $z\in\RRR^{d}$, the following statements are valid.
\begin{itemize}
\item [1)] $S_{n}(\omega)$ is nonvoid for $n\in\NNN, \omega\in\Omega$. 
\item [2)] Using notations 
\eqref{auxiliary events AVAR} and \eqref{compact interval} 
$$
S_{n}(\omega)\subseteq\Theta\times\mathcal{I}^{\alpha}_{\xi^{G},\eta}\quad\mbox{for}~n\in\NNN, \eta > 0, \omega\in A_{n,\eta}^{\alpha,\xi^{G}}. 
$$
Moreover, under assumption \eqref{unique solution erweitert} the unique solution $\big(\theta^{\alpha,*},-\Flinks_{\theta^{\alpha,*}}(\alpha)\big)$ always belongs to every parameter set $\Theta\times\mathcal{I}^{\alpha}_{\xi^{G},\eta}$ $(\eta > 0)$.
\end{itemize}
\end{proposition}
Proposition \ref{stochastic equicontinuity} tells us that up to probabilities $\MP(\Omega\setminus A_{n,\eta}^{\alpha,\xi^{G}})$ we may focus on M-estimators according to the problems
\begin{equation}
\label{HilfsSAA}
\min_{(\theta,x)\in\Theta\times\mathcal{I}^{\alpha}_{\xi^{G},\eta}}\frac{1}{n}~\sum_{j=1}^{n}\widehat{G}_{\alpha}\big((\theta,x),Z_{j}\big)\quad(n\in\NNN).
\end{equation}
However, in applying Theorem \ref{deviation m-estimation} we want to impose most of the conditions directly on the genuine objective $G$ instead of the auxiliary ones $\widehat{G}^{\alpha}_{|(\Theta\times\mathcal{I}^{\alpha}_{\xi^{G},\eta})\times\RRR^{d}}$.
\medskip

We start with introducing for any nonvoid interval $\mathcal{I}\subseteq\RRR$ 
the function classes 
\begin{align}
&\label{neue Funktionsklassen}
\FFF^{\Theta}_{\alpha, \mathcal{I}} :=\big\{\widehat{G}_{\alpha}\big((\theta,x),\cdot\big)\mid (\theta,x)\in \Theta\times\mathcal{I}\big\},
\end{align}
and, if $x^{\alpha,*}\in\mathcal{I}$,
\begin{align}
&\label{neue Funktionsklassen2}
\FFF^{\Theta}_{\alpha,\mathcal{I},\delta} := \big\{\widehat{G}_{\alpha}\big((\theta,x),\cdot\big) - \widehat{G}_{\alpha}\big((\theta^{\alpha,*},x^{\alpha,*}),\cdot\big) \mid (\theta,x)\in\cV_{\mathcal{I},\delta}\big\}\quad(\delta > 0),
\end{align}
where $\cV_{\mathcal{I},\delta} := \{(\theta,x)\in \Theta\times\mathcal{I}\mid \|(\theta - \theta^{\alpha,*},x - x^{\alpha,*})\|_{m+1}\leq\delta\}$ for $(\overline{\theta},\overline{x})\in\Theta\times\mathcal{I}$.
\smallskip

For application of Theorem \ref{deviation m-estimation} we have to verify at least the condition \ref{assump:A3} for the function classes $\FFF^{\Theta}_{\alpha,\mathcal{I}}$ and requirement \ref{assump:A5} for the classes $\FFF^{\Theta}_{\alpha,\mathcal{I}}$ as well as $\FFF^{\Theta}_{\alpha,\mathcal{I},\delta}$. The aim is to conclude the desired properties from assumptions on the genuine function classes $\FFF^{\Theta}$ and $\FFF^{\Theta}_{\delta}$.
\medskip
 
First of all, we want to fix an analogue of \ref{assump:A3} for the function classes $\FFF^{\Theta}_{\alpha,\mathcal{I}}$ in case that the genuine class $\FFF^{\Theta}$ fulfills the following version of this separability condition. 
\begin{enumerate}[label=(A \arabic*$\alpha$), ref=(A \arabic*$\alpha$)]
\setcounter{enumi}{5}
\item \label{assump:A3alpha}For arbitraty $\delta > 0$ there exist some at most countable subset $\cC(\cU^{\alpha}_{\delta})$ of the set $\cU^{\alpha}_{\delta} := \{\theta\in\Theta\mid \|\theta - \theta^{\alpha,*}\|_{m}\leq\delta\}$ and $(\MP^{Z})^{n}$-null sets $N_{n}^{\delta}$ $(n\in\NNN)$ such that 
$$
\inf_{\vartheta\in\cC(\cU^{\alpha}_{\delta})}\left\{\EEE[|G(\theta,Z_{1}) - G(\vartheta,Z_{1})|] + \max_{j\in\{1,\ldots,n\}}\big|G(\theta,z_{j}) - G(\vartheta,z_{j})\big|\right\} = 0\quad
$$
if $\theta\in\cU^{\alpha}_{\delta}$, $n\in\NNN$ and $(z_{1},\ldots,z_{n})\in\RRR^{d n}\setminus N_{n}^{\delta}$.
\end{enumerate}

%
We may observe the following basic inequalities for the functions from the classes $\FFF^{\Theta}_{\alpha,\mathcal{I}}$
\begin{equation}
\label{basic inequalities}
\big|\widehat{G}_{\alpha}\big((\theta,x),z\big) - \widehat{G}_{\alpha}\big((\vartheta,y),z\big)\big|\leq \big[\big|G(\theta,z) - G(\vartheta,z)\big| + (2 - \alpha)~|x-y|\big]/(1- \alpha)
\end{equation}
for $(\theta,x), (\vartheta,y)\in\Theta\times\mathcal{I}$ and $z\in\RRR^{d}$. Then in view of \eqref{basic inequalities} the analogue of \ref{assump:A3} for the classes $\FFF^{\Theta}_{\alpha,\mathcal{I}}$ is already implied by \ref{assump:A3alpha} for $\FFF^{\Theta}$.
\begin{lemma}
\label{countable parameter subsets}
Let \ref{assump:A1} - (A 3) be fulfilled, and let $\mathcal{I}\subseteq\RRR$ denote a nondegenerated interval with $x^{\alpha,*}\in\mathcal{I}$. Then under \ref{assump:A3alpha} the restriction of $\widehat{G}_{\alpha}$ to $(\Theta\times\mathcal{I})\times\RRR^{d}$ satisfies the separability property analogously to  \ref{assump:A3}.
\end{lemma}
The proof of Lemma \ref{countable parameter subsets} is postponed to the Subsection \ref{Proof of Lemma countable parameter subsets}.
\medskip

Next, let us consider the implicit requirement from \ref{assump:A5} that the function classes $\FFF^{\Theta}_{\alpha,\mathcal{I}}$ should have finite uniform entropy 
integrals. We borrow Lemma 4.3 from \cite{Kraetschmer2023a} which describes the relationship of the 
uniform entropy integrals $J(\FFF^{\Theta}_{\alpha,\mathcal{I}},C_{\FFF^{\Theta}_{\alpha,\mathcal{I}}},\varepsilon)$ with uniform entropy integrals  $J(\FFF^{\Theta},C_{\FFF^{\Theta}},\varepsilon)$. In particular the finiteness of $J(\FFF^{\Theta}_{\alpha,\mathcal{I}},C_{\FFF^{\Theta}_{\alpha,\mathcal{I}}},\varepsilon)$ is already guaranteed if $J(\FFF^{\Theta},C_{\FFF^{\Theta}},\varepsilon)$ if finite.
\begin{lemma}
\label{relationship}
Let $\mathcal{I}\subseteq\RRR$ be a nondegenerated compact interval fulfilling the property $\sup \mathcal{I} = |\inf \mathcal{I}|\vee |\sup \mathcal{I}| > 0$. If $\xi$ is a square $\MP^{Z}$-integrable positive envelope of $\FFF^{\Theta}$, then 
$
C_{\FFF^{\Theta}_{\alpha,\mathcal{I}}} := 2 (2 - \alpha)\sqrt{\xi^{2} + (\sup \mathcal{I})^{2}}/(1-\alpha)
$
is a positive envelope of $\FFF^{\Theta}_{\alpha,\mathcal{I}}$ satisfying 
\begin{align*}
J(\FFF^{\Theta}_{\alpha,\mathcal{I}},C_{\FFF^{\Theta}_{\alpha,\mathcal{I}}},\varepsilon)\leq \sqrt{2}~ J(\FFF^{\Theta},\xi,\varepsilon) + 4\varepsilon\sqrt{\ln(2/\varepsilon)} + \sqrt{2\ln(2)}~\varepsilon\quad\mbox{for}~\varepsilon \in ]0,\exp(-1)].
\end{align*}
\end{lemma}
We also want to find an explicit relationship between the uniform entropy integrals $J(\FFF^{\Theta}_{\alpha,\mathcal{I},\delta},C_{\FFF^{\Theta}_{\alpha,\mathcal{I},\delta}},\varepsilon)$ and $J(\FFF^{\Theta}_{\delta},C_{\FFF^{\Theta}_{\delta}},\varepsilon)$. This will be the subject of the following result.
\begin{lemma}
\label{relationship2}
Let $\mathcal{I}\subseteq\RRR$ be a nondegenerated compact interval containing $x^{\alpha,*}$ and satisfying $\sup \mathcal{I} = |\inf \mathcal{I}|\vee |\sup \mathcal{I}| > 0$. 
Furthermore let $\delta > 0$. If $\xi_{\delta}$ is a square $\MP^{Z}$-integrable positive envelope of $\FFF^{\Theta}_{\delta}$, then 
$
C_{\FFF^{\Theta}_{\alpha,\mathcal{I},\delta}} := [\xi_{\delta} + (2-\alpha)~\delta]/(1-\alpha)
$
is a square $\MP^{Z}$-integrable positive envelope of $\FFF^{\Theta}_{\alpha,\mathcal{I},\delta}$, and  
\begin{align*}
J(\FFF^{\Theta}_{\alpha,\mathcal{I},\delta},C_{\FFF^{\Theta}_{\alpha,\mathcal{I},\delta}},\varepsilon)\leq 4~ J(\FFF^{\Theta}_{\delta},\xi_{\delta},\varepsilon/4) + 2\varepsilon\sqrt{\ln(8/\varepsilon)}\quad\mbox{for}~ \varepsilon \in ]0,1].
\end{align*}
\end{lemma}
The proof of Lemma \ref{relationship2} may be found in Subsection \ref{Beweis relationship2}.
\medskip

Finally, we need a suitable second growth condition at the unique solution $(\theta^{\alpha,*},x^{\alpha,*})$. 
\begin{enumerate}[label=(A \arabic*$\alpha$), ref=(A \arabic*$\alpha$)]
\setcounter{enumi}{4}
\item \label{assump:A7alpha}
There exists $M_{1,\alpha} > 0$ such that 
\begin{align*}
\EEE\big[\widehat{G}_{\alpha}\big((\theta,x),Z_{1}\big)\big]- \EEE\big[\widehat{G}_{\alpha}\big((\theta^{\alpha,*},x^{\alpha,*}),Z_{1}\big)\big]
\geq M_{1,\alpha}~\|(\theta - \theta^{\alpha,*},x - x^{\alpha,*})\|_{m+1}^{2}
\end{align*}
for $(\theta,x)\in\Theta\times\RRR$. 
\end{enumerate}
\begin{remark}
\label{Isolierung}
Let $(\theta^{\alpha,*},x^{\alpha,*})$ be the unique minimizer of the problem \eqref{auxiliary program} and let $\xi^{G}$ the positive envelope of $\FFF^{\Theta}$ from \ref{assump:A2}. Then by definition of $\widehat{G}_{\alpha}$ along with 
Jensen's inequality, we may observe
$$
\inf_{\theta\in\Theta}\EEE\big[\widehat{G}_{\alpha}\big((\theta,x),Z_{1}\big)\big]\geq\max\big\{- x, \big(x-\EEE[\xi^{G}(Z_{1})]\big)^{+}/(1-\alpha) - x\big\}\quad\mbox{for}~x\in\RRR.
$$
This implies that the mapping $(\theta,x)\mapsto \EEE\big[\widehat{G}_{\alpha}\big((\theta,x),Z_{1}\big)\big]$ on $\Theta\times\RRR$ has bounded weak lower level sets. Furthermore by Vitalis theorem (see \cite[Proposition 21.4]{Bauer2001}) this mapping may be verified to be continuous w.r.t. the Euclidean metric under \ref{assump:A1} and \ref{assump:A2}. Thus in this case their weak lower level sets are even compact, in particular every minimizing sequence converges to its unique minimizer $(\theta^{\alpha,*},x^{\alpha,*})$. Hence for every $\delta > 0$
\begin{align*}
M^{\alpha}(\delta) := \inf_{(\theta,x)\in\Theta\times\RRR\atop \|(\theta - \theta^{\alpha,*},x - x^{\alpha,*})\|_{m+1} > \delta}\Big(\EEE\big[\widehat{G}_{\alpha}\big((\theta,x),Z_{1}\big)\big]- \EEE\big[\widehat{G}_{\alpha}\big((\theta^{\alpha,*},x^{\alpha,*}),Z_{1}\big)\big]\Big) > 0.
\end{align*}
This implies, following the argumentation in the proof of Lemma \ref{lokal gleich global}, that \ref{assump:A7alpha} may be derived by any local version of it at $(\theta^{\alpha,*},x^{\alpha,*})$, in the same way as \ref{assump:A7} follows from \ref{assump:A7'}.
\par
If $\theta^{\alpha,*}$ belongs to the topological interior of $\Theta$ then we may formulate exactly as is in Remark \ref{well-separated} criteria to provide property \ref{assump:A7alpha}. 
\end{remark}
\smallskip

Now, we are prepared to derive upper estimates for the errors of the stochastic minimizers of \eqref{optimization approximativ}. Actually, it is an application of Proposition \ref{stochastic equicontinuity} along with Theorem \ref{deviation m-estimation}. The relevant part consists of upper estimations of the deviation probabilities 
\eqref{auxiliary probabilities} with auxiliary events. Setting 
$\cV_{\delta}^{\eta} := \cV_{\mathcal{I}^{\alpha}_{\xi^{G},\eta},\delta}$ for $\delta, \eta > 0$, these events are defined by
\begin{align}
& \label{Hilfsereignis 3}
\overline{\Omega}^{\alpha,\eta}_{n,a} := \big\{\sup_{(\theta,x)\in\Theta\times\mathcal{I}^{\alpha}_{\xi^{G},\eta}}\big|(\MP_{n} - \MP)\big(\widehat{G}_{\alpha}(\theta,\cdot)\big)\big| \leq a/\sqrt{n}\big\}\quad(a, \eta > 0)\\
&\label{Hilfsereignis 4}
\overline{\Omega}^{\alpha,n,\eta}_{\delta,b} := \Big\{\sup_{(\theta,x)\in\cV_{\delta}^{\eta}}
\Big|(\MP_{n} - \MP)\Big(\widehat{G}_{\alpha}\big((\theta,x),\cdot\big) - \widehat{G}_{\alpha}\big((\theta^{\alpha,*},x^{\alpha,*}),\cdot\big)\Big)\Big| \leq b~\delta^{\beta}/\sqrt{n}\Big\} 
\end{align}
for $n\in\NNN$ and $a,b,\delta,\eta > 0$. These events take over the role of the auxiliary events in \eqref{Hilfsereignis 1}, \eqref{Hilfsereignis 2}. Lemma \ref{countable parameter subsets} ensures that 
every $\overline{\Omega}^{\alpha,\eta}_{n,a}$ and $\overline{\Omega}^{\alpha,n,\eta}_{\delta,b}$ belong to $\cF$ under condition \ref{assump:A3alpha} in the same way as measurability of the auxiliary events in \eqref{Hilfsereignis 1}, \eqref{Hilfsereignis 2} is shown in the proof of Lemma \ref{messbare Hilfsereignisse} below.
\begin{theorem}
\label{allgemeine Konvergenzraten m-estimation I}
Let (\ref{unique solution erweitert}) and assumptions \ref{assump:A1},\ref{assump:A2}, \ref{assump:A7alpha}, \ref{assump:A3alpha}, \ref{assump:A5},  be fulfilled with 
\begin{itemize}
\item 
$G$ being Borel measurable and lower semicontinuous in the parameter $\theta$;
\item $\xi^{G}$ denoting the square $\MP^{Z}$-integrable mapping from \ref{assump:A2},
\item constant $M_{1,\alpha} >0$ from \ref{assump:A7alpha},
\item constants $\overline{M}_{1}, \overline{M}^{1} > 0$, $\beta\in ]0,1]$ as in \ref{assump:A5},
\item family $\big(\xi^{G}_{\delta}\big)_{\delta\in ]0,\Delta(\Theta)])}$ of envelopes required in \ref{assump:A5}.
\end{itemize}
Furthermore 
let $a,b,\varepsilon,\eta > 0$ with $a~2^{\beta + 1/2}\leq b\varepsilon^{\beta + 1/2}$, and let the auxiliary event $A^{\alpha,\xi^{G}}_{n,\eta}$ as well as the interval $\mathcal{I}^{\alpha}_{\xi^{G},\eta}$ be as defined in \eqref{auxiliary events AVAR} and \eqref{compact interval}. 
Then, using notations 
\begin{itemize}
\item $\mathfrak{g}(t), \widehat{\eta}$ from \eqref{frakg}, \eqref{etadachepsilon} respectively,
\item for the auxiliary events introduced in \eqref{Restmenge} as well as in  \eqref{Hilfsereignis 3} and \eqref{Hilfsereignis 4},
\item $\delta_{nk} := 2^{k}/n^{1/(4-\beta)}\quad(k\in\mathbb{Z}, n\in\NNN)$,
\item $\Delta_{\alpha,\eta}~:=~$ diameter of $\Theta\times \mathcal{I}^{\alpha}_{\xi^{G},\eta}$ w.r.t. the Euclidean metric,
\item $\widehat{n}_{\alpha,\eta} := \frac{[\overline{M}_{1} + (2-\alpha)~\Delta_{\alpha,\eta}^{1-\beta}]^{2}\cdot\Delta_{\alpha,\eta}^{2\beta}}{2~(1-\alpha)^{2}}\vee\Big(\frac{[\overline{M}_{1} + (2-\alpha)~\Delta_{\alpha,\eta}^{1-\beta}]^{2}\cdot\Delta_{\alpha,\eta}^{2\beta + 1}}{(1-\alpha)^{2}}\Big)^{\frac{4-2\beta}{3 - 2\beta}}$,
\end{itemize}
the following upper estimations for the deviation probabilities \eqref{deviations} hold.
\begin{itemize}
\item [1)] For $n\in\NNN$ with $n\geq \widehat{n}_{\alpha,\eta}$ 
\begin{align*}
\MP\big(\big\{n^{1/(4 - 2\beta)}~ \|\widehat{\theta}^{\alpha}_{n} - \theta^{\alpha,*}\|_{m} > \varepsilon\big\}\big)\leq P^{\eta,\varepsilon}_{1} + P^{\eta,\varepsilon}_{n 2} +  \MP\big(\Omega\setminus \overline{\Omega}^{\alpha,\eta}_{n,a}\big) + \MP\big(\Omega\setminus A_{n,\eta}^{\alpha,\xi^{G}}\big),
\end{align*}
where 
\begin{align*} 
&P^{\eta,\varepsilon}_{1} := \frac{2^{2-\beta}~640~\sqrt{2}~[\overline{M}_{1} + (2-\alpha)~\Delta_{\alpha,\eta}^{1-\beta}]~(4~ \overline{M}^{1} + 1,7)}{(1-\alpha)~M_{1,\alpha}~\varepsilon^{3/2 - \beta}~\sqrt{\varepsilon\wedge 2}},\\
&P^{\eta,\varepsilon}_{n 2} := 
\sum_{k= K_{\varepsilon} + 1\atop \delta_{n (k-1)}\leq\Delta_{\alpha,\eta}}^{\infty}\hspace*{-0.35cm}\MP\big(\Omega\setminus\overline{\Omega}^{\alpha,n,\eta}_{\delta_{nk}\wedge\Delta_{\alpha,\eta},\sqrt{2}^{k} b}\big).
\end{align*}
\item[2)] If $n\in\NNN$ and $t>0$ such that $n\geq \widehat{n}_{\alpha,\eta}$ and
\begin{align*}
\frac{\varepsilon^{\frac{3}{2} - \beta}~\sqrt{\varepsilon\wedge 2}}{2^{4-\beta}} > 
\frac{\big(\overline{M}_{1} + (2-\alpha)~\Delta_{\alpha,\eta}^{1-\beta}\big)~\big[1 + 128\sqrt{2} (t + 1)~(4 \overline{M}^{1} + 1,7)\big]}{(1- \alpha)~M_{1,\alpha}},
\end{align*}
then
\begin{align*}
\MP\big(\big\{n^{1/(4 - 2\beta)}~ \|\widehat{\theta}_{n} - \theta^{\alpha,*}\|_{m} > \varepsilon\big\}\big)&\leq  P^{\eta,t,\varepsilon}_{1} + P^{\eta,\varepsilon}_{n 2} + P^{\eta,\varepsilon}_{n 3} \\ 
&\quad+ \MP\big(\Omega\setminus \overline{\Omega}^{\alpha,\eta}_{n,a}\big) + \MP\big(\Omega\setminus A_{n,\eta}^{\alpha,\xi^{G}}\big),
\end{align*}
with $P^{\eta,\varepsilon}_{n2}$ as in statement 1), and
\begin{align*}
&P^{\eta,t,\varepsilon}_{1} :=
\frac{2^{4 - \beta}~\big(\overline{M}_{1} + (2-\alpha)~\Delta_{\alpha,\eta}^{1-\beta}\big)~\exp\left(-\dfrac{\varepsilon^{3/2}~\sqrt{\varepsilon\wedge 2}~(1-\alpha)\cdot~M_{1,\alpha}\cdot\mathfrak{g}(t)}{2^{4 - \beta}~\big(\overline{M}_{1} + (2-\alpha)~\Delta_{\alpha,\eta}^{1-\beta}\big)}\right)}{(1-\alpha)~(3/2-\beta)~M_{1,\alpha}~\ln(2)~\varepsilon^{3/2-\beta}~\sqrt{\varepsilon\wedge 2}} \\
& 
P^{\eta,\varepsilon}_{n3} := \sum_{k= K_{\varepsilon} + 1\atop \delta_{n (k-1)}\leq\Delta_{\alpha,\eta}}^{\infty}\hspace*{-0.35cm}
\MP\big(\Omega\setminus B_{n}^{\xi^{G}_{\delta_{nk}\wedge\Delta_{\alpha,\eta}} + (2-\alpha)(\delta_{nk}\wedge\Delta_{\alpha,\eta})}\big).
\end{align*}
\end{itemize}
Moreover, value $P^{\eta}_{n2}$ and the probability $\MP\big(\Omega\setminus\overline{\Omega}^{\alpha,\eta}_{n,a}\big)$ satisfy the following inequalities for any $n\in\NNN$.

\begin{itemize}
\item [3)] 
$
\sum\limits_{k= K_{\varepsilon} + 1\atop \delta_{n (k-1)}\leq\Delta_{\alpha,\eta}}^{\infty}\hspace*{-0.35cm}\MP\big(\Omega\setminus\overline{\Omega}^{\alpha,n,\eta}_{\delta_{nk}\wedge\Delta_{\alpha,\eta},\sqrt{2}^{k} b}\big)
\leq 
\dfrac{128~(4~\overline{M}^{1} + 1)~[\overline{M}_{1} + (2-\alpha)~\Delta_{\alpha,\eta}^{1-\beta}]}{(1-\alpha)~(\sqrt{2} - 1)~b~\sqrt{\varepsilon}},
$
\item[4)] $\MP\big(\Omega\setminus\overline{\Omega}^{\alpha,\eta}_{n,a}\big)\leq 128\sqrt{2}~(2 - \alpha)~(4 \overline{M}^{1}+ 1)~\sqrt{\|\xi^{G}\|_{\MP^{Z},2}^{2} + (\sup\mathcal{I}^{\alpha}_{\xi^{G},\eta})^{2}}/[a (1-\alpha)]$.
\end{itemize}
The inequalities in statements 3), 4) may be improved in the following way.
\begin{itemize}
\item[5)] If $n\in\NNN$ with $n\geq [\overline{M}_{1} + (2-\alpha)~\Delta_{\alpha,\eta}^{1-\beta}]^{2}~\Delta_{\alpha,\eta}^{2}/[2~(1 - \alpha)^{2}]$, then for $t, \varepsilon > 0$ with 
$b > [\overline{M}_{1} + (2-\alpha)~\Delta_{\alpha,\eta}^{1-\beta}]~[1 + 64~\sqrt{2}~(t + 1)~(4~\overline{M}^{1} + 1)]/(\sqrt{\varepsilon} - \sqrt{\varepsilon}~\alpha)$
\begin{align*}
&
\sum_{k= K_{\varepsilon} + 1\atop \delta_{n (k-1)}\leq\Delta_{\alpha,\eta}}^{\infty}\hspace*{-0.35cm}\MP\big(\Omega\setminus\overline{\Omega}^{\alpha,n,\eta}_{\delta_{nk}\wedge\Delta_{\alpha,\eta},\sqrt{2}^{k} b}\big)\\
&\leq 
\dfrac{2~\sqrt{2}~[\overline{M}_{1} + (2-\alpha)~\Delta_{\alpha,\eta}^{1-\beta}]}{(1-\alpha)~\sqrt{\varepsilon}~\ln(2)~b~\mathfrak{g}(t)}~\exp\left(-\dfrac{(1-\alpha)~\sqrt{\varepsilon}~b~\mathfrak{g}(t)}{\sqrt{2}~[\overline{M}_{1} + (2-\alpha)~\Delta_{\alpha,\eta}^{1-\beta}]}\right) + P^{\eta,\varepsilon}_{n3},
\end{align*}
where $P^{\eta,\varepsilon}_{n3}$ is as in statement 2).
\item [6)] If $n\in\NNN$ with $n\geq 4~ (2-\alpha)^{2}~[\|\xi^{G}\|_{\MP,2}^{2} + (\sup\mathcal{I}^{\alpha}_{\xi^{G},\eta})^{2}]/(1-\alpha)^{2}$, and if $t > 0$ such that 
$a > 2~(2-\alpha)~\sqrt{\|\xi^{G}\|_{\MP,2}^{2} + (\sup\mathcal{I}^{\alpha}_{\xi^{G},\eta})^{2}}~[1 + 64\sqrt{2}~(t + 1)~(4 \overline{M}^{1} + 1)]/(1 - \alpha)$, then
$$
\MP(\Omega\setminus\overline{\Omega}^{\alpha,\eta}_{n,a})\leq 
\exp\Big(-\frac{ a~(1 - \alpha)\cdot\mathfrak{g}(t)}{2~(2 - \alpha)~\sqrt{\|\xi^{G}\|_{\MP^{Z},2}^{2} + (\sup\mathcal{I}^{\alpha}_{\xi^{G},\eta})^{2}}}\Big) + 
\MP\big(\Omega\setminus B_{n}^{\xi^{G}}\big).
$$
\end{itemize}
As a further consequence the sequence $\big(n^{1/(4- 2\beta)}[\widehat{\theta}_{n}^{\alpha} - \theta^{\alpha,*}]\big)_{n\in\NNN}$ is uniformly tight.
\end{theorem}
\begin{proof}
Since $G$ is assumed to be $\cB(\Theta)\otimes\cB(\RRR)$-measurable and lower semicontinuous in the parameter $\theta$, the mappings
\begin{align*}
(\Theta\times\mathcal{I}^{\alpha}_{\xi^{G},\eta})\times\RRR^{dn}\rightarrow\RRR,~\big((\theta,x),(z_{1},\ldots,z_{n})\big)\mapsto 
\frac{1}{n}\sum_{j=1}^{n}\widehat{G}_{\alpha}\big((\theta,x),z_{j}\big)\quad(n\in\NNN)
\end{align*}
are Borel measurable and lower semicontinuous in the parameter $(\theta,x)$. Furthermore $\Theta\times\mathcal{I}^{\alpha}_{\xi^{G},\eta}$, $\RRR^{dn}$ endowed with the topologies induced by the respective Euclidean norms are Polish spaces. Hence by compactness of $\Theta\times\mathcal{I}^{\alpha}_{\xi^{G},\eta}$ together with some specific measurable selection theorem (e.g. Corollary 1 in \cite{BrownPurves1973}) we may find for $n\in\NNN$ some Borel measurable mapping $(\widetilde{\theta}^{\alpha}_{n},\widetilde{x}^{\alpha}_{n}): \RRR^{dn}\rightarrow \Theta\times\mathcal{I}^{\alpha}_{\xi^{G},\eta}$ such that 
$(\widetilde{\theta}^{\alpha}_{n},\widetilde{x}^{\alpha}_{n})(Z_{1},\ldots,Z_{n})$ is a stochastic minimizer of the SAA problem \eqref{HilfsSAA} w.r.t. the sample size $n$. Moreover, by Lemma \ref{completing minimizers} there exists for every $n\in\NNN$ some random variable $\widehat{x}_{n}^{\alpha}$ on $\OFP$ such that $(\widehat{\theta}^{\alpha}_{n},\widehat{x}^{\alpha}_{n})$ is a stochastic minimizer of the auxiliary problem \eqref{auxiliary SAA} w.r.t. the sample size $n$. Now, we may conclude from Proposition \ref{stochastic equicontinuity} that
\begin{align*}
(\widehat{\theta}_{n},\widehat{x}_{n}) := \eins_{A^{\alpha,\xi^{G}}_{n,\eta}}\cdot (\widehat{\theta}^{\alpha}_{n},\widehat{x}^{\alpha}_{n}) + \eins_{\Omega\setminus A^{\alpha,\xi^{G}}_{n,\eta}}\cdot (\widetilde{\theta}^{\alpha}_{n},\widetilde{x}^{\alpha}_{n})(Z_{1},\ldots,Z_{n})
\end{align*}
defines another stochastic solution of \eqref{HilfsSAA} w.r.t. the sample size $n$. We may observe
\begin{align}
&\nonumber
\MP\big(\big\{n^{1/(4 - 2\beta)}~ \|\widehat{\theta}^{\alpha}_{n} - \theta^{\alpha,*}\|_{m} > \varepsilon\big\}\big)\\ 
&\nonumber 
\leq
\MP\big(\big\{n^{1/(4 - 2\beta)}~ \|(\widehat{\theta}^{\alpha}_{n} - \theta^{\alpha,*},\widehat{x}^{\alpha}_{n} - x^{\alpha,*})\|_{m + 1} > \varepsilon\big\}\big)\\ 
&\label{Kreuz0}\leq 
\MP\big(\big\{n^{1/(4 - 2\beta)}~ \|(\widehat{\theta}_{n} - \theta^{\alpha,*},\widehat{x}_{n} - x^{\alpha,*})\|_{m + 1} > \varepsilon\big\}\big) + \MP(\Omega\setminus A^{\alpha,\xi^{G}}_{n,\eta}).
\end{align}  
\par
For abbreviation we set $\mathcal{I} := \mathcal{I}^{\alpha}_{\xi^{G},\eta}$. It is a compact interval containing 
$x^{\alpha,*}$ by Proposition \ref{stochastic equicontinuity}, and satisfying $\sup\mathcal{I} = |\inf\mathcal{I}|\vee |\sup\mathcal{I}|$. Note that $(\theta^{\alpha,*},x^{\alpha,*})$ is the unique minimizer of 
the stochastic program
\begin{align*}
\inf_{(\theta,x)\in\Theta\times\mathcal{I}}\EEE\big[\widehat{G}_{\alpha}\big((\theta,x),Z_{1}\big)\big].
\end{align*}
We want to apply Theorem \ref{deviation m-estimation} to the function class 
$\FFF^{\Theta}_{\alpha,\mathcal{I}}$ defined in \eqref{neue Funktionsklassen}. First of all, the restriction of $\widehat{G}_{\alpha}$ to $\Theta\times\mathcal{I}$ satisfies the analogue of \ref{assump:A1} due to \ref{assump:A1} for $G$ along with \eqref{basic inequalities}. In Lemma \ref{relationship} it has been already pointed out how to construct a square $\MP^{Z}$-integrable positive envelope $C_{\FFF^{\Theta}_{\alpha,\mathcal{I}}}$ of $\FFF^{\Theta}_{\alpha,\mathcal{I}}$ out of the positive envelope $\xi^{G}$ from \ref{assump:A2}. In addition we have
\begin{align}
&\label{Kreuz00}
\|C_{\FFF^{\Theta}_{\alpha,\mathcal{I}}}\|_{\MP^{Z},2} = 2 (2 - \alpha)~\sqrt{\EEE[\xi^{G}(Z_{1})^{2}] + (\sup\mathcal{I})^{2}}/(1-\alpha).
\end{align}
Moreover, Lemma \ref{countable parameter subsets} ensures that the analogue of \ref{assump:A3} is fulfilled under \ref{assump:A3alpha}. Condition \ref{assump:A7alpha} provides the second order growth condition analogously to \ref{assump:A7}. So it is left to look at the analogue of property \ref{assump:A5}. 
\smallskip

For this purpose consider the family $\big(\xi^{G}_{\delta}\big)_{\delta\in ]0,\Delta(\Theta)}$ of positive square $\MP^{Z}$-integrable envelopes $\xi^{G}_{\delta}$ of $\FFF^{\Theta}_{\delta}$ from assumption \ref{assump:A5}. Since $\FFF^{\Theta}_{\delta}$ coincides with $\FFF^{\Theta}_{\Delta(\Theta)}$ for $\delta > \Delta(\Theta)$, we may 
choose $\xi^{G}_{\delta} := \xi^{G}_{\Delta(\Theta)}$ as a square $\MP^{Z}$-integrable positive envelope of $\FFF^{\Theta}_{\delta}$ whenever $\delta > \Delta(\Theta)$. Then in view of Lemma \ref{relationship2} a square $\MP^{Z}$-integrable positive envelope of the function class $\FFF^{\Theta}_{\alpha,\mathcal{I},\delta}$ from \eqref{neue Funktionsklassen2} is provided by $C_{\FFF^{\Theta}_{\alpha,\mathcal{I},\delta}} := 
[\xi^{G}_{\delta} + (2-\alpha)~\delta]/(1-\alpha)$ for $\delta > 0$.
\par
By Lemma \ref{relationship} along with \ref{assump:A5} we may observe 
\begin{align*}
J\big(\FFF^{\Theta}_{\alpha,\mathcal{I}},C_{\FFF^{\Theta}_{\alpha,\mathcal{I}}},1/8\big)
\leq 
\sqrt{2}~J\big(\FFF^{\Theta},\xi^{G},1/8\big) + \sqrt{\ln(16)}/2 + \sqrt{2\ln(2)}/8
\leq \sqrt{2}~\overline{M}^{1} + 1,
\end{align*}
and by \ref{assump:A5} together with Lemma \ref{relationship2}
\begin{align}
&\label{Kreuz1}
\|C_{\FFF^{\Theta}_{\alpha,\mathcal{I},\delta}}\|_{\MP^{Z},2}\leq [\overline{M}_{1} + (2-\alpha)~\Delta_{\alpha,\eta}^{1-\beta}]~\delta^{\beta}/(1-\alpha)\\
&\nonumber
J\big(\FFF^{\Theta}_{\alpha,\mathcal{I},\delta},C_{\FFF^{\Theta}_{\alpha,\mathcal{I},\delta}},1/8)
\leq 
4~J\big(\FFF^{\Theta}_{\delta},C_{\FFF^{\Theta}_{\delta\wedge\Delta(\Theta)}},1/32\big) +\sqrt{\ln(64)}/4
\leq 
4~\overline{M}^{1}+ 1
\end{align}
for $\delta\in ]0,\Delta_{\alpha,\eta}]$. In particular
\begin{align}
\label{Kreuz2}
J\big(\FFF^{\Theta}_{\alpha,\mathcal{I},\delta},C_{\FFF^{\Theta}_{\alpha,\mathcal{I},\delta}},1/8\big)~\vee~J\big(\FFF^{\Theta}_{\alpha,\mathcal{I}},C_{\FFF^{\Theta}_{\alpha,\mathcal{I}}},1/8\big)
\leq 4 \overline{M}^{1} + 1\quad\mbox{for}~\delta\in ]0,\Delta_{\alpha,\eta}].
\end{align}
Hence the analogue of \ref{assump:A5} is fulfilled with $\beta$ from \ref{assump:A5} and the family of positive envelopes $(C_{\FFF^{\Theta}_{\alpha,\mathcal{I},\delta}})_{\delta\in ]0,\Delta_{\alpha,\eta}]}$ associated with the family 
$(\FFF^{\Theta}_{\alpha,\mathcal{I},\delta})_{\delta\in ]0,\Delta_{\alpha,\eta}]}$ of function classes $\FFF^{\Theta}_{\alpha,\mathcal{I},\delta}$, where $\overline{M}_{1,\alpha} := [\overline{M}_{1} + (2-\alpha)~\Delta_{\alpha,\eta}^{1-\beta}]/(1-\alpha)$ and $\overline{M}^{1,\alpha} := 4 \overline{M}^{1} + 1$ take the role of $\overline{M}_{1}$ and $\overline{M}^{1}$ respectively. In particular notation \eqref{etadachepsilon} reads as follows
\begin{align}
\label{neues eta dach}
64\sqrt{2}~(\overline{M}^{1,\alpha} + 0,7) = 64\sqrt{2}~(4 \overline{M}^{1} + 1,7).
\end{align}
\smallskip

All in all we are in the position to apply Theorem \ref{deviation m-estimation} to the function classes $\FFF^{\Theta}_{\alpha,\mathcal{I}}$ and $\FFF^{\Theta}_{\alpha,\mathcal{I},\delta}$. Then statement 1) follows from 
statement 1) in Theorem \ref{deviation m-estimation} along with \eqref{Kreuz0} and \eqref{neues eta dach}. 
In view of \eqref{Kreuz0} and \eqref{neues eta dach} the statement 2) may be obtained easily from statement 2) in Theorem \ref{deviation m-estimation}.
\par
Next, concerning statements 3) - 6) we apply Remark \ref{erste Restwahrscheinlichkeiten} to $\FFF^{\Theta}_{\alpha,\mathcal{I}}$ and $\FFF^{\Theta}_{\alpha,\mathcal{I},\delta}$. Statement 4) may be derived from statement 1) there along with \eqref{Kreuz00}. Furthermore, drawing on the inclusion $\Omega\setminus B_{n}^{C_{\FFF^{\Theta}_{\alpha,\mathcal{I}}}}\subseteq \Omega\setminus B_{n}^{\xi^{G}}$ and \eqref{Kreuz00}, statement 6) comes from statement 4) in Remark \ref{erste Restwahrscheinlichkeiten}. Finally, statements 3), 5) may be concluded from statements 3), 6) in Remark \ref{erste Restwahrscheinlichkeiten}.
\smallskip

Implicitly we have already shown that the function class $\FFF^{\Theta}_{\alpha,\mathcal{I}}$ also meets the requirements of Corollary \ref{allgemeine Konvergenzraten zweiter Teil}. Hence 
the sequence $\big(n^{1/(4 - 2\beta)}~[(\widehat{\theta}_{n} - \theta^{\alpha,*},\widehat{x}_{n}- x^{\alpha,*})]\big)_{n\in\NNN}$ is uniformly tight. Then $\big(n^{1/(4 - 2\beta)}~[\widehat{\theta}_{n} - \theta^{\alpha,*}]\big)_{n\in\NNN}$ is uniformly tight due to \eqref{Kreuz0} because $\MP\big(\Omega\setminus A_{n,\eta}^{\alpha,\xi^{G}}\big)\to 0$ by law of large numbers.
This completes the proof.
\end{proof}
\begin{remark}
Analogously to Remark \ref{Simplifications} we may invoke classical concentration inequalities to provide in Theorem \ref{allgemeine Konvergenzraten m-estimation I} upper estimations for the probabilities $\MP\big(\Omega\setminus A^{\alpha,\xi^{G}}_{n,\eta}\big)$.
\end{remark}
We want to discuss the assumptions of Theorem \ref{allgemeine Konvergenzraten m-estimation I} if $G$ has some specific representation. Firstly, we look at representation (H).
\begin{example}
\label{EXH}
Let $G$ have representation (H) with $G(\theta,\cdot)$ being Borel measurable. Then by continuity of the paths it is known to be $\cB(\Theta)\otimes\cB(\RRR^{d})$-measurable (e.g. \cite[Lemma 6.7.3]{Pfanzagl1994}), and \ref{assump:A2} is fulfilled due to Proposition \ref{Hoelder-Bedingung auxiliary}. Furthermore, Remark \ref{Remark Hoelder-Bedingung} shows how to verify \ref{assump:A5}. Finally, \ref{assump:A3alpha} is an easy consequence of (H) along with the dominated convergence theorem because $\Theta$ encloses some at most countable dense subset.
\end{example}
Next, we want to consider the case when $G$ has representation (PH).
\begin{example}
\label{EXPH}
Let $G$ be representable by (PH), being also lower semicontinuous in the parameter $\theta$ and satisfying condition (*) in Remark \ref{Remark PH}. Continuity of the paths implies that the mappings $G^{1},\ldots,G^{r}$ in the representation (PH) are $\cB(\Theta)\otimes\cB(\RRR^{d})$-measurable. The constraints in the representation (PH) are defined in terms of functions which are $\cB(\Theta)\otimes\cB(\RRR^{d})$-measurable. Putting all together $G$ is also $\cB(\Theta)\otimes\cB(\RRR^{d})$-measurable. In Remark \ref{Remark PH} it is discussed how to meet the requirements \ref{assump:A1}, \ref{assump:A2}, \ref{assump:A5}. The condition \ref{assump:A3alpha} follows easily from (*).
\end{example}

\section{Concentration inequalities for locally small increments of compound empirical processes}
\label{MainResult}
Let \eqref{unique solution} and assumptions \ref{assump:A1} - \ref{assump:A6} be fulfilled with $\beta\in ]0,1]$ from \ref{assump:A5}, and $\theta^{*}$ denoting the unique solution of the stochastic programm \eqref{composite optimization} according to \eqref{unique solution}. As pointed out in Proposition \ref{variational inequalities} the main task in deriving the upper estimations of the deviation probabilities in Theorems \ref{simpliest version rough}, \ref{simpliest version} is to find suitable concentration inequalities for locally small increments of the compound empirical processes $\GGG_{n}^{G,H}$. More precisely, recalling
$\cU_{\delta} := \{\theta\in\Theta\mid \|\theta - \overline{\theta}\|_{m}\leq\delta\}$, we are interested in upper estimations of the following deviation probabilities
\begin{align}
\label{deviation probabilities}
\MP\Big(\Big\{\sup_{\theta\in\cU_{\delta}}\big|\GGGH_{n}(\theta,\cdot) - \GGGH_{n}(\theta^{*},\cdot)\big| > \varepsilon~n^{-\beta/(4 - 2\beta)}\Big\}\Big)\quad(\varepsilon > 0)
\end{align} 
for small $\delta > 0$. Assumptions \ref{assump:A1} - \ref{assump:A6} allow to find upper bounds of the probabilities in \eqref{deviation probabilities}. They will be described in terms of the sample size and explicit constants. 
\begin{theorem}
\label{Main Theorem}
Let assumptions \ref{assump:A1} - \ref{assump:A4}, \ref{assump:A3} - \ref{assump:A6} be fulfilled with 
\begin{itemize}
\item the positive envelope $\xi^{G}$ of $\FFF^{\Theta}$ as in \ref{assump:A2},
\item the square $\MP^{Z}$-integrable mappings $L_{1}, L_{2}$ from \ref{assump:A4},
\item the constants $\overline{M}_{1}, \overline{M}^{1} > 0$ as well as $\beta\in ]0,1]$ from \ref{assump:A5},
\item the family $\big(\xi^{G}_{\delta}\big)_{\delta\in ]0,\delta_{1}]}$ of positive envelopes as in \ref{assump:A5},
\item the constants $\delta_{1} > 0, K_{\delta_{1}}\geq 0$ from \ref{assump:A6}.
\end{itemize}
Finally, let $a, b, \varepsilon > 0$. Then, using notations \eqref{Hilfsereignis 1} - 
\eqref{neue envelope}, 
the following statements hold if $n\in\NNN$ with  $n\geq\max\{\overline{M}_{b}^{2}~\Delta(\Theta)^{2\beta}/2, \eins_{]0,\infty[}(\EEE[L_{2}(Z_{1})])\cdot a^{2}/\delta_{1}^{2}\}$.
\begin{itemize}
\item [1)]
For arbitrary $\delta\in ]0,\Delta(\Theta)]$
\begin{align*}
&\MP^{*}\Big(\Big\{\sup_{\theta\in\cU_{\delta}}\big|\GGGH_{n}(\theta,\cdot) - \GGGH_{n}(\theta^{*},\cdot)\big| > \varepsilon\cdot n^{-\beta/(4 - 2\beta)}\Big\}\cap\Omega_{n,a}\cap\Omega^{n}_{\delta,b}\Big)\\
&\leq 
\big[\overline{M}_{b}~\overline{\eta}_{n}(a,b,\delta) + \eins_{]0,\infty[}(\EEE[L_{2}(Z_{1})])\cdot M_{a,b}\big]~\delta^{\beta}~n^{\beta/(4 - 2\beta)}/\varepsilon.
\end{align*}  
\item [2)] For every 
$\delta\in ]0,\Delta(\Theta)]$, and arbitrary $\varepsilon, t > 0$
\begin{align*}
&\MP^{*}\Big(\Big\{\sup_{\theta\in\cU_{\delta}}\big|\GGGH_{n}(\theta,\cdot) - \GGGH_{n}(\theta^{*},\cdot)\big| > \varepsilon\cdot n^{-\beta/(4 - 2\beta)}\Big\}\cap\Omega_{n,a}\cap\Omega^{n}_{\delta,b}\Big)\\
&\leq
\exp\left(\eins_{]0,\infty[}\big(\EEE[L_{2}(Z_{1})]\big)\cdot \frac{M_{a,b}\cdot\mathfrak{g}(t)}{\overline{M}_{b}} \right)~\exp\left(-\frac{\varepsilon\cdot\mathfrak{g}(t)}{n^{\beta/(4 - 2\beta)}~\delta^{\beta}~\overline{M}_{b}}\right)
+ \MP\big(\Omega\setminus B_{n}^{\xi^{G}_{a,b,\delta}}\big)
\end{align*}
in case of $\delta < \Big(\varepsilon/\big[\eins_{]0,\infty[}\big(\EEE[L_{2}(Z_{1})]\big)\cdot M_{a,b} + \overline{M}_{b}\cdot \eta_{t,n}(a,b,\delta)\big]\Big)^{1/\beta}\cdot n^{-1/(4 - 2\beta)}$. Here $\eta_{t,n}(a,b,\delta) := 1 + 2 (t + 1)~\overline{\eta}_{n}(a,b,\delta)$ for $t, \delta > 0$.
\end{itemize}
In both statements $\MP^{*}$ denotes the outer probability w.r.t. $\MP$.
\end{theorem}
The proof of Theorem \ref{Main Theorem} may be found in the Subsection \ref{proof of Theorem Main Theorem rough}.
\smallskip


\section{Proofs}
\label{Beweise}
Recall the notions and notations from empirical process theory introduced in the Section \ref{deviation probabilities compound SAA}. 
For $n\in\NNN$ we introduce the random function $\MMM_{n}: \Theta\times\Omega\rightarrow\RRR$ via 
\begin{equation}
\label{random function}
\MMM_{n}(\theta,\omega) = \frac{1}{n}\sum_{j=1}^{n}H\Big(\theta,\frac{1}{n}~\sum_{k=1}^{n}G\big(\theta,Z_{k}(\omega)\big),Z_{j}(\omega)\Big).
\end{equation} 
\medskip

The proofs of the main results are repeated applications of Theorems 2.1, 2.2 from \cite{Kraetschmer2023a}. We present present now strengthened versions of them tailored to our situation.
\begin{theorem}
\label{upper estimation}
Let $\Gamma\subseteq\RRR^{k}$ be compact, and let $\overline{G}^{\Gamma}: \Gamma\times\RRR^{d}\rightarrow\RRR$ satisfy the following properties
\begin{itemize}
\item [1)] $\overline{G}^{\Gamma}(\gamma,\cdot)$ is Borel measurable for every $\gamma\in\Gamma$. 
\item [2)] The associated function class $\overline{\FFF}^{\Gamma} := \{\overline{G}^{\Gamma}(\gamma,\cdot)\mid \gamma\in\Gamma\}$ has a square $\MP^{Z}$-integrable positive envelope 
$\xi$ with finite uniform entropy integral $J(\overline{\FFF}^{\Gamma},\xi,1/2)$.
\item [3)] There exist some at most countable subset $\overline{\Gamma}\subseteq\Gamma$ and $(\MP^{Z})^{n}$-null sets $N_{n}$ ($n\in\NNN$) such that 
$$
\inf_{\gamma\in\overline{\Gamma}}\left\{\big|\EEE[\overline{G}^{\Gamma}(\gamma,Z_{1})] - \EEE[\overline{G}^{\Gamma}(\overline{\gamma},Z_{1})]\big| + \max_{j\in\{1,\ldots,n\}}\big|\overline{G}^{\Gamma}(\gamma,z_{j}) - \overline{G}^{\Gamma}(\overline{\gamma},z_{j})\big|\right\} = 0
$$
for $\overline{\gamma}\in\Gamma$, $n\in\NNN$ and $(z_{1},\ldots,z_{n})\in\RRR^{d n}\setminus N_{n}$.
\end{itemize}
Furthermore, set $B_{n}^{\xi} := \Big\{\frac{1}{n}\sum_{j=1}^{n}\xi(Z_{j})^{2}\leq  2\EEE[\xi(Z_{1})^{2}]\Big\}$, and let $\mathfrak{g}(t)\in\RRR$ be defined as in \eqref{frakg} for $t > 0$. Then $\sup_{\gamma\in\Gamma}\big|(\MP_{n} - \MP)\big(\overline{G}^{\Gamma}(\gamma,\cdot)\big)\big|$ is a random variable on $\OFP$ such that
\begin{align*}
\EEE\Big[\sup_{\gamma\in\Gamma}\big|(\MP_{n} - \MP)\big(\overline{G}^{\Gamma}(\gamma,\cdot)\big)\big|\Big]\leq 16\sqrt{2}~\|\xi\|_{\MP^{Z},2}~J(\overline{\FFF}^{\Gamma},\xi,1/2)/\sqrt{n},
\end{align*}
is valid for $n\in\NNN$, and
\begin{align*}
\MP\Big(\Big\{\sup_{\gamma\in\Gamma}\big|(\MP_{n} - \MP)\big(\overline{G}^{\Gamma}(\gamma,\cdot)\big)\big| > \varepsilon\Big\}\Big)
\leq
\exp\left(\frac{-\sqrt{n}\varepsilon}{\|\xi\|_{\MP^{Z},2}}\cdot \mathfrak{g}(t)\right)
+ \MP\big(\Omega\setminus B_{n}^{\xi}\big)
\end{align*}
holds for $t > 0$ and arbitrary $n\in\NNN$ with $\varepsilon > \eta_{t,n}$ as well as 
$n\geq \|\xi\|_{\MP^{Z},2}^{2}/2$, where 
$$
\eta_{t,n} := \|\xi\|_{\MP^{Z},2}~\big[1 + 32\sqrt{2} (1+t) ~J(\overline{\FFF}^{\Gamma},\xi,1/4)\big]/\sqrt{n}.
$$
\end{theorem}
\begin{proof}
Let the at most countable subset $\overline{\Gamma}$ of $\Gamma$ and the $(\MP^{Z})^{n}$-null sets $N_{n}$ ($n\in\NNN$) be from condition 3). 
Fix $n\in\NNN$. By 3) we may find for any $\overline{\gamma}\in\Gamma$ some sequence $(\gamma_{l})_{l\in\NNN}$ in $\overline{\Gamma}$ such that
\begin{equation*}
\lim_{l\to\infty}\big|\big(\MP_{n} - \MP\big)\big(\overline{G}^{\Gamma}(\gamma_{l},\cdot)\big)_{|\omega}\big| = \big|\big(\MP_{n} - \MP\big)\big(\overline{G}^{\Gamma}(\overline{\gamma},\cdot)\big)_{|\omega}\big|
\end{equation*}
for $\omega\in\{(Z_{1},\ldots,Z_{n})\in\RRR^{d n}\setminus N_{n}\}$. Hence 
\begin{align}
\label{Stern1}
\sup_{\gamma\in\Gamma}\big|(\MP_{n} - \MP)\big(\overline{G}^{\Gamma}(\gamma,\cdot)\big)\big| 
=
\sup_{\gamma\in\overline{\Gamma}}\big|(\MP_{n} - \MP)\big(\overline{G}^{\Gamma}(\gamma,\cdot)\big)\big|\quad\MP-\mbox{a.s.}.
\end{align}
In particular, $\sup_{\gamma\in\Gamma}\big|(\MP_{n} - \MP)\big(\overline{G}^{\Gamma}(\gamma,\cdot)\big)\big|$ is a random variable on $\OFP$ due to completeness of $\OFP$. 
Now, in view of \eqref{Stern1} the remaining part of Theorem \ref{upper estimation} follows immediately from Theorems 2.1, 2.2 in \cite{Kraetschmer2023a}.
\end{proof}
For ease of reference we may fix the measurability of the auxiliary events $\Omega_{n,a}$ and $\Omega^{n}_{\delta,b}$ via Theorem \ref{upper estimation}.
\begin{lemma}
\label{messbare Hilfsereignisse}
Let \eqref{unique solution} and assumptions \ref{assump:A1}, \ref{assump:A2}, \ref{assump:A3}, \ref{assump:A6} be fulfilled. Then $\Omega_{n,a}, \Omega^{n}_{\delta,b}\in\cF$ for $a,b > 0$ and $\delta\in ]0,\Delta(\Theta)], n\in\NNN$.
\end{lemma}
\begin{proof}
The function classes $\FFF^{\Theta}$ and $\FFF^{\Theta}_{\delta}$ $(\delta\in ]0,\Delta(\Theta)])$ meet the requirements of Theorem \ref{upper estimation}.
\end{proof}
\subsection{Proof of Proposition \ref{variational inequalities}}
\label{Proof of Proposition variational inequalities}
Let $n\in\NNN$ and $\varepsilon, \gamma\in ]0,\infty[$. Without loss of generality we may assume $\varepsilon/n^{\gamma} < \Delta(\Theta)$. Fix any $\overline{\Omega}\in\cF$, and let $\omega\in\Omega$ 
with $\varepsilon/n^{\gamma} < \|\widehat{\theta}_{n}(\omega) - \theta^{*}\|_{m}\leq \Delta(\Theta)$. Then the inequalities $2^{k-1}/n^{\gamma} < \|\widehat{\theta}_{n}(\omega) - \theta^{*}\|_{m}\leq \Delta(\Theta)\wedge (2^{k}/n^{\gamma})$ hold for some $k\in\mathbb{Z}$ satisfying $k\geq K_{\varepsilon} +1$ and $2^{k-1}\leq \Delta(\Theta)~n^{\gamma}$. By definition the realizations of $\widehat{\theta}_{n}$ minimize the paths of the random function $\MMM_{n}$ defined in \eqref{random function}.
Hence by \ref{assump:A7}
\begin{align*}
&\nonumber
\sup_{\theta\in\cU_{\Delta(\Theta)\wedge (2^{k}/n^{\gamma})}}\big|\MMM_{n}(\theta,\omega) - \MMM_{n}(\theta^{*},\omega) - \psi_{H,\Theta}(\theta) + \psi_{H,\Theta}(\theta^{*})\big|\\
&\nonumber 
\geq 
\MMM_{n}(\theta^{*},\omega) - \MMM_{n}\big(\widehat{\theta}_{n}(\omega),\omega\big) - \psi_{H,\Theta}(\theta^{*}) + \psi_{H,\Theta}\big(\widehat{\theta}_{n}(\omega)\big)\\
&
\geq 
\psi_{H,\Theta}\big(\widehat{\theta}_{n}(\omega)\big) - \psi_{H,\Theta}(\theta^{*})
>
M_{1}~2^{2(k-1)}/n^{2\gamma}.
\end{align*}
Finally, note that $\GGG_{n}^{G,H}(\theta,\cdot) =\sqrt{n}[\MMM_{n}(\theta,\cdot) - \psi_{H,\Theta}(\theta)]$ is valid for $\theta\in\Theta$. Now, 
we may derive the statement of Proposition \ref{variational inequalities} easily.
\hfill$\Box$
\subsection{Proof of Theorem \ref{Main Theorem}}
\label{proof of Theorem Main Theorem rough}
Let us retake assumptions and notations from the display of Theorem \ref{Main Theorem}. 
In addition $\MMM_{n}$ denotes the random function defined by \eqref{random function}.
\medskip

As a first observation, recalling the events $\Omega_{n,a}, \Omega^{n}_{\delta,b}$ from \eqref{Hilfsereignis 1} and \eqref{Hilfsereignis 2}, 
$$
\Big(\theta,\frac{1}{n}\sum_{j=1}^{n}G\big(\theta,Z_{j}(\omega)\big),\frac{1}{n}\sum_{j=1}^{n}G\big(\theta^{*},Z_{j}(\omega)\big) \Big)\in \mathcal{W}_{a,b}^{n,\delta}\quad\mbox{for}~\theta\in\cU_{\delta}, \omega\in\Omega_{n,a}~\cap~\Omega^{n}_{\delta,b},
$$
where $\mathcal{W}_{a,b}^{n,\delta}$ denotes the set of all vectors $(\theta,t,s)\in\cU_{\delta}\times\RRR^{2}$ satisfying
\begin{align*}
|t -\overline{\psi}(\theta)|\vee |s -\overline{\psi}(\theta^{*})|\leq a/\sqrt{n}\quad\mbox{and}\quad 
|t -s - \overline{\psi}(\theta) + \overline{\psi}(\theta^{*})|\leq b~\delta^{\beta}/\sqrt{n}.
\end{align*}
Then we may conclude 
\begin{align}
&\nonumber
\sup_{\theta\in U_{\delta}}
\big|\MMM_{n}(\theta,\omega) - \MMM_{n}(\theta^{*},\omega) - \psi_{H,\Theta}(\theta) + \psi_{H,\Theta}(\theta^{*})\big|\\ 
&\nonumber 
\leq 
\sup_{(\theta,t,s)\in \mathcal{W}^{n,\delta}_{a,b}}\big|(\MP_{n} - \MP))\big(H(\theta,t,\cdot) - H(\theta^{*},s,\cdot)_{|\omega}\big)\big|\\
&\label{zentrale Aufspaltung 1} 
\quad+ \sup_{(\theta,t,s)\in \mathcal{W}^{n,\delta}_{a,b}}\big|\EEE[H(\theta,t,Z_{1}) - H(\theta^{*},s,Z_{1})] - \psi_{H,\Theta}(\theta) + \psi_{H,\Theta}(\theta^{*})\big|
\end{align}
for $\omega\in\Omega_{n,a}\cap\Omega^{n}_{\delta,b}$. We continue by providing an upper estimation of the second summand in \eqref{zentrale Aufspaltung 1}. 
\begin{lemma}
\label{schweres Lemma}
Let \ref{assump:A1}, \ref{assump:A4}, \ref{assump:A5}, \ref{assump:A6} be fulfilled with $\overline{M}_{1} > 0$ as well as $\beta \in ]0,1]$ from \ref{assump:A5}, and   
$\delta_{1} > 0$, $K_{\delta_{1}} \geq 0$ being as in \ref{assump:A6}. Then for $a, b > 0$ and $\delta\in ]0,\Delta(\Theta)]$
$$
\sup_{(\theta,t,s)\in \mathcal{W}^{n,\delta}_{a,b}}\big|\EEE[H(\theta,t,Z_{1}) - H(\theta^{*},s,Z_{1})] - \psi_{H,\Theta}(\theta) + \psi_{H,\Theta}(\theta^{*})\big|
\leq 
\eins_{]0,\infty[}(\EEE[L_{2}(Z_{1})])~\frac{M_{a,b}~\delta^{\beta}}{\sqrt{n}}
$$
holds if $n\in\NNN$ with $\eins_{]0,\infty[}(\EEE[L_{2}(Z_{1})])\cdot a\leq\sqrt{n}~\delta_{1}$, where $M_{a,b}$ 
denotes the constant defined in \eqref{Mab}. 
\end{lemma}
The proof of Lemma \ref{schweres Lemma} is provided in Subsection \ref{Beweis schweres Lemma}.
\medskip

Let us turn over to the first summand on the right hand side of \eqref{zentrale Aufspaltung 1}. The aim is to apply Theorem \ref{upper estimation} to the function class $\FFF^{n,\delta}_{a,b} := \{H(\theta,t,\cdot) - H(\theta^{*},s,\cdot)\mid (\theta,t,s)\in\mathcal{W}^{n,\delta}_{a,b}\}$ for $a, b, \delta > 0$ and $n\in\NNN$. For this purpose we have to study the covering numbers which are involved in the associated uniform entropy integrals.
\begin{lemma}
\label{uniform covering numbers}
Let \ref{assump:A1}, \ref{assump:A4} and \ref{assump:A5} be satisfied with $L_{1},L_{2}$ denoting the square $\MP^{Z}$-integrable mappings from \ref{assump:A4}, and $ \overline{M}_{1} > 0$ as well as $\beta\in ]0,1]$ being as in assumption \ref{assump:A5}. Furthermore, let $a, b > 0$, $n\in\NNN$, $\delta\in ]0,\Delta(\Theta)]$, and let $\xi^{G}_{\delta}$ stand for the positive envelope of $\FFF^{\Theta}_{\delta}$ as in \ref{assump:A5}. Then the mapping $\xi^{G}_{a,b,\delta}$ introduced in \eqref{neue envelope} is a positive envelope of $\FFF^{n,\delta}_{a,b}$,  
and 
\begin{align*}
N\big(\eta \|\xi^{G}_{a,b,\delta}\|_{\MQ,2},
\FFF^{n,\delta}_{a,b},L^{2}(\MQ)\big)
\leq 
\frac{144~(a/\sqrt{n} + \overline{M}_{1} \delta^{\beta})^{2}}{[(b + \overline{M}_{1})\cdot \delta^{\beta}\cdot\eta]^{2} }~\sup_{\MQ\in\cM_{\textrm{\tiny fin}}}N\big(\eta \|\xi^{G}_{\delta}\|_{\MQ,2}/4,\FFF^{\Theta}_{\delta},L^{2}(\MQ)\big) 
\end{align*}
for $\MQ\in\cM_{\textrm{\tiny fin}}$ and $\eta\in ]0,1[$. In particular 
\begin{align*}
J(\FFF^{n,\delta}_{a,b},\xi^{G}_{a,b,\delta},\varepsilon)
&\leq
4 ~J(\FFF^{\Theta}_{\delta},\xi^{G}_{\delta},\varepsilon/4) + 2~\sqrt{2}\varepsilon~
\sqrt{\ln\left(\frac{12}{\varepsilon}~\Big[\frac{a}{b \sqrt{n} \delta^{\beta}}\vee 1\Big]\right)}
\end{align*}
holds for $\varepsilon\in ]0,1]$.
\end{lemma}
The proof of Lemma \ref{uniform covering numbers} may be found in Subsection \ref{Beweis uniform covering numbers}. 
\smallskip

As a last providing step we want to verify the conditions 1) - 3) in Theorem \ref{upper estimation} for the function classes  
$\FFF^{n,\delta}_{a,b}$. So let $a,b > 0$, $\delta\in ]0,\Delta(\Theta)]$ and $\overline{n}\in\NNN$.
\smallskip

The parameter set $\mathcal{W}^{\overline{n},\delta}_{a,b}$ is a bounded subset of $\RRR^{m+2}$. By Remark \ref{Stetigkeiten} the mapping $\theta\mapsto\EEE[G(\theta,Z_{1})]$ on $\Theta$ is continuous under \ref{assump:A1} and \ref{assump:A2}. Then it may be verified easily that $\mathcal{W}^{\overline{n},\delta}_{a,b}$ is closed w.r.t. the Euclidean metric, and thus compact. 
\par
Next, the mapping $H(\theta,t,\cdot) - H(\theta^{*},s,\cdot)$  is Borel measurable for $(\theta,t,s)\in\mathcal{W}^{\overline{n},\delta}_{a,b}$ due to \ref{measurability H}. This means that $\FFF^{\overline{n},\delta}_{a,b}$ meets condition 1) in Theorem \ref{upper estimation}. 
\par
Condition 2) in Theorem \ref{upper estimation} is already known to be fulfilled with 
the square $\MP^{Z}$-integrable mapping $\xi^{G}_{a,b,\delta}$ from \eqref{neue envelope} due to Lemma \ref{uniform covering numbers} along with \ref{assump:A4} and \ref{assump:A5}. 
\par
Concerning condition 3)  we may find by property \ref{assump:A3}  an at most countable subset $\cC(\cU_{\delta})\subseteq\cU_{\delta}$ and some $(\MP^{Z})^{n}$-null sets $N_{n}^{\cU_{\delta}}$ ($n\in\NNN$) such that 
$$
\inf_{\vartheta\in\cC(\cU_{\delta})}\left\{\EEE[|G(\theta,Z_{1}) - G(\vartheta,Z_{1})|] + \max_{j\in\{1,\ldots,n\}}\big|G(\theta,z_{j}) - G(\vartheta,z_{j})\big|\right\} = 0\
$$
for $\theta\in\cU_{\delta}$, $n\in\NNN$ and $(z_{1},\ldots,z_{n})\in\RRR^{d n}\setminus N_{n}^{\cU_{\delta}}$. Every subset of $\RRR^{2}$ is separable w.r.t. the Euclidean metric so that for each $\vartheta\in\mathcal{C}(\mathcal{U}_{\delta})$ there is some at most countable subset $\mathcal{M}_{\overline{n}}(\vartheta)$ of $\{(t,s)\in\RRR^{2}\mid (\vartheta,t,s)\in \mathcal{W}^{\overline{n},\delta}_{a,b}\}$ which is dense w.r.t. the Euclidean metric. Then $\overline{\mathcal{W}}^{\overline{n},\delta}_{a,b} := \{(\vartheta,t,s)\mid \vartheta\in\mathcal{C}(\mathcal{U}_{\delta}), (t,s)\in\mathcal{M}_{\overline{n}}(\vartheta)\}$ is an at most countable subset of $\mathcal{W}^{\overline{n},\delta}_{a,b}$. 
\smallskip

Our intention is to verify condition 3) in Theorem \ref{upper estimation} for $\FFF^{\overline{n},\delta}_{a,b}$ with at most countable subset $\overline{\mathcal{W}}^{\overline{n},\delta}_{a,b}$ of $\mathcal{W}^{\overline{n},\delta}_{a,b}$ and $(\MP^{Z})^{n}$-null sets $N_{n}^{\cU_{\delta}}$ $(n\in\NNN)$. For this purpose fix $(\overline{\theta},\overline{t},\overline{s})$ from $\mathcal{W}^{\overline{n},\delta}_{a,b}$, $n\in\NNN$ and $(\overline{z}_{1},\ldots,\overline{z}_{n})\in\RRR^{d n}\setminus N_{n}^{\cU_{\delta}}$. Then there exists some sequence $(\vartheta_{l})_{l\in\NNN}$ in $\mathcal{C}(\mathcal{U}_{\delta})$ such that $\big(G(\vartheta_{l},\overline{z}_{j})\big)_{l\in\NNN}$ converges to $G(\overline{\theta},\overline{z}_{j})$ for every $j\in\{1,\ldots,n\}$, and 
$\EEE[|G(\vartheta_{l},Z_{1}) - G(\overline{\theta},Z_{1})|]\to 0$ holds. Choose $t_{l} := \overline{\psi}(\vartheta_{l}) - \overline{\psi}(\overline{\theta}) + \overline{t}$ for $l\in\NNN$. Hence $(\vartheta_{l},t_{l},\overline{s})\in \mathcal{W}^{\overline{n},\delta}_{a,b}$ for $l\in\NNN$, and thus we may select for any $l\in\NNN$ some sequence $\big((t_{lq},s_{lq})\big)_{q\in\NNN}$ in $\mathcal{M}_{\overline{n}}(\vartheta_{l})$ which converges to $(t_{l},\overline{s})$ w.r.t. the Euclidean metric. Note 
that $\{(\vartheta_{l},t_{lq},s_{lq})\mid l,q\in\NNN\}$ is a subset of $\overline{\mathcal{W}}^{\overline{n},\delta}_{a,b}$.

Now, by \ref{assump:A4} along with the choices of the numbers $t_{l},t_{lq}, s_{lq}$ we end up with
\begin{align}
&\nonumber 
\big|[H(\vartheta_{l},t_{lq},z) - H(\theta^{*},s_{lq},z)] - [H(\overline{\theta},\overline{t},z) - H(\theta^{*},\overline{s},z)]\big|^{2}\\
&\label{overall inequality}\leq 
4 L_{1}(z)^{2}~|G(\vartheta_{l},z) - G(\overline{\theta},z)|^{2} + 4~ L_{2}(z)^{2}~[|t_{lq} -\overline{t}|^{2} + |s_{lq} - \overline{s}|^{2}]\quad\mbox{for}~z\in\RRR^{d},
\end{align}
and thus for every $j\in\{1,\ldots,n\}$
\begin{align*}
&
\limsup_{l\to\infty}\limsup_{q\to\infty}\big|[H(\vartheta_{l},t_{lq},\overline{z}_{j}) - H(\theta^{*},s_{lq},\overline{z}_{j})] - [H(\overline{\theta},\overline{t},\overline{z}_{j}) - H(\theta^{*},\overline{s},\overline{z}_{j})]\big|^{2} \nonumber\\
&\leq 
\limsup_{l\to\infty}~ 4 [L_{1}(\overline{z}_{j})^{2} + L_{2}(\overline{z}_{j})^{2}]~[|G(\vartheta_{l},\overline{z}_{j}) - G(\overline{\theta},\overline{z}_{j})|^{2} +|t_{l} -\overline{t}|^{2}] = 0. 
\end{align*}
In the last step we used $\overline{\psi}(\vartheta_{l})\to\overline{\psi}(\overline{\theta})$ which implies $t_{l}\to \overline{t}$. 
\smallskip

Moreover, $L_{1}(Z_{1})~|G(\vartheta_{l},Z_{1}) - G(\overline{\theta},Z_{1})|\to 0$ in 
$\MP$-probability because by assumption $\EEE[|G(\vartheta_{l},Z_{1}) - G(\overline{\theta},Z_{1})|]\to 0$. In addition the sequence $\big(L_{1}(Z_{1})~|G(\vartheta_{l},Z_{1}) - G(\overline{\theta},Z_{1})|\big)_{l\in\NNN}$ is dominated by 
the square integrable random variable $2~L_{1}(Z_{1})~\xi^{G}_{\delta}(Z_{1})$. Then by Vitalis theorem (see \cite[Proposition 21.4]{Bauer2001})
$$
\lim_{l\to\infty}\EEE\big[L_{1}(Z_{1})~|G(\vartheta_{l},Z_{1}) - G(\overline{\theta},Z_{1})|\big] = 0.
$$
Hence in view of \eqref{overall inequality}
\begin{align*}
&
\limsup_{l\to\infty}\limsup_{q\to\infty}\EEE\big[\big|[H(\vartheta_{l},t_{lq},Z_{1}) - H(\theta^{*},s_{lq},Z_{1})] - [H(\overline{\theta},\overline{t},Z_{1}) - H(\theta^{*},\overline{s},Z_{1})]\big|\big]\\
&\leq
\limsup_{l\to\infty}\big\{2 ~\EEE[L_{1}(Z_{1})~|G(\vartheta_{l},Z_{1}) - G(\overline{\theta},Z_{1})|] + 2~\EEE[L_{2}(Z_{1})]~|t_{l} -\overline{t}|\big\} = 0. 
\end{align*}
Now, we have verified condition 3) in Theorem \ref{upper estimation} with at most countable subset $\overline{\mathcal{W}}^{\overline{n},\delta}_{a,b}$ and $(\MP^{Z})^{n}$-nulls set $N_{n}^{\cU_{\delta}}$ ($n\in\NNN$).
\par
Summarizing, $\FFF^{\overline{n},\delta}_{a,b}$ meets all the requirements of Theorem \ref{upper estimation}.
\medskip

Now, we are in the position to show Theorem \ref{Main Theorem}. 
\medskip

\noindent
\textbf{Proof of Theorem \ref{Main Theorem}:}\\
Fix $\delta\in ]0,\Delta(\Theta)]$ and any $n\in\NNN$ satisfying the inequalities $n\geq \overline{M}_{b}^{2} \Delta(\Theta)^{2\beta}/2$ and $n\geq\eins_{]0,\infty[}\big(\EEE[L_{2}(Z_{1})]\big)\cdot a^{2}/\delta_{1}^{1}$ with $\overline{M}_{b}$ as in \eqref{Mb}. Let us also set $r_{\beta,n} := n^{1/{(4 - 2\beta)}}$ and $\mathfrak{l}_{2} := \eins_{]0,\infty[}\big(\EEE[L_{2}(Z_{1})]\big)$. We have already verified that the function class $\FFF^{n,\delta}_{a,b}$ meets all the requirements of Theorem \ref{upper estimation}. In particular, $\sup_{f\in\FFF^{n,\delta}_{a,b}}\big|(\MP_{n} - \MP)(f)\big|$ is a random variable on $\OFP$. Note in addition that $\GGG_{n}^{G,H}(\theta,\cdot) = \sqrt{n}~[\MMM_{n}(\theta,\cdot) - \psi_{H,\Theta}(\theta)]$ holds for $\theta\in\cU_{\delta}$. Then, combining \eqref{zentrale Aufspaltung 1} with Lemma \ref{schweres Lemma}, we obtain 
\begin{align}
&\nonumber 
\MP^{*}\Big(\Big\{\sup_{\theta\in U_{\delta}}
\big|\GGG_{n}^{G,H}(\theta,\cdot) - \GGG_{n}^{G,H}(\theta^{*},\cdot) \big|> \varepsilon\cdot r_{\beta,n}^{-\beta}\Big\}~\cap \Omega_{n,a}\cap\Omega^{n}_{\delta,b}\Big)\\
&\label{Aufspaltung 2}\leq
\MP\Big(\Big\{\sup_{f\in\FFF^{n,\delta}_{a,b}}\big|(\MP_{n} - \MP)(f)\big|> \frac{\varepsilon}{r_{\beta,n}^{2}}- \mathfrak{l}_{2}\cdot\frac{M_{a,b}\cdot\delta^{\beta}}{\sqrt{n}}\Big\}\Big).
\end{align}
Here $M_{a,b}$ denotes the constant defined in \eqref{Mab}. In order to find suitable upper estimations of the right hand side of \eqref{Aufspaltung 2} we want to apply Theorem \ref{upper estimation} to $\FFF^{n,\delta}_{a,b}$.
\smallskip   

Drawing on \ref{assump:A5}, the square $\MP^{Z}$-integrable positive envelope $\xi^{G}_{a,b,\delta}$ of $\FFF^{n,\delta}_{a,b}$ satisfies
\begin{align}
\label{Vor sample size}
\frac{\|\xi^{G}_{a,b,\delta}\|_{\MP^{Z},2}^{2}}{2}
\leq 
\frac{\overline{M}_{b}^{2}\cdot\delta^{2\beta}}{2}
\leq 
\frac{\overline{M}_{b}^{2}\cdot\Delta(\Theta)^{2\beta}}{2}\leq n.
\end{align}
Thus 
by Lemma \ref{uniform covering numbers} along with \ref{assump:A5}, recalling notation \eqref{etantab},
\begin{align}
\label{estimation entropy integral}
16\sqrt{2}~J(\FFF^{n,\delta}_{a,b},\xi^{G}_{a,b,\delta},1/2)
\leq 
\overline{\eta}_{n}(a,b,\delta).
\end{align}
Putting together \eqref{Vor sample size} and \eqref{estimation entropy integral}, the application of Theorem \ref{upper estimation} to $\FFF^{n,\delta}_{a,b}$ yields
$$
\EEE\Big[\sup_{f\in\FFF^{n,\delta}_{a,b}}\big|(\MP_{n} - \MP)(f)\big|\Big]\leq 
\overline{M}_{b}~\overline{\eta}_{n}(a,b,\delta)~\delta^{\beta}/\sqrt{n}.
$$
so that in view of \eqref{Aufspaltung 2} the inequality in statement 1) may be derived via Markov's inequality.
\smallskip 

Let us now turnover to statement 2). 
\smallskip

Additionally, with $\eta_{t,n}(a,b,\delta) := 1 + 2 (t+ 1)~\overline{\eta}_{n}(a,b,\delta)$, we assume 
$$
\delta < \big(\varepsilon/[\mathfrak{l}_{2}\cdot M_{a,b} + \overline{M}_{b}\cdot \eta_{t,n}(a,b,\delta)]\big)^{1/\beta}\cdot n^{-1/(4 - 2\beta)}.
$$ 

By Lemma \ref{uniform covering numbers}  we may observe
\begin{align}
\label{Stern}
1 + 32\sqrt{2} (t + 1)~J(\FFF^{n,\delta}_{a,b},\xi^{G}_{a,b,\delta},1/4)
\leq 
\eta_{t,n}(a,b,\delta).
\end{align}
Thus, 
combining \eqref{Vor sample size} and \eqref{Stern}, we obtain
\begin{align*}
\varepsilon/r_{\beta,n}^{2} - \mathfrak{l}_{2}\cdot M_{a,b}\cdot \delta^{\beta}/\sqrt{n} > \|\xi^{G}_{a,b,\delta}\|_{\MP^{Z},2}~\big[1 + 32\sqrt{2} (t + 1)~J(\FFF^{n,\delta}_{a,b},\xi^{G}_{a,b,\delta},1/4)
\big]/\sqrt{n}.
\end{align*}
Then a further application of Theorem \ref{upper estimation} together with \eqref{Vor sample size} yields
\begin{align*}
&
\MP\Big(\Big\{\sup_{f\in \FFF^{n,\delta}_{a,b}}\big|(\MP_{n} - \MP)(f)\big|> \varepsilon/r_{\beta,n}^{2} - \mathfrak{l}_{2}\cdot M_{a,b} \delta^{\beta}/\sqrt{n}\Big\}\Big)\\
&\leq
\exp\left(-\frac{(\varepsilon - n^{\beta/(4 - 2\beta)}\cdot\mathfrak{l}_{2}\cdot M_{a,b}~\delta^{\beta} )\cdot\mathfrak{g}(t)}{n^{\beta/(4 - 2\beta)}~\|\xi^{G}_{a,b,\delta}\|_{\MP^{Z},2}}\right)
+ \MP\big(\Omega\setminus B_{n}^{\xi^{G}_{a,b,\delta}}\big)\\
&\leq\exp\big(\mathfrak{l}_{2}\cdot M_{a,b}\cdot\mathfrak{g}(t)/\overline{M}_{b}\big) 
\cdot \exp\left(-\frac{\varepsilon \cdot\mathfrak{g}(t)}{n^{\beta/(4 - 2\beta)}~\overline{M}_{b}~\delta^{\beta}}\right)
+ \MP\big(\Omega\setminus B_{n}^{\xi^{G}_{a,b,\delta}}\big).
\end{align*}
This shows statement 2) due to \eqref{Aufspaltung 2} and completes the proof. 
\hfill$\Box$
\subsubsection{Proof of Lemma \ref{schweres Lemma}}
\label{Beweis schweres Lemma}
If $\EEE[L_{2}(Z_{1})] = 0$, then $H(\theta,t,Z_{1})- H\big(\theta,\overline{\psi}(\theta),Z_{1}\big) = 0$ $\MP$-a.s. for every $\theta\in\Theta$. Thus the statement follows immediately. So let us assume that $\EEE[L_{2}(Z_{1})] \not= 0$ which means that $\EEE[L_{2}(Z_{1})]$ is strictly positive.
\par 
Firstly, applying change of variable formula several times 
we may observe by \ref{assump:A6}
\begin{align*}
&
\EEE\big[H(\theta,t,Z_{1})- H\big(\theta,\overline{\psi}(\theta),Z_{1}\big)\big] \\
&=
- \int_{0}^{[t - \overline{\psi}(\theta)]^{-}}\hspace*{-0.75cm}m_{\Theta\times\RRR}\big(\theta,\overline{\psi}(\theta) - u\big)~du + \int_{0}^{[t - \overline{\psi}(\theta)]^{+}}\hspace*{-0.75cm}m_{\Theta\times\RRR}\big(\theta,\overline{\psi}(\theta) + u\big)~du
\end{align*}
for $(\theta,t)\in\Theta\times\RRR$. Then routine calculations yield
\begin{align}
&\nonumber
\big|\EEE\big[H(\theta,t,Z_{1})- H\big(\theta,\overline{\psi}(\theta),Z_{1}\big) - H(\theta^{*},s,Z_{1}) + H\big(\theta^{*},\overline{\psi}(\theta^{*}),Z_{1}\big)\big]\big|\\
&\nonumber\leq
\big|\EEE\big[H(\theta,t,Z_{1})- H\big(\theta,\overline{\psi}(\theta),Z_{1}\big) - H\big(\theta^{*},t - \overline{\psi}(\theta) + \overline{\psi}(\theta^{*}),Z_{1}\big) + H\big(\theta^{*},\overline{\psi}(\theta^{*}),Z_{1}\big)\big]\big|\\
&\nonumber\quad+
\big|\EEE\big[H(\theta^{*},s,Z_{1}) - H\big(\theta^{*},t - \overline{\psi}(\theta) +\overline{\psi}(\theta^{*}),Z_{1}\big)\big]\big|\\
&\nonumber\leq 
\int_{0}^{[t - \overline{\psi}(\theta)]^{-}}\hspace*{-0.75cm}|m_{\Theta\times\RRR}(\theta,\overline{\psi}(\theta) - u)- m_{\Theta\times\RRR}(\theta^{*},\overline{\psi}(\theta^{*}) - u)|~du\\ 
&\nonumber\quad+
\int_{0}^{[t - \overline{\psi}(\theta)]^{+}}\hspace*{-0.75cm}|m_{\Theta\times\RRR}\big(\theta,\overline{\psi}(\theta) + u\big)- m_{\Theta\times\RRR}\big(\theta^{*},\overline{\psi}(\theta^{*}) + u\big)|~du\\
&\label{Auftakt}
\quad + 
\EEE\big[\big|H(\theta^{*},s,Z_{1}) - H\big(\theta^{*},t - \overline{\psi}(\theta) +\overline{\psi}(\theta^{*}),Z_{1}\big)\big|\big]
\end{align}
for $(\theta,t), (\theta^{*},s)\in\Theta\times\RRR$.
\smallskip

By definition of the sets $\mathcal{W}^{n,\delta}_{a,b}$ we may conclude from \ref{assump:A4}
\begin{align}
&\nonumber
\EEE\big[\big|H(\theta^{*},s,Z_{1}) - H\big(\theta^{*},t - \overline{\psi}(\theta) +\overline{\psi}(\theta^{*}),Z_{1}\big)\big|\big]
\\
&\nonumber\leq 
\EEE[L_{2}(Z_{1})]~|t - s - \overline{\psi}(\theta) + \overline{\psi}(\theta^{*})|
\\
&\label{einfache Ungleichung}
\leq \EEE[L_{2}(Z_{1})]~b~\delta^{\beta}/\sqrt{n}
\quad\mbox{for}~\delta\in ]0,\Delta(\Theta)], (\theta,t,s)\in \mathcal{W}^{n,\delta}_{a,b}.
\end{align}
\smallskip

Now, fix any $\delta\in ]0,\Delta(\Theta)]$, 
and choose $n\in\NNN$ with $a/ \sqrt{n}\leq \delta_{1}$. By \ref{assump:A5} 
we know that the chosen positive envelope $\xi^{G}_{\delta}$ of $\FFF^{\Theta}_{\delta}$ 
satisfies the inequality $\|\xi^{G}_{\delta}\|_{\MP^{Z},2}\leq \overline{M}_{1}\cdot\delta^{\beta}$. This implies 
\begin{equation*}
\label{wichtige Ungleichung}
|\EEE[G(\theta,Z_{1})] - \EEE[G(\theta^{*},Z_{1})]|\leq \|\xi^{G}_{\delta}\|_{\MP^{Z},2}\leq \overline{M}_{1}\cdot\delta^{\beta}\quad\mbox{for}~ \theta\in\cU_{\delta}.
\end{equation*} 
Then, invoking assumption \ref{assump:A6} with constants $\delta_{1} > 0, K_{\delta_{1}}\geq 0$,
\begin{align*}
\big|m_{\Theta\times\RRR}\big(\theta,\overline{\psi}(\theta) + u\big) - m_{\Theta\times\RRR}\big(\theta^{*},\overline{\psi}(\theta^{*}) + u\big)\big|
&\nonumber\leq 
K_{\delta_{1}}~\big\|(\theta - \theta^{*}, \overline{\psi}(\theta) - \overline{\psi}(\theta^{*})\big)\big\|_{m+1}\\
&\nonumber\leq 
K_{\delta_{1}}~\big[\|\theta - \theta^{*}\|_{m} + |\overline{\psi}(\theta) - \overline{\psi}(\theta^{*})|\big]\\
&\nonumber\leq 
K_{\delta_{1}}~\big[\delta + \overline{M}_{1}\cdot\delta^{\beta}\big]\\
&
\leq
K_{\delta_{1}}~\big[\Delta(\Theta)^{1-\beta} + \overline{M}_{1}\cdot\big]~\delta^{\beta}
\end{align*} 
for $(\theta,t,u)\in\cU_{\delta}\times\RRR^{2}$ with $|\overline{\psi}(\theta) - t|\leq a/\sqrt{n}\leq\delta_{1}$ and $|u|\leq |\overline{\psi}(\theta) - t|$. 
Thus
\begin{align}
&\nonumber
\int_{0}^{[t - \overline{\psi}(\theta)]^{-}}\hspace*{-0.75cm}|m_{\Theta\times\RRR}(\theta,\overline{\psi}(\theta) - u)- m_{\Theta\times\RRR}(\theta^{*},\overline{\psi}(\theta^{*}) - u)|~du\\ 
&\nonumber\quad+
\int_{0}^{[t - \overline{\psi}(\theta)]^{+}}\hspace*{-0.75cm}|m_{\Theta\times\RRR}\big(\theta,\overline{\psi}(\theta) + u\big)- m_{\Theta\times\RRR}\big(\theta^{*},\overline{\psi}(\theta^{*}) + u\big)|~du\\
&\label{erste Konsequenz}\leq 
|\overline{\psi}(\theta) - t|~K_{\delta_{1}}~\big[\Delta(\Theta)^{1-\beta} + \overline{M}_{1}\cdot\big]~\delta^{\beta}
\leq K_{\delta_{1}}~a~(\Delta(\Theta)^{1-\beta} + \overline{M}_{1})~\delta^{\beta}/\sqrt{n}
\end{align}
for $(\theta,t)\in\cU_{\delta}\times\RRR$ with $|\overline{\psi}(\theta) - t|\leq a/\sqrt{n}$.
Combining \eqref{Auftakt} with \eqref{einfache Ungleichung} and \eqref{erste Konsequenz}, we end up with 
\begin{align*}
&\nonumber
\big|\EEE\big[H(\theta,t,Z_{1})- H\big(\theta,\overline{\psi}(\theta),Z_{1}\big) - H(\theta^{*},s,Z_{1}) + H\big(\theta^{*},\overline{\psi}(\theta^{*}),Z_{1}\big)\big]\big|\\
&\leq 
[K_{\delta_{1}}~a~(\Delta(\Theta)^{1-\beta} + \overline{M}_{1})~+ \EEE[L_{2}]~b
]~\delta^{\beta}/\sqrt{n} 
\end{align*}
for $(\theta,t,s)\in\mathcal{W}^{n,\delta}_{a,b}$.
This completes the proof.
\hfill$\Box$
\subsubsection{Proof of Lemma \ref{uniform covering numbers}}
\label{Beweis uniform covering numbers}
For $(\theta,t,s)\in\mathcal{W}^{n,\delta}_{a,b}$ we may observe by \ref{assump:A5} 
\begin{align*}
|t - s|
\leq
|t - s- \overline{\psi}(\theta) + \overline{\psi}(\theta^{*})| + |\overline{\psi}(\theta) - \overline{\psi}(\theta^{*})|
\leq 
\frac{b~\delta^{\beta}}{\sqrt{n}} + \EEE[\xi^{G}_{\delta}(Z_{1})]
\leq 
(b + \overline{M}_{1})~\delta^{\beta}.
\end{align*}
Hence by \ref{assump:A4} it is easy to verify $\xi^{G}_{a,b,\delta}$ as a positive envelope of $\FFF^{n,\delta}_{a,b}$.
\par
Furthermore by \ref{assump:A5} along with definition of $\mathcal{W}^{n,\delta}_{a,b}$ we obtain for $(\theta,t,s)\in\mathcal{W}^{n,\delta}_{a,b}$
\begin{align*}
&
|t - \overline{\psi}(\theta^{*})|\leq |t - \overline{\psi}(\theta)| + 
|\overline{\psi}(\theta) - \overline{\psi}(\theta^{*})|
\leq 
a/\sqrt{n} + \EEE[\xi^{G}_{\delta}(Z_{1})] 
\leq 
a/\sqrt{n} + \overline{M}_{1}~\delta^{\beta},\\
&
|s - \overline{\psi}(\theta^{*})|\leq a/\sqrt{n}.
\end{align*}
In particular $$
t,s\in\mathcal{I}^{n,\delta}_{a,b} := \big[~\overline{\psi}(\theta^{*}) - a/\sqrt{n} - \overline{M}_{1}~\delta^{\beta},\overline{\psi}(\theta^{*}) + a/\sqrt{n} + \overline{M}_{1}~\delta^{\beta}\big],
$$ 
and thus $\FFF^{n,\delta}_{a,b}\subseteq \widehat{\FFF}^{n,\delta}_{a,b} := \big\{H(\theta,t,\cdot) - H(\theta^{*},s,\cdot)\mid \theta\in\cU_{\delta}, t, s\in \mathcal{I}^{n,\delta}_{a,b}\big\}$.
\smallskip

Fix $\eta \in ]0,1[$ and $\MQ\in\cM_{\textrm{\tiny fin}}$ with support $\textrm{supp}(\MQ)$. We define a new Borel probability measure $\overline{Q}$ on $\RRR^{d}$  by $\overline{\MQ}(A) := \|(L_{1}\vee 1)\cdot \eins_{A}\|_{\MQ,2}^{2}/\|L_{1}\vee 1\|_{\MQ,2}^{2}$. It is absolutely continuous w.r.t. $\MQ$, and thus belongs also to $\cM_{\textrm{\tiny fin}}$. Next, let $\theta,\vartheta\in\cU_{\delta}$ such that the inequality $\|G(\theta,\cdot) - G(\vartheta,\cdot)\|_{\overline{\MQ},2}\leq \eta~\|\xi^{G}_{\delta}\|_{\overline{\MQ},2}/4$ holds, and let $t_{1}, t_{2}, s_{1}, s_{2}\in \mathcal{I}^{n,\delta}_{a,b} $ satisfy the inequalities $|t_{1} - t_{2}|, |s_{1} - s_{2}| \leq \eta~ [b + \overline{M}_{1}]~\delta^{\beta}/6$. Then by \ref{assump:A4} we may observe
\begin{align*}
&
\big\|\big[H(\theta,t_{1},\cdot) - H(\theta^{*},s_{1},\cdot)\big] - \big[H(\vartheta,t_{2},\cdot) - H(\theta^{*},s_{2},\cdot)\big]\big\|_{\MQ,2}^{2}\\
&\leq 
4 \big\|H(\theta,t_{1},\cdot) - H(\vartheta,t_{2},\cdot) \big\|_{\MQ,2}^{2}  
+
4 \big\|H(\theta^{*},s_{1},\cdot) - H(\theta^{*},s_{2}\cdot)\big\|_{\MQ,2}^{2}\\
&= 
4 \sum_{z\in\textrm{supp}(\MQ)}\big\{L_{1}(z)^{2}~|G(\theta,z) -  G(\vartheta,z)|^{2} + L_{2}(z)^{2}\big(|t_{1} - t_{2}|^{2} + |s_{1} - s_{2}|^{2}\big)~\big\}~\MQ(\{z\})\\
&\leq 
4~\|L_{1}\vee 1\|_{\MQ,2}^{2}~\big\|G(\theta,\cdot) -  G(\vartheta,\cdot)\big\|_{\overline{\MQ},2}^{2} +  4 \|L_{2}\|_{\MQ,2}^{2}~\big(|t_{1} - t_{2}|^{2} + |s_{1} - s_{2}|^{2}\big)\\
&\leq 
\eta^{2}~ \|L_{1}\vee 1\|_{\MQ,2}^{2}~
\|\xi^{G}_{\delta}\|
_{\overline{\MQ},2}^{2}/4 
+ 8\eta^{2}~ \|L_{2}\|_{\MQ,2}^{2}~[b + \overline{M}_{1}]^{2}~\delta^{2\beta}/36\\
&= 
\eta^{2}~ 
\|(L_{1}\vee 1)\cdot \xi^{G}_{\delta}\|
_{\MQ,2}^{2}/4 
+ 8\eta^{2}~ \|L_{2}\|_{\MQ,2}^{2}~[b + \overline{M}_{1}]^{2}~\delta^{2\beta}/36
\leq 
\eta^{2}~\|\xi^{G}_{a,b,\delta}\|_{\MQ,2}^{2}/4.
\end{align*}
Hence, putting all together, we may conclude 
\begin{align}
&\nonumber
N
\big(\eta 
\|\xi^{G}_{a,b,\delta}\|_{\MQ,2},
\FFF^{n,\delta}_{a,b},L^{2}(\MQ)\big)\leq 
N
\big(\eta 
\|\xi^{G}_{a,b,\delta}\|_{\MQ,2}/2,
\widehat{\FFF}^{n,\delta}_{a,b},L^{2}(\MQ)\big)\\
&\nonumber\leq
\sup_{\MQ\in\cM_{\textrm{\tiny fin}}} N\big(\eta 
\|\xi^{G}_{\delta}\|_{\MQ,2}/4,
\FFF^{\Theta}_{\delta},L^{2}(\MQ)
\big)\cdot
N\big(\eta ~[b + \overline{M}_{1}]~\delta^{\beta}/6,\mathcal{I}^{n,\delta}_{a,b},|\cdot|\big)^{2},
\end{align}
where 
$N\big(\eta ~[b + M]~\delta^{\beta}/6,\mathcal{I}^{n,\delta}_{a,b},|\cdot|\big)$ stands for the minimal number to cover $\mathcal{I}^{n,\delta}_{a,b}$ by intervals of the form $[t - \eta~[b + \overline{M}_{1}]~\delta^{\beta}/6,t + \eta~ [b + \overline{M}_{1}]~\delta^{\beta}/6]$ with $t\in \mathcal{I}^{n,\delta}_{a,b}$. It may be easily checked that 
\begin{align*}
N\big(\eta ~[b + \overline{M}_{1}]~\delta^{\beta}/6,\mathcal{I}^{n,\delta}_{a,b},|\cdot|\big)
&\leq \frac{6~ (\sup\mathcal{I}^{n,\delta}_{a,b} - \inf\mathcal{I}^{n,\delta}_{a,b})}{\eta~[b + \overline{M}_{1}]~\delta^{\beta}}
= 
\frac{12~ (a/\sqrt{n} + \overline{M}_{1}\delta^{\beta})}{\eta~[b + \overline{M}_{1}]~\delta^{\beta}}
\end{align*} 
holds. Now, the upper estimate of the covering number $N\big(\eta \|\xi^{G}_{a,b,\delta}\|_{\MQ,2},\FFF^{n,\delta}_{a,b},L^{2}(\MQ)\big)$ that is claimed in the statement may be derived immediately. 
\medskip

For the remaining part of the proof fix $\varepsilon\in ]0, 1]$. Note that $(x + s)/(y + s)\leq (x/y)\vee 1$ holds for $s, x, y > 0$. Hence 
\begin{align*}
\frac{a/\sqrt{n} + M_{\delta_{1}}~\delta^{\beta}}{ (b + M_{\delta_{1}})~\delta^{\beta}} 
= 
\frac{a/(\sqrt{n}~\delta^{\beta}) + M_{\delta_{1}}}{ b + M_{\delta_{1}}} \leq \frac{a}{b~\sqrt{n}~\delta^{\beta}}\vee 1.
\end{align*}
Using the change of variable formula, and invoking  \eqref{Integralabschaetzung}, we end up with
\begin{align*}
&
J(\FFF^{n,\delta}_{a,b},\xi^{G}_{a,b,\delta},\varepsilon)\\
&\leq
\int_{0}^{\varepsilon}\sup_{\MQ\in \cM_{\textrm{\tiny fin}}}\sqrt{\ln\big(2 N\big(\eta~\|\xi^{G}_{\delta}\|_{\MQ,2}/4,\FFF^{\Theta}_{\delta},L^{2}(\MQ)\big)\big) }~d\eta
\\
&\quad +
\int_{0}^{\varepsilon}\sqrt{2\ln\Big(12 \big([a/(b~\sqrt{n}~\delta^{\beta})]\vee 1\big)/\eta\Big)}~d\eta\\
&\leq
4 J(\FFF^{\Theta}_{\delta},\xi^{G}_{\delta},\varepsilon/4) + 
\varepsilon \int_{0}^{1}\sqrt{2\ln\Big(12 \big([a/(b~\sqrt{n}~\delta^{\beta})]\vee 1 \big)/[\varepsilon~ u]\Big)}~du\\
&\leq 
4 J(\FFF^{\Theta}_{\delta},\xi^{G}_{\delta},\varepsilon/4) + 2~\sqrt{2}~\varepsilon~\sqrt{\ln\Big(12 \big([a/(b~\sqrt{n}~\delta^{\beta})]\vee 1] \big)/\varepsilon\Big)}.
\end{align*}
This completes the proof.
\hfill$\Box$

\subsection{Proof of Theorem  \ref{simpliest version rough} and Theorem \ref{simpliest version} }
\label{Beweis Adaption}
We start with recalling the result by Proposition \ref{variational inequalities}
\begin{align}
&\nonumber 
\MP\big(\big\{n^{1/(4- 2\beta)}~\|\widehat{\theta}_{n} - \theta^{*}\|_{m} > \varepsilon\cdot \big\}\cap\Omega_{n,a}\big)\\
&\label{variational inequality 1}\leq 
\sum_{k= K_{\varepsilon} + 1\atop \delta_{n (k-1)}\leq\Delta(\Theta)}^{\infty}\hspace*{-0.25cm}\MP^{*}\Big(\Big\{\sup_{\theta\in U_{\Delta(\Theta)\wedge\delta_{nk}}}\big|\GGG_{n}^{G,H}(\theta,\cdot) - \GGG_{n}^{G,H}(\theta^{*},\cdot) \big| > 
\frac{M_{1}~2^{2 (k - 1)}}{n^{\beta/(4 - 2\beta)}} 
\Big\}\cap\Omega_{n,a}\Big)
\end{align}
for $n\in\NNN$.
\smallskip

Now, let $n\in\NNN$ with $n\geq n(a,b,\beta)$, and consider 
$k\in\mathbb{Z}$ such that $k\geq K_{\varepsilon} + 1$ as well as $\delta_{n (k-1)}\leq\Delta(\Theta)$.
\par
First of all
\begin{align*}
\big[\ln\big(a/[\sqrt{2}^{k}~b]\big) - \ln\big(\sqrt{n}~[\delta_{nk}\wedge\Delta(\Theta)]^{\beta}\big)\big]^{+}
&\leq 
\big[\ln(a/b) - \ln(2^{k-1})/2 -\ln\big(\sqrt{n}~\delta_{n (k-1)}^{\beta}\big)\big]^{+}\\
&\leq 
\big[\ln(a/b) - (\beta + 1/2)~\ln(2^{k-1}) \big]^{+}\\
&\leq 
\big[\ln(a/b) - (\beta + 1/2)~\ln(\varepsilon/2) \big]^{+}.
\end{align*}
Then by $a~2^{\beta + 1/2}\leq b~\varepsilon^{\beta + 1/2}$
\begin{align}
\label{zweite Gleichung}
\overline{\eta}_{n}\big(a,\sqrt{2}^{k}~b,\delta_{nk}\wedge\Delta(\Theta)\big)\leq\widehat{\eta},
\end{align}
where we used notations \eqref{etantab},\eqref{etadachepsilon}. Furthermore
\begin{align}
\label{erste Gleichung}
\overline{M}_{\sqrt{2}^{k} b}\leq \sqrt{2}^{k^{+}}~\overline{M}_{b}\quad\mbox{and}\quad
M_{a,\sqrt{2}^{k} b}\leq \sqrt{2}^{k^{+}}~M_{a,b}
 \end{align}
In particular by $n\geq n(a,b,\beta)$
\begin{align}
\frac{\overline{M}_{\sqrt{2}^{k} b}^{2}~\Delta(\Theta)^{2\beta}}{2}
&\nonumber\leq 
\frac{\overline{M}_{b}^{2}~\Delta(\Theta)^{2\beta}}{2}
\vee [n^{\frac{1}{4-2\beta}}~\overline{M}_{b}^{2}~\Delta(\Theta)^{2\beta + 1}]\\ 
&\label{Ungleichung Stichprobenumfang}\leq n(a,b,\beta)^{\frac{3-2\beta}{4-2\beta}}~n^{\frac{1}{4- 2\beta}}\leq n.
\end{align}
Then by definition of $n(a,b,\beta)$ we may apply statement 1) in Theorem \ref{Main Theorem} which together with \eqref{erste Gleichung} and \eqref{zweite Gleichung} yields
\begin{align}
&
\MP\Big(\Big\{\sup_{\theta\in U_{\delta_{nk}\wedge\Delta(\Theta)}}\big|\GGG_{n}^{G,H}(\theta,\cdot) - \GGG_{n}^{G,H}(\theta^{*},\cdot) \big| > 
\frac{M_{1}~2^{2 (k - 1)}}{n^{\beta/(4 - 2\beta)}} 
\Big\}\cap \Omega_{n,a}\cap\Omega^{n}_{\delta_{nk}\wedge\Delta(\Theta),\sqrt{2}^{k}~b}\Big)
\nonumber\\
&\leq 
\frac{4~\big[\overline{M}_{\sqrt{2}^{k} b}~\overline{\eta}_{n}\big(a,\sqrt{2}^{k} b,\delta_{nk}\wedge\Delta(\Theta)\big) + \eins_{]0,\infty[}(\EEE[L_{2}(Z_{1})])\cdot M_{a,\sqrt{2}^{k} b}\big]}{M_{1}~2^{k (2 - \beta)}}\nonumber\\
&\leq 
\frac{4~\big[\overline{M}_{b}~\widehat{\eta} + \eins_{]0,\infty[}(\EEE[L_{2}(Z_{1})])\cdot M_{a,b}\big]}{M_{1}~2^{k (2-\beta) - k^{+}/2}}.\label{vierte Gleichung}
\end{align}
Moreover, invoking well known formulas of geometric series, we obtain 
\begin{align}
\sum_{k=K_{\varepsilon} + 1}^{\infty}\frac{1}{2^{k(2 - \beta) - k^{+}/2}}\leq \frac{1}{2^{K_{\varepsilon}(2-\beta) - K_{\varepsilon}^{+}/2}~(2^{3/2 - \beta} - 1)}\leq\frac{5}{2}~\left(\frac{2}{\varepsilon}\right)^{3/2-\beta}~\sqrt{\frac{2}{\varepsilon\wedge 2}}.\label{sechste Gleichung}
\end{align}
Combining \eqref{variational inequality 1} with \eqref{vierte Gleichung} and \eqref{sechste Gleichung}, we may finish the proof of Theorem \ref{simpliest version rough} easily.
\bigskip

Let us turn over to Theorem \ref{simpliest version}. Let $n\in\NNN$ with $n\geq \eta(a,b,\beta)$, and let $t > 0$ such that  
$$
2^{K_{\varepsilon} (2 - \beta) - K_{\varepsilon}^{+}/2} > 4 ~\big[M_{a,b}~\eins_{]0,\infty[}\big(\EEE[L_{2}(Z_{1})]\big) + \overline{M}_{b}~\big(1 + 2 (t + 1)~\widehat{\eta}\big)\big]/M_{1}.
$$ 
Furthermore let $k\in\mathbb{Z}$ with $k\geq K_{\varepsilon} + 1$ and $\delta_{n (k-1)}\leq\Delta(\Theta)$. Then firstly by \eqref{Ungleichung Stichprobenumfang} the sample size $n$ satisfies
$n\geq\max\big\{\overline{M}_{\sqrt{2}^{k}~b}^{2}~\Delta(\Theta)^{2\beta}/2,\eins_{]0,\infty[}\big(\EEE[L_{2}(Z_{1})]\big)~a^{2}/2\big\}$, and in view of \eqref{erste Gleichung} along with \eqref{zweite Gleichung}
\begin{align*}
&
\frac{M_{1}~2^{2(k-1)}}{M_{a,\sqrt{2}^{k} b}~\eins_{]0,\infty[}\big(\EEE[L_{2}(Z_{1})]\big) + \overline{M}_{\sqrt{2}^{k} b}\big[1 + 2 (t + 1)~\overline{\eta}_{n}\big(a,\sqrt{2}^{k} b, \delta_{nk}\wedge\Delta(\Theta)\big)\big]}
\\
&
\geq 
2^{k\beta}~\frac{M_{1}~2^{K_{\varepsilon} (2 - \beta) - K_{\varepsilon}^{+}/2}}{4~\big[M_{a,b} ~\eins_{]0,\infty[}\big(\EEE[L_{2}(Z_{1})]\big) + \overline{M}_{b}\big(1 + 2 (t + 1)~\widehat{\eta}\big)\big]} > 2^{k\beta}
\end{align*}
Hence by statement 2) in Theorem \ref{Main Theorem} along with \eqref{erste Gleichung}
\begin{align}
&
\MP\Big(\Big\{\sup_{\theta\in U_{\delta_{nk}\wedge\Delta(\Theta)}}\big|\GGG_{n}^{G,H}(\theta,\cdot) - \GGG_{n}^{G,H}(\theta^{*},\cdot) \big| > 
\frac{M_{1}~2^{2 (k - 1)}}{n^{\beta/(4 - 2\beta)}} 
\Big\}\cap \Omega_{n,a}\cap\Omega^{n}_{\delta_{nk}\wedge\Delta(\Theta),\sqrt{2}^{k}~b}\Big)
\nonumber\\
&\leq 
\exp\left(\eins_{]0,\infty[}\big(\EEE[L_{2}(Z_{1})]\big)~\frac{M_{a,\sqrt{2}^{k} b}~\mathfrak{g}(t)}{\overline{M}_{\sqrt{2}^{k} b}}\right)~\exp\left(-\frac{2^{k(2-\beta)}~M_{1}~\mathfrak{g}(t)}{4~\overline{M}_{\sqrt{2}^{k}~b}}\right)\nonumber\\
&\quad 
+ \MP\big(\Omega\setminus B_{n}^{\xi^{G}_{a,\sqrt{2}^{k} b,(\delta_{nk}\wedge\delta)}}\big)\nonumber\\
&\leq 
\exp\left(\eins_{]0,\infty[}\big(\EEE[L_{2}(Z_{1})]\big)~\frac{\sqrt{2}^{(K_{\varepsilon} + 1)^{-}}M_{a,b}~\mathfrak{g}(t)}{b~\|L_{2}\|_{\MP^{Z},2}}\right)~\exp\left(-\frac{2^{k(2-\beta) - k^{+}/2}~M_{1}~\mathfrak{g}(t)}{4~\overline{M}_{b}}\right)\nonumber\\
&\quad 
+ \MP\big(\Omega\setminus B_{n}^{\xi^{G}_{a,\sqrt{2}^{k} b,[\delta_{nk}\wedge\Delta(\Theta)]}}\big)\label{fuenfte Gleichung}
\end{align}
Invoking the change of variable formula, we further obtain
\begin{align*}
&\exp\left(
-\frac{2^{k(2-\beta) - k^{+}/2}~M_{1}~\mathfrak{g}(t)}{4~\overline{M}_{b}}
\right)\\ 
&\leq 
\int_{k-1}^{k}\exp\left(-\frac{2^{u (2-\beta) - u^{+}/2}~M_{1}~\mathfrak{g}(t)}{4~\overline{M}_{b}}\right)~du\\
&= 
\bcswitch
\int_{2^{(k-1) (2-\beta)}}^{2^{k (2-\beta)} } \exp\left(-\frac{y~M_{1}~\mathfrak{g}(t)}{4~\overline{M}_{b}}\right)~\frac{1}{(2-\beta)\ln(2) y}~dy
&k\leq 0\\[0.5cm]
\int_{2^{(k-1) (3/2-\beta)}}^{2^{k (3/2-\beta) }} \exp\left(-\frac{y~M_{1}~\mathfrak{g}(t)}{4~\overline{M}_{b}}\right)~\frac{1}{(3/2-\beta)\ln(2) y}~dy&k > 0
\ecswitch
\\
&\leq 
\bcswitch
\frac{1}{(2-\beta)\ln(2) 2^{K_{\varepsilon}(2-\beta)}}\int_{2^{(k-1) (2-\beta)}}^{2^{k (2-\beta)} } \exp\left(-\frac{y~M_{1}~\mathfrak{g}(t)}{4~\overline{M}_{b}}\right)~dy&k\leq 0\\[0.5cm]
\frac{2}{(3 - 2\beta)\ln(2) 2^{K_{\varepsilon}(3/2-\beta)}}~\int_{2^{(k-1) (3/2-\beta)}}^{2^{k (3/2-\beta) }} \exp\left(-\frac{y~M_{1}~\mathfrak{g}(t)}{4~\overline{M}_{b}}\right)~dy&k > 0 
\ecswitch
\end{align*}
Note further that
$$
\frac{1}{(2-\beta)\ln(2) 2^{K_{\varepsilon}(2-\beta)}}\vee \frac{2}{(3 - 2\beta)\ln(2) 2^{K_{\varepsilon}(3/2-\beta)}}\leq \frac{2}{(3 - 2\beta)\ln(2) 2^{K_{\varepsilon}(2-\beta) - K_{\varepsilon}^{+}/2}}
$$
holds.  
Hence
\begin{align}
&\sum_{k = K_{\varepsilon} +1}^{\infty}\exp\left(
-\frac{2^{k(2-\beta) - k^{+}/2}~M_{1}~\mathfrak{g}(t)}{4~\overline{M}_{b}}
\right)\nonumber\\
&\leq 
\frac{2}{(3 - 2\beta)\ln(2) 2^{K_{\varepsilon}(2-\beta) - K_{\varepsilon}^{+}/2}}
~\int_{2^{K_{\varepsilon}(2 - \beta) - K_{\varepsilon}^{+}/2}}^{\infty}\exp\left(-\frac{y~M_{1}~\mathfrak{g}(t)}{4~\overline{M}_{b}}\right)~dy\nonumber\\
&= 
\frac{4~\overline{M}_{b}}{(3/2 - \beta)\ln(2)~M_{1}~\mathfrak{g}(t)~2^{K_{\varepsilon}(2-\beta) - K_{\varepsilon}^{+}/2}}~\exp\left(-\frac{2^{K_{\varepsilon}(2 - \beta) - K_{\varepsilon}^{+}/2}~M_{1}~\mathfrak{g}(t)}{4~\overline{M}_{b}}\right)\label{siebte Gleichung}
\end{align}
Putting \eqref{variational inequality 1}, \eqref{fuenfte Gleichung} and \eqref{siebte Gleichung} together we may complete the proof of Theorem \ref{simpliest version} easily.
\hfill$\Box$
\subsection{Proof of Remark \ref{erste Restwahrscheinlichkeiten}}
\label{Restwahrscheinlichkeiten}
Firstly, $\Omega_{n,a}, \Omega^{n}_{\delta,b}\in\cF$ for $a,b > 0$, $n\in\NNN$ and $\delta\in ]0,\Delta(\Theta)]$ due to Lemma \ref{messbare Hilfsereignisse}. Secondly, we may observe that the function classes $\FFF^{\Theta}$ and $\FFF^{\Theta}_{\delta}$ meet the requirements of Theorem \ref{upper estimation} due to assumptions \ref{assump:A1}, \ref{assump:A2}, \ref{assump:A3} and \ref{assump:A5}. Thirdly, for $\lambda, \varepsilon > 0$ and any function class $\FFF$ with positive envelope $C_{\FFF}$, the inequality  $J(\FFF,C_{\FFF},\lambda\cdot\varepsilon)\leq \lambda\cdot J(\FFF,C_{\FFF},\varepsilon)$ is valid (see \cite[Lemma 3.5.3]{GineNickl2016}). Then $J(\FFF^{\Theta},\xi^{G},1/4)\vee  
J(\FFF^{\Theta}_{\delta},\xi^{G}_{\delta},1/4)\leq 2~\overline{M}^{1}$ holds due to \ref{assump:A5}. Hence we may apply directly Theorem \ref{upper estimation} to the function classes $\FFF^{\Theta}$ and $\FFF^{\Theta}_{\delta}$ to conclude statements 4), 5). Moreover, $J(\FFF^{\Theta},\xi^{G},1/2)\leq  4~\overline{M}^{1}$ and  
$J(\FFF^{\Theta}_{\delta},\xi^{G}_{\delta},1/2)\leq 4~\overline{M}^{1}$ by \ref{assump:A5}. Therefore, statements 1), 2) follow immediately from Theorem \ref{upper estimation}, invoking Markov's inequality along with assumption \ref{assump:A5}. 
\smallskip

Concerning the statement 3) we obtain by statement 2) 
\begin{align*}
\sum_{k= K_{\varepsilon} + 1\atop \delta_{n (k-1)}\leq\Delta(\Theta)}^{\infty}\hspace*{-0.35cm}\MP\big(\Omega\setminus\Omega^{n}_{\delta_{nk}\wedge\Delta(\Theta),\sqrt{2}^{k} b}\big)
\leq 
\sum_{k= K_{\varepsilon} + 1}^{\infty}\frac{64~\sqrt{2}~\overline{M}_{1}~\overline{M}^{1}} 
{\sqrt{2}^{k}~b} 
= 
\frac{1}{\sqrt{2}^{K_{\varepsilon}}}~\frac{64~\sqrt{2}~\overline{M}_{1}~\overline{M}^{1}} 
{(\sqrt{2}- 1)~b}. 
\end{align*}
Then the statement 3) follows immediately from $K_{\varepsilon}\geq \ln(\varepsilon)/\ln(2) - 1$. For statement 6) let $n\in\NNN$ with $n\geq \overline{M}_{1}^{2}~\Delta(\Theta)^{2\beta}$. Then statement 5) implies
\begin{align*}
\sum_{k= K_{\varepsilon} + 1\atop \delta_{n (k-1)}\leq\Delta(\Theta)}^{\infty}\hspace*{-0.35cm}\MP\big(\Omega\setminus\Omega^{n}_{\delta_{nk}\wedge\Delta(\Theta),\sqrt{2}^{k} b}\big)
\leq 
\sum_{k= K_{\varepsilon} + 1}^{\infty}\exp\left(-\frac{\sqrt{2}^{k}~b~\mathfrak{g}(t)}{\overline{M}_{1}}\right) + \sum_{k= K_{\varepsilon} + 1\atop \delta_{n (k-1)}\leq\delta}^{\infty}\hspace*{-0.35cm}\MP\big(\Omega\setminus B_{n}^{\xi^{G}_{\delta_{nk}\wedge\Delta(\Theta)}}\big)
\end{align*}
for $t > 0$ such that $b > \overline{M}_{1}~[1 + 64\sqrt{2}~(t + 1)~\overline{M}^{1}]/\sqrt{\varepsilon}$. Moreover, invoking the change of variable formula
\begin{align*}
\sum_{k= K_{\varepsilon} + 1}^{\infty}\exp\left(-\frac{\sqrt{2}^{k}~b~\mathfrak{g}(t)}{\overline{M}_{1}}\right) 
&\leq 
\sum_{k= K_{\varepsilon} + 1}^{\infty}\int_{k-1}^{k}\exp\left(-\frac{\sqrt{2}^{u}~b~\mathfrak{g}(t)}{\overline{M}_{1}}\right)~du\\
&\leq 
\sum_{k= K_{\varepsilon} + 1}^{\infty}\int_{\sqrt{2}^{k-1}}^{\sqrt{2}^{k}}\exp\left(-\frac{y~b~\mathfrak{g}(t)}{\overline{M}_{1}}\right)~\frac{2}{y~\ln(2)}~dy\\
&\leq 
\frac{2}{\sqrt{2}^{K_{\varepsilon}}~\ln(2)}~\int_{\sqrt{2}^{K_{\varepsilon}}}^{\infty}\exp\left(-\frac{y~b~\mathfrak{g}(t)}{\overline{M}_{1}}\right)~dy\\ 
&\leq 
\frac{2~\overline{M}_{1}}{\sqrt{2}^{K_{\varepsilon}}~\ln(2)~b~\mathfrak{g}(t)}~\exp\left(-\frac{\sqrt{2}^{K_{\varepsilon}}~b~\mathfrak{g}(t)}{\overline{M}_{1}}\right)
\end{align*}
Now, in view of $\sqrt{2}^{K_{\varepsilon}}\geq \sqrt{\varepsilon/2}$ statement 6) may be derived easily. This completes the proof. 
\hfill$\Box$
\subsection{Proof of Proposition \ref{important covering numbers}}
\label{Beweis von Proposition important covering numbers}
By representation (PH) we have with 
\begin{align}
&\nonumber 
G(\theta,z) - G(\theta^{*},z)\\ 
&\label{Eingangsrepraesentation}=
\sum_{i=1}^{r}\Big([G^{i}(\theta,z) - G^{i}(\theta^{*},z)]\cdot\eins_{B_{i}(\theta)}(z) + G^{i}(\theta^{*},z)\cdot [\eins_{B_{i}(\theta)}(z) - \eins_{B_{i}(\theta^{*})}(z)]\Big)
\end{align}
for $\theta\in\Theta$ and $z\in\RRR^{d}$. 
The processes $\big(G^{i}(\theta,\cdot)\big)_{\theta\in\Theta}$ have H\"older continuous paths according to representation (PH). Hence
by triangle inequality
\begin{align*}
|G(\theta,z) - G(\theta^{*},z)| 
\leq 
\sum_{i=1}^{r}\big(C^{i}(z)~\|\theta - \theta^{*}\|_{m}^{\beta_{i}} + |G^{i}(\theta^{*},z)|\cdot \eins_{B_{i}(\theta)\Delta B_{i}(\theta^{*})}(z)\big)
\end{align*}
for $\theta\in\Theta$ and $z\in\RRR^{d}$. Then it may be concluded easily that $\xi^{G}_{\delta}$ is a positive envelope of $\FFF^{\Theta}_{\delta}$ for $\delta > 0$. 
\smallskip

In the next step we want to derive for $\MQ\in\cM_{\textrm{\tiny fin}}$ upper estimates of the covering numbers $N\big(\eta~\|\xi^{G}_{\delta}\|_{\MQ,2},\FFF^{\Theta}_{\delta},L^{2}(\MQ)\big)$ based on representation \eqref{Eingangsrepraesentation}. For this purpose let us introduce the functions $\overline{G}^{1},\ldots,\overline{G}^{r}, \overline{G}:\Theta\times\RRR^{d}\rightarrow\RRR$ via
\begin{align*}
\overline{G}^{i}(\theta,z) := G^{i}(\theta,z) - G^{i}(\theta^{*})~\mbox{for}~i\in\{1,\ldots,r\}\quad\mbox{and}\quad \overline{G}(\theta,z) = \sum_{i=1}^{r}\overline{G}^{i}(\theta,z)~\eins_{B_{i}(\theta)}(z).
\end{align*}
Note that $\overline{G}^{i}(\theta^{*},\cdot) \equiv 0$ is square $\MP^{Z}$-integrable and 
$|\overline{G}^{i}(\theta,z) - \overline{G}^{i}(\vartheta,z)|\leq C^{i}(z) \|\theta - \vartheta\|_{m}^{\beta_{i}}$ for $\theta, \vartheta\in\Theta$ as well as $i\in\{1,\ldots,r\}$. Then we may apply directly from \cite{Kraetschmer2023a} Proposition 2.8 together with formulas (5.12) and (5.13) in its proof to the function classes $\overline{F}^{\Theta}_{\delta} := \{\overline{G}(\theta,\cdot)\mid\theta\in\cU_{\delta}\}$ ($\delta > 0$). Hence for any $\delta > 0$ a positive square $\MP^{Z}$-integrable envelope $C_{\overline{F}^{\Theta}_{\delta}}$ of $\overline{F}^{\Theta}_{\delta}$ is defined by $C_{\overline{F}^{\Theta}_{\delta}} = 2 \sum_{i=1}^{r}\delta^{\beta_{i}}~C^{i}(z)$, and
\begin{align}
&\nonumber
\sup_{\MQ\in \cM_{\textrm{\tiny fin}}}N\big(\eta~\|C_{\overline{\FFF}^{\Theta}_{\delta}}\|_{\MQ,2},\overline{\FFF}^{\Theta}_{\delta},L^{2}(\MQ)\big)\\
&\label{Abschaetzung1}\leq 
\prod_{i=1}^{r}9^{m}~16^{(m+2) s_{i}}~e^{1 + (m + 2)~s_{i}}~([m + 2]~s_{i} + 1)~\big(4 /\eta\big)^{2 (m + 2) s_{i} + m/\beta_{i}}
\end{align}
for $\eta\in ]0,r[$.
\smallskip

Next let us introduce the auxiliary function classes 
$$
\FFF^{\Theta,i}_{\delta} := \big\{G^{i}(\theta^{*},\cdot)\cdot [\eins_{B_{i}(\theta)} - \eins_{B_{i}(\theta^{*})}]\mid\theta\in\cU_{\delta}\big\}\quad (i\in\{1,\ldots,r\}, \delta > 0).
$$
Concerning upper estimations for the covering numbers of these classes we need some further preparation from the theory of empirical process theory. 
To recall, define for a collection $\cB$ of subsets of $\RRR^{d}$, and $z_{1},\dots,z_{n}\in\RRR^{d}$
$$
\Delta_{n}(\cB,z_{1},\dots,z_{n})~:=~\mbox{cardinality of}~\left\{B~\cap~\{z_{1},\dots,z_{n}\}\mid B\in\cB\right\}.
$$
Then 
$$
V(\cB)~:=~\inf~\Big\{n\in\NNN\mid \max_{z_{1},\dots,z_{n}\in\RRR^{d}}\Delta_{n}(\cB,z_{1},\dots,z_{n}) < 2^{n}\Big\}\quad(\inf\emptyset~:=~\infty)
$$
is known as the \textit{index} of $\cB$ (see \cite{vanderVaartWellner1996}, p. 135). In case of finite index, $\cB$ is known as a so called \textit{VC-class} (see \cite{vanderVaartWellner1996}, p. 135). The concept of VC-classes may be carried over from sets to functions in the following way. A set $\FFF$ of Borel measurable real valued functions on $\RRR^{d}$ is defined to be a \textit{VC-subgraph class} or a \textit{VC-class} if 
the corresponding collection $\big\{\{(z,t)\in\RRR^{d}\times\RRR\mid h(z) > t\}~\mid~h\in\FFF\big\}$ of subgraphs is a {\textit VC-class} (\cite{vanderVaartWellner1996}, p. 141). Its {\textit VC-index} $V(\FFF)$ coincides with the index of the subgraphs. The significance of VC-subgraph classes stems from the fact that for every VC-subgraph class $\FFF$ and any $\MP^{Z}$-integrable positive envelope $C_{\FFF}$ of $\FFF$
\begin{equation}
\label{VC-subgraph covering numbers}
\sup_{\MQ\in\cM_{\textrm{\tiny fin}}}N\big(\eta \|C_{\FFF}\|_{\MQ,2},\FFF,L^{2}(\MQ)\big)\leq e~V(\FFF)~\big(4 e^{1/2}/\eta\big)^{2 [V(\FFF) - 1]}\quad\mbox{if}~\eta\in ]0,1[
\end{equation} 
(see e.g. \cite[Corollary 5.3]{Kraetschmer2023a}). 
A simple example of VC-subgraph classes is provided by function classes whose linear hulls have finite dimensions.
\begin{example}
\label{finite dimension} 
Let $\FFF$ be a class of real-valued Borel measurable mappings on $\RRR^{d}$ such that the dimension $D$ of its linear hull is finite. 
Then, by Lemma 2.6.15 of \cite{vanderVaartWellner1996}, the set $\FFF$ is a VC-subgraph class with VC-index $V(\FFF)\leq D + 2$. Hence, combining 
\eqref{VC-subgraph covering numbers} with \eqref{Integralabschaetzung} and Lemma  3.5.3 from \cite{GineNickl2016} we end up with
\begin{align*}
J(\FFF,C_{\FFF},\delta)\leq\delta~ J(\FFF,C_{\FFF},1)\leq 2~\delta\sqrt{\ln(2~e~[D + 2]) + (D + 1)~\ln(16\cdot e)}
\end{align*}
for every positive envelope $C_{\FFF}$ of $\FFF$ and any $\delta\in ]0,1[$.
\end{example}
\par 
Now, it is already known that the class $\FFF_{\delta,i} := \{\eins_{B_{i}(\theta)}\mid\theta\in\cU_{\delta}\}$ is a VC-subgraph class with $V(\FFF_{\delta,i})\leq (m + 2) s_{i} + 1$ for $i\in\{1,\ldots,r\}, \delta > 0$ (see \cite[Lemma 5.4]{Kraetschmer2023a}). Then by permance properties of VC-subgraph classs, namely \cite[Lemma 9.9, (v) along with (vi)]{Kosorok2008}, we obtain $\FFF^{\Theta,i}_{\delta}$ as a VC-subgraph class with index with index 
$$
V(\FFF^{\Theta,i}_{\delta}) \leq 2 V(\FFF_{\delta,i}) - 1\leq 2 (m + 2) s_{i} + 1
$$
for $i\in\{1,\ldots,r\}$ and $\delta > 0$. Since $C_{\FFF^{\Theta,i}_{\delta}} := |G^{i}(\theta^{*},\cdot)|~\eins_{\widehat{B}_{i\delta}} + (\delta\wedge 1)^{2}$ defines a positive envelope of $\FFF^{\Theta,i}_{\delta}$ we end up with
\begin{align}
\sup_{\MQ\in \cM_{\textrm{\tiny fin}}}N\big(\eta~\|C_{\FFF^{\Theta,i}_{\delta}}\|_{\MQ,2},\FFF^{\Theta,i}_{\delta},L^{2}(\MQ)\big)
&\label{Abschaetzung2}\leq e~[2 (m + 2) s_{i} + 1]~\big(4~e^{1/2}/\eta\big)^{4 (m + 2) s_{i}}
\end{align}
for $i\in\{1,\ldots,r\}$, $\delta > 0$ and $\eta\in ]0,1[$. Next, fix $\MQ\in\cM_{\textrm{\tiny fin}}$, $\delta, \eta > 0$. Let $h^{0}, \overline{h}^{0}\in \overline{\FFF}^{\Theta}_{\delta}$ and $h^{i},\overline{h}^{i}\in\FFF^{\Theta,i}_{\delta}$ such that $\|h^{i}-\overline{h}^{i}\|_{\MQ,2}\leq \eta \|C_{\FFF^{\Theta,i}_{\delta}}\|_{\MQ,2}/(2r + 2)$ for $i = 1,\ldots,r$, and $\|h^{0}-\overline{h}^{0}\|_{\MQ,2}\leq \eta \|C_{\overline{\FFF}^{\Theta}_{\delta}}\|_{\MQ,2}/(2r + 2)$. Then 
by $\sqrt{\sum_{i=0}^{r}t_{i}}\geq\sum_{i=0}^{r}\sqrt{t_{i}}/(r +1)$ for $t_{0},\ldots,t_{r}\geq 0$
\begin{align*}
\|\sum_{i=0}^{r} h^{i} - \sum_{i=0}^{r}\overline{h}^{i}\|_{\MQ,2}
\leq 
\sum_{i=0}^{r}\|h^{i}-\overline{h}^{i}\|_{\MQ,2} 
&\leq 
\frac{\eta}{2(r + 1)} \big[\|C_{\overline{\FFF}^{\Theta}_{\delta}}\|_{\MQ,2} + \sum_{i=1}^{r}\|C_{\FFF^{\Theta,i}_{\delta}}\|_{\MQ,2}\big]\\
&\leq 
\frac{\eta}{2} \|C_{\overline{\FFF}^{\Theta}_{\delta}} + \sum_{i=1}^{r}C_{\FFF^{\Theta,i}_{\delta}}\|_{\MQ,2} 
\leq \eta \|\xi^{G}_{\delta}\|_{\MQ,2}.
\end{align*}
Thus in view of representation \eqref{Eingangsrepraesentation}
\begin{align}
&\nonumber 
N\big(\eta \|\xi^{G}_{\delta}\|_{\MQ,2},\FFF^{\Theta}_{\delta},L^{2}(\MQ)\big)\\
\label{random entropies Summen}
&\leq 
N\big(\eta \|C_{\overline{\FFF}^{\Theta}_{\delta}}\|_{\MQ,2}/(2r + 2),\overline{\FFF}^{\Theta}_{\delta},L^{2}(\MQ)\big)\cdot\prod_{i=1}^{r}N\big(\eta \|C_{\FFF^{\Theta,i}_{\delta}}\|_{\MQ,2}/(2r + 2),\FFF^{\Theta,i}_{\delta},L^{2}(\MQ)\big)
\end{align}
for $\MQ\in\cM_{\textrm{\tiny fin}}$ and $\delta, \eta > 0$. 
\medskip

Combining \eqref{Abschaetzung1} and \eqref{Abschaetzung2} with \eqref{random entropies Summen}, we obtain for $\delta > 0$ and $\varepsilon\in ]0,1]$ by change of variable formula
\begin{align*}
J(\FFF^{\Theta}_{\delta},\xi^{G}_{\delta},\varepsilon) 
&= 
\varepsilon \int_{0}^{1}\hspace*{-0.1cm}\sup_{\MQ\in\cM_{\textrm{\tiny fin}}}\sqrt{\ln\Big(2 N\big(\varepsilon \eta\|\xi^{G}_{\delta}\|_{\MQ,2},\FFF^{\Theta}_{\delta},L^{2}(\MQ)\big)\Big)}~d\eta\\
&\leq \varepsilon \sum_{i=1}^{2}\int_{0}^{1}\sqrt{v_{i} \ln(K_{i,\varepsilon}/\eta)}~d\eta,
\end{align*}
where
$$
v_{1} := 2 (m + 2)\sum_{i=1}^{r}s_{i} + m \sum_{i=1}^{r}1/\beta_{i}\quad\mbox{and}\quad v_{2} := 4 (m + 2)\sum_{i=1}^{r}s_{i},
$$
and
\begin{align*} 
&
K_{1,\varepsilon} := 
\frac{8~ (r + 1)~\big\{2\cdot 9^{r\cdot m}16^{(m + 2)~\sum_{i=1}^{r} s_{i}}~e^{r + (m + 2)~\sum_{i=1}^{r} s_{i}} ~\prod_{i=1}^{r}([m + 2]~s_{i} + 1)\big\}^{1/v_{1}}}{\varepsilon},\\
&
K_{2,\varepsilon} := \frac{8~(r + 1)~\big\{e^{r + 2 (m + 2)\sum_{i=1}^{r}s_{i}}~\prod_{i=1}^{r}(2[m + 2]~s_{i} + 1)\big\}^{1/v_{2}}}{\varepsilon}.
\end{align*}
Now, we may finish the proof of Proposition \ref{important covering numbers} via 
\eqref{Integralabschaetzung} by routine calculations.
\hfill$\Box$

\subsection{Proofs of results from Section \ref{estimates under Average Value at Risk}}
\subsubsection{Proof of Lemma \ref{completing minimizers}}
\label{Proof of Lemma completing minimizers}
Since $G$ is measurable w.r.t. $\cB(\Theta)\otimes\cB(\RRR^{d})$, the process
$$
\overline{X}^{\alpha}_{n}: \RRR\times\Omega\rightarrow\RRR,~(x,\omega)\mapsto \frac{1}{n}\sum_{j=1}^{n}\frac{\big[G\big(\widehat{\theta}^{\alpha}_{n}(\omega),Z_{j}(\omega)\big) + x\big]^{+}}{1 - \alpha} - x
$$
is continuous in $x$ and measurable in $\omega$. Furthermore
\begin{align*}
\lim_{x\to -\infty}\overline{X}^{\alpha}_{n}(x,\omega) = \lim_{x\to \infty}\overline{X}^{\alpha}_{n}(x,\omega) = \infty\quad\mbox{for}~\omega\in\Omega.
\end{align*}
Hence the paths of $\overline{X}^{\alpha}_{n}$ always attain their minimum. Putting all together, and noting that $\RRR$ with the usual topology is a Polish space, we may find by a version of the measurable selection theorem (e.g. \cite[Theorem 6.7.22]{Pfanzagl1994}) a random variable $\widehat{x}^{\alpha}_{n}$ on $\OFP$ whose realizations minimize the paths of $\overline{X}^{\alpha}_{n}$. Then obviously, $(\widehat{\theta}^{\alpha}_{n},\widehat{x}^{\alpha}_{n})$ is a stochastic minimizer of the auxiliary problems \eqref{auxiliary SAA} w.r.t. the sample size $n$. This completes the proof.
\hfill$\Box$
\subsubsection{Proof of Lemma \ref{countable parameter subsets}}
\label{Proof of Lemma countable parameter subsets}
Let $\varepsilon > 0$ and 
$\mathcal{J}_{\varepsilon} := ]0,\varepsilon]\cap(\mathbb{Q}\cup\{\varepsilon\})$. According to \ref{assump:A3alpha} we may find for any $\overline{\varepsilon}\in\mathcal{J}_{\varepsilon}$ some at most countable subset $\cC(\cU^{\alpha}_{\overline{\varepsilon}})\subseteq \cU^{\alpha}_{\overline{\varepsilon}}$ and $(\MP^{Z})^{n}$-null sets $N_{n}^{\overline{\varepsilon}}$ $(n\in\NNN)$
such that $\theta^{\alpha,*}\in\cC(\cU_{\overline{\varepsilon}}^{\alpha})$, and
\begin{equation}
\label{Anwendung A6}
\inf_{\vartheta\in\cC(\cU^{\alpha}_{\overline{\varepsilon}})}\left\{\EEE\big[\big|G(\theta,Z_{1}) - G(\vartheta,Z_{1})\big|\big] + \max_{j\in\{1,\ldots,n\}}\big|G(\theta,z_{j}) - G(\vartheta,z_{j})\big|\right\} = 0
\end{equation}
if $\theta\in\cU^{\alpha}_{\overline{\varepsilon}}$, $n\in\NNN$ and $(z_{1},\ldots,z_{n})\in\RRR^{d n}\setminus N_{n}^{\overline{\varepsilon}}$.
Now set
\begin{align*}
\cK_{\varepsilon} := \cV_{\RRR,\varepsilon}~\cap~
\bigcup_{\overline{\varepsilon}\in\mathcal{J}_{\varepsilon}}\cC(\cU^{\alpha}_{\overline{\varepsilon}})\times[\mathcal{I}\cap(\mathbb{Q}\cup\{x^{\alpha,*}\})]\quad\mbox{and}\quad \widehat{N}_{n}^{\varepsilon} := \bigcup_{\overline{\varepsilon}\in\mathcal{J}_{\varepsilon}} N_{n}^{\overline{\varepsilon}}~(n\in\NNN).
\end{align*}
Note that $\cK_{\varepsilon}$ is an at most countable subset of $\cV_{\mathcal{I},\varepsilon}$, whereas $\widehat{N}_{n}^{\varepsilon}$ is a $(\MP^{Z})^{n}$-null set for 
$n\in\NNN$. It remains to show for $(\theta,x)\in\cV_{\mathcal{I},\varepsilon}$, $n\in\NNN$ and $(z_{1},\ldots,z_{n})\in\RRR^{d N}\setminus \widehat{N}_{n}^{\varepsilon}$
\begin{align*}
&\inf_{(\vartheta,y)\in\cK_{\varepsilon}}\hspace*{-0.1cm}\big\{\EEE\big[\big|\widehat{G}_{\alpha}\big((\theta,x),Z_{1}) - \widehat{G}_{\alpha}\big((\vartheta,y),Z_{1}\big)\big|\big] +\max_{j\in\{1,\ldots,n\}}\big|\widehat{G}_{\alpha}\big((\theta,x),z_{j}) - \widehat{G}_{\alpha}\big((\vartheta,y),z_{j}\big)\big|\big\}\\ 
&= 0.
\end{align*}
For this purpose fix $(\overline{\theta},\overline{x})\in\cV_{\mathcal{I},\varepsilon}$, $n\in\NNN$ and $(\overline{z}_{1},\ldots,\overline{z}_{n})\in\RRR^{d n}\setminus \widehat{N}_{n}^{\varepsilon}$. If $\overline{x} = x^{\alpha,*}$, note that $\overline{\theta}\in\cU^{\alpha}_{\varepsilon}$ and $(\overline{z}_{1}\ldots,z_{n})\in\RRR^{d}\setminus N_{n}^{\varepsilon}$ so that we may choose by \ref{assump:A3alpha} a sequence $(\vartheta_{l})_{l\in\NNN}$ in 
$\cC(U^{\alpha}_{\varepsilon})$ fulfilling
\begin{align*}
\EEE[|G(\overline{\theta},Z_{1}) - G(\vartheta_{l},Z_{1})|] 
+ 
\max_{j\in\{1,\ldots,n\}}\big|G(\overline{\theta},\overline{z}_{j}) - G(\vartheta_{l},\overline{z}_{j})\big|\to 0\quad\mbox{for}~l\to\infty.
\end{align*}
Then $\big((\vartheta_{l},x^{\alpha,*})\big)_{l\in\NNN}$ is a sequence in $\cK_{\varepsilon}$ which in view of \eqref{basic inequalities} satisfies
\begin{align*}
&
\lim_{l\to\infty}\EEE\big[\big|\widehat{G}_{\alpha}\big((\overline{\theta},\overline{x}),Z_{1}) - \widehat{G}_{\alpha}\big((\vartheta_{l},x^{\alpha,*}),Z_{1}\big)\big|\big]\\ 
&= 
\lim_{l\to\infty}\max_{j\in\{1,\ldots,n\}}\big|\widehat{G}_{\alpha}\big((\overline{\theta},\overline{x}),\overline{z}_{j}) - \widehat{G}_{\alpha}\big((\vartheta_{l},x^{\alpha,*}),\overline{z}_{j}\big)\big|= 0.
\end{align*}
So let $\overline{x}\not=x^{\alpha,*}$. Firstly, we may choose for any $l\in\NNN$ some $\widehat{x}_{l}\in \mathcal{I}\cap \mathbb{Q}$ satisfying
\begin{align*}
|\overline{x} - \widehat{x}_{l}| < |\overline{x} - x^{\alpha,*}|/l\quad\mbox{and}\quad 
|\widehat{x}_{l} - x^{\alpha,*}| < |\overline{x} - x^{\alpha,*}|.
\end{align*}
In particular there is some $\overline{\varepsilon}_{l}\in\mathcal{J}_{\varepsilon}\setminus\{\varepsilon\}$ such that $\|\overline{\theta} - \theta^{\alpha,*}\|_{m} < \overline{\varepsilon}_{l} < \sqrt{\varepsilon^{2} - |\widehat{x}_{l} - x^{\alpha,*}|^{2}}$. Since $(\overline{z}_{1},\ldots,\overline{z}_{n})\in\RRR^{d}\setminus N_{n}^{\overline{\varepsilon}_{l}}$, we may find by \ref{assump:A3alpha} a some $\widehat{\vartheta}_{l}\in \cC(U^{\alpha}_{\overline{\varepsilon}_{l}})$ such that 
\begin{align*}
\EEE[|G(\overline{\theta},Z_{1}) - G(\widehat{\vartheta}_{l},Z_{1})|] + \max_{j\in\{1,\ldots,n\}}\big|G(\overline{\theta},\overline{z}_{j}) - G(\overline{\vartheta}_{l},\overline{z}_{j})\big| < 1/l.
\end{align*}
By choice $(\widehat{\vartheta}_{l},\widehat{x}_{l})$ belongs to $\cK_{\varepsilon}$ for $l\in\NNN$, and in view of \eqref{basic inequalities} we end up with
\begin{align*}
&\EEE\big[\big|\widehat{G}_{\alpha}\big((\overline{\theta},\overline{x}),Z_{1}) - \widehat{G}_{\alpha}\big((\widehat{\vartheta}_{l},\widehat{x}_{l}),Z_{1}\big)\big|\big] +\max_{j\in\{1,\ldots,n\}}\big|\widehat{G}_{\alpha}\big((\overline{\theta},\overline{x}),\overline{z}_{j}) - \widehat{G}_{\alpha}\big((\widehat{\vartheta}_{l},\widehat{x}),\overline{z}_{j}\big)\big|\\
&
\leq \frac{1 + 2~(2-\alpha)~|\overline{x} - x^{\alpha,*}|}{l (1- \alpha)}\to 0\quad\mbox{for}~l\to\infty.
\end{align*}
This completes the proof.
\hfill$\Box$
\subsubsection{Proof of Lemma \ref{relationship2}}
\label{Beweis relationship2}
In view of \eqref{basic inequalities} the Borel measurable mapping $C_{\FFF^{\Theta}_{\alpha,\mathcal{I},\delta}}$ may be verified easily as a positive envelope of $\FFF^{\Theta}_{\alpha,\mathcal{I},\delta}$. By construction it is square $\MP^{Z}$-integrable because $\xi_{\delta}$ fulfills this property.
\smallskip

Next, let $\MQ\in\cM_{\textrm{\tiny fin}}, \eta > 0$ and $\theta, \vartheta\in\cU_{\delta}$ as well as $x,y\in \mathcal{I}_{\delta} := \mathcal{I}~\cap~[x^{\alpha,*} - \delta,x^{\alpha,*} +\delta]$ with 
$$
\|G(\theta,\cdot) - G(\vartheta,\cdot)\|_{\MQ,2}\leq \eta~\|\xi_{\delta}\|_{\MQ,2}/4\quad\mbox{and}\quad |x - y|\leq \eta~\delta/4.
$$
Then by \eqref{basic inequalities} again, we may observe
\begin{align*}
& 
\big\|\big[\widehat{G}_{\alpha}\big((\theta,x)\cdot\big) - \widehat{G}_{\alpha}\big((\theta^{\alpha,*},x^{\alpha,*})\cdot\big)\big] - \big[\widehat{G}_{\alpha}\big((\vartheta,y)\cdot\big) - \widehat{G}_{\alpha}\big((\theta^{\alpha,*},x^{\alpha,*})\cdot\big)\big]\big\|_{\MQ,2}^{2}\\
& 
\leq 
4~\big[\|G(\theta,\cdot) - G(\vartheta,\cdot)\|_{\MQ,2}^{2} + (2-\alpha)^{2}~|x-y|^{2}\big]/(1-\alpha)^{2}\\
&
\leq 
\eta^{2}~\big[\|\xi_{\delta}\|_{\MQ,2}^{2} + (2-\alpha)^{2}~\delta^{2}\big]/[4~(1-\alpha)^{2}]
\leq 
\eta^{2}~\|C_{\FFF^{\Theta}_{\alpha,\mathcal{I},\delta}}\|_{\MQ,2}^{2}/4.
\end{align*}
Hence, denoting by $N\big(\eta~\delta/2,\mathcal{I}_{\delta},|\cdot|\big)$ the minimal number to cover the set $\mathcal{I}_{\delta}$ by intervals $[x - \eta~\delta/2, x + \eta~\delta/2]$ with $x\in\mathcal{I}_{\delta}$, we may conclude from the observation that $\FFF^{\Theta}_{\alpha,\mathcal{I},\delta}$ is a subset of $\big\{\widehat{G}_{\alpha}\big((\theta,x),\cdot\big) - \widehat{G}_{\alpha}\big((\theta^{\alpha,*},x^{\alpha,*}),\cdot\big)\mid (\theta,x)\in\cU_{\delta}\times\mathcal{I}_{\delta}\big\}$
\begin{align*}
N\big(\eta~\|C_{\FFF^{\Theta}_{\alpha,\mathcal{I},\delta}}\|_{\MQ,2},\FFF^{\Theta}_{\alpha,\mathcal{I},\delta},L^{2}(\MQ)\big)
&\leq 
N\big(\eta~\|\xi_{\delta}\|_{\MQ,2}/4,\FFF^{\Theta}_{\delta},L^{2}(\MQ)\big)\cdot 
N\big(\eta~\delta/4,\mathcal{I}_{\delta},|\cdot|\big)\\
&\leq 
N\big(\eta~\|\xi_{\delta}\|_{\MQ,2}/4,\FFF^{\Theta}_{\delta},L^{2}(\MQ)\big)\cdot 4~\big(\sup \mathcal{I}_{\delta} - \inf \mathcal{I}_{\delta}\big)/(\eta~\delta)\\
&\leq 
N\big(\eta~\|\xi_{\delta}\|_{\MQ,2}/4,\FFF^{\Theta}_{\delta},L^{2}(\MQ)\big)\cdot 8/\eta.
\end{align*} 
Now, invoking subadditivity of $\sqrt{\cdot}$ together with the change of variable formula, we end up with 
\begin{align*}
J(\FFF^{\Theta}_{\alpha,\mathcal{I},\delta},C_{\FFF^{\Theta}_{\alpha,\mathcal{I},\delta}},\varepsilon)
\leq 
4~J(\FFF^{\Theta}_{\delta},\xi_{\delta},\varepsilon/4) + \varepsilon~\int_{0}^{1}\sqrt{\ln\big(8/[\varepsilon \eta]\big)}~d\eta\quad\mbox{for}~\varepsilon\in ]0,1].
\end{align*}
The proof may be finished immediately by applying inequality \eqref{Integralabschaetzung}.
\hfill$\Box$
\section*{Acknowledgments}
The author thanks two anonymous referees and the associated editor for useful comments and suggestions which have helped to improve an earlier draft.

\section*{Declarations}

\textbf{Conflicts of interest:} The author has no competing interests to declare that are relevant to the content of this article.

\medskip\noindent
\textbf{Data Availability Statement:} No datasets were generated or analysed during the study

\medskip\noindent
\textbf{Funding:} The authors did not receive support from any organization for the submitted work.



\begin{thebibliography}{0}
\bibitem{Bauer2001} Bauer, H. (2001) \textsl{Measure and integration theory}. de Gruyter, Berlin.
\bibitem{BonnansShapiro2000} Bonnans, J. F. and Shapiro, A. (2000). \textit{Pertubation analysis of optimization problems}. Springer, New York.
\bibitem{BrownPurves1973} Brown, L. D. and Purves, R. (1973). {\em Measurable selections of extrema}. Annals of Statistics \textbf{1}, 902--912.
\bibitem{DentchevaEtAl2017} Dentcheva, D., Penev, S. and Ruszczynski, A. (2017). \textit{Statistical estimation of composite risk functionals and risk optimization problems}. Annals of the Institute of Statistical Mathematics \textbf{69}, 737--760.
\bibitem{EichhornRoemisch2007} Eichhorn, A. and R\"omisch, W. (2007) \textit{Stochastic integer programming: Limit theorems and confidence intervals}. Mathematics of Operations Research 32, 118--135.
\bibitem{ErmolievNorkin2013} Ermoliev, Y. M. and Norkin, V. I. (2013). \textit{Sample Average Approximation method for compound stochastic optimization problems}. SIAM J. OPTIM. \textbf{23}, 2231--2263.
\bibitem{GineNickl2016} Gine, E. and Nickl, R.  (2016). \textit{Mathematical Foundations of Infinite-Dimensional Statistical Models}. Cambridge University Press, Cambridge.
\bibitem{GoetzeEtAl2021} Götze, F., Sambale, H. and Sinulis, A. (2021). \textit{Concentration inequalities for polynomials in $\alpha$-sub-exponential random variables}. {\em Electron. J. Probab.} \textbf{26}, 1--22.
\bibitem{Guigues2017} Guigues, V. (2017). \textit{Multistep stochastic mirror descent for risk-averse convex stochastic programs based on extended polyhedral risk measures}. Mathematical Programming \textbf{163}, 169--212.
\bibitem{GuiguesEtAl2017} 
Guigues, V., Juditsky, A. and Nemirovski, A. (2017). \textit{Non-asymptotic confidence bounds for optimal value of a stochastic program}. 
Optimization Methods and Software 32, 1033--1058.
\bibitem{GuiguesKraetschmerShapiro2018} Guigues, V., Kr\"atschmer, V. and Shapiro, A. (2018). \textit{A central limit theorem and hypotheses testing for risk-averse stochastic programs}. SIAM J. OPTIM. \textbf{28}, 1337--1366.
\bibitem{KainaRueschendorf2007} Kaina, M. and R\"uschendorf, L. (2009). \textit{On convex risk measures on $L^{p}-$spaces}. Math. Methods Oper. Res. \textbf{69}, 475 -- 495. 
\bibitem{Kosorok2008} Kosorok, M. R. (2008). \textsl{Introduction to empirical processes and semiparametric inference}. Springer, New York.
\bibitem{Kraetschmer2023a} Kr\"atschmer, V. (2024a). \textit{Nonasymptotic upper estimates for errors of the sample average approximation method to solve risk averse stochastic programs}. SIAM Journal on Optimization \textbf{34}, 1264--1294.
\bibitem{Kraetschmer2020} Kr\"atschmer, V. (2024b). \textit{First order asymptotics of the sample average approximation method to solve risk averse stochastic programs}. Mathematical Programming \textbf{208}, 209--242.
%
\bibitem{KuchibhotlaChakrabortty2022} Kuchibhotla, A. K. and Chakrabortty, A. (2022). \textit{Moving beyond sub-Gaussianity in high-dimensional statistics: applications in covariance estimation and linear regression}. {\em Information and Inference: A Journal of the IMA} \textbf{11}, 1389--1456.
\bibitem{OstrovskiiBach2021} Ostrovskii, D. M. and Bach, F. (2021). \textit{Finite-sample analysis of M-estimators using self-concordance}. {\em Electronic Journal of Statistics} \textbf{15}, 326--391.
\bibitem{Pfanzagl1994} Pfanzagl, J. (1994). \textsl{Parametric Statistical Theory}, de Gruyter, Berlin, New York.
\bibitem{Pflug1999} Pflug, G. Ch. (1999). Stochastic programs and statistical data. {\em Ann. Oper. Res.} 85, 59--78.
\bibitem{ShapiroEtAl} Shapiro, A., Dentcheva, D. and Ruszczynski, A. (2014). \textsl{Lectures on stochastic programming.} MOS-SIAM Ser. Optim., Philadelphia (2nd ed.).
\bibitem{vanderVaart1998} van der Vaart, A.W. (1998). \textsl{Asymptotic statistics}. Cambridge University Press, Cambridge.
\bibitem{vanderVaartWellner1996} van der Vaart, A.W. and Wellner, J.A. (1996). \textsl{Weak convergence and empirical processes}. Springer, New York.
\end{thebibliography}
\end{document}